\documentclass[12pt]{amsart}

\usepackage[T1]{fontenc}
\usepackage{amssymb}
\usepackage{bm}
\usepackage{ytableau}
\usepackage{amsmath}
\usepackage{quiver}
\usepackage{amssymb,amsfonts,amsthm,amsmath,calligra}
\usepackage{slashed}
\usepackage{yfonts}
\usepackage{mathrsfs,pifont}
\usepackage{float}
\usepackage{comment}

\usepackage[all]{xy}

\usepackage{tikz}

\usepackage{graphicx}
\usepackage{xcolor}

\usepackage{amssymb,amsfonts,amsthm,amsmath}
\usepackage[all]{xy}
\usepackage{slashed}
\usepackage{yfonts}
\usepackage{mathrsfs,pifont}
 \usepackage{mathpple}
\usepackage{boondox-cal}
\usepackage{amsmath,amsfonts, amscd,amsthm, graphicx}
  \usepackage{amsrefs}
\usepackage[all]{xy}
\usepackage{color}
\usepackage{tikz-cd}

\usepackage{hyperref}
\def\vsk#1>{\vskip#1\baselineskip}

\def\XXint#1#2#3{{\setbox0=\hbox{$#1{#2#3}{\int}$}
\vcenter{\hbox{$#2#3$}}\kern-.5\wd0}}

\numberwithin{equation}{section}

\newcommand{\mbf}{\mathbf}
\newcommand{\mbb}{\mathbb}
\newcommand{\mf}{\mathfrak}
\newcommand{\mc}{\mathcal}

\theoremstyle{plain}
\newtheorem{thm}{Theorem}[section]

\newtheorem{thm-defn}{Theorem/Definition}[section]
\newtheorem{lem}[thm]{Lemma}
\newtheorem{lem-defn}[thm]{Lemma/Definition}
\newtheorem{prop}[thm]{Proposition}

\newtheorem{prop-defn}[thm]{Proposition-Definition}

\newtheorem{thm-alg}[thm]{Theorem/Algorithm}

\allowdisplaybreaks[4]  

\begin{document}

 \title{From quantum difference equation to Dubrovin connection of affine type $A$ quiver varieties}
  \author{Tianqing Zhu}
  \address{Department of Mathematics, Columbia University and Yau Mathematical Sciences Center, Tsinghua University}
\email{tz2611@columbia.edu}

  \date{}

  \maketitle
\begin{abstract} This is the continuation of the article \cite{Z23}. In this article we analyze in detail the quantum difference equation of the equivariant $K$-theory of the affine type $A$ quiver varieties. We obtain a convenient presentation of the algebraic quantum difference operator $\mbf{B}_{\mc{L}}^s(z)$. We also give the detailed analysis of the connection matrix for the quantum difference equation in the nodal limit $p\rightarrow0$. Using these two results, we prove that the degeneration limit of the quantum difference equation is the Dubrovin connection for the quantum cohomology of the affine type $A$ quiver varieties, and the monodromy representation for the Dubrovin connection is generated by the monodromy operators $\mbf{B}_{\mbf{m}}$.
\end{abstract}

\tableofcontents
\section{\textbf{Introduction}}
This paper is a continuation of the paper \cite{Z23}. We do the explicit calculation for the quantum difference equation for the affine type $A$ quiver varieties $M(\mbf{v},\mbf{w})$ with stability condition $\bm{\theta}=(1,1,\cdots,1)$.

The quantum difference equation is the difference-equation analog of the quantum differential equation in \cite{MO12}. It is the difference equation solved by the capping operator $\mbf{J}(u,z)\in K_{T}(X\times X)_{loc}[[q^{\mc{L}}]]_{\mc{L}\in\text{Pic}(X)}$ described in \cite{OS22} and \cite{O15}. The capping operator is the operator that gives the $K$-theoretic counting of the stable quasi-map from $\mbb{P}^1$ to the fixed Nakajima quiver varieties $M(\mbf{v},\mbf{w})$. The difference equation can be written as:
\begin{align*}
\Psi(q^{\mc{L}}z)=\mbf{M}_{\mc{L}}(z)\Psi(z)\in K_{T}(X)_{loc}\otimes\mbb{C}((\text{Pic}(X))).
\end{align*}

It has been shown in section $7$ in \cite{O15} that the fundamental solution normalised at $z=0$ for the quantum difference equation is given by the capping operator $\mbf{J}(u,z)$. The solution $\Psi(z)$ around different regions also gives rise to the analytic continuation of $\mbf{J}(u,z)$ over the Picard torus $\text{Pic}(X)\otimes\mbb{C}^{\times}$, which relates the capping operator of quiver varieties of different stability conditions. 

To solve the quantum difference equation, we need to first figure out the expression of the quantum difference operator $\mbf{M}_{\mc{L}}(z)$. The first systematic construction is given by Okounkov and Smirnov \cite{OS22}. They use the quantum algebra defined by the RTT formalism of the $K$-theoretic stable envelope to express the quantum difference operator $\mbf{M}_{\mc{L}}(z)$, and the formula can be written as the product of the monodromy operators $\mbf{B}_{w}(z)$:
\begin{align*}
\mbf{M}_{\mc{L}}(z)=\mc{L}\prod_{w\in Walls}\mbf{B}_{w}(z).
\end{align*}

Each monodromy operator $\mbf{B}_{w}(z)\in\widehat{U_q(\mf{g}_{w})}$ lies in the completion of the slope subalgebra $\widehat{U_q(\mf{g}_{w})}$ and can be solved via the ABRR equation in \cite{OS22} and \cite{ABRR97}. 

In \cite{Z25}, it has been shown that the geometric quantum loop group $U_{q}(\hat{\mf{g}}_{Q})$ constructed via the $K$-theoretic stable envelope is isomorphic to the quantum affine algebra of the corresponding quiver type $U_{q}(\hat{\mf{g}}_{Q})$. 

On the algebraic side, one can construct the algebraic analog of the quantum difference equation on the quantum affine algebra $U_{q}(\hat{\mf{g}}_{Q})$ \cite{Z23}. For the affine type $A$ case, the construction is based on the slope factorization of the quantum toroidal algebra $U_{q,t}(\hat{\hat{\mf{sl}}}_{n})$ given by \cite{N15}. The formula for the slope subalgebra gives the explicit form of the quantum difference operator and the monodromy operator. It was also shown in \cite{Z23}\cite{Z26} that the algebraic quantum difference equation coincides with the geometric Okounkov-Smirnov quantum difference equation for the Nakajima quiver varieties.

On the other hand, one interesting application of the main result of the paper is that one can compare the computation of the monodromy representation with respect to the computation given by the Okounkov-Smirnov quantum difference equation without using the isomorphism given in \cite{Z25}\cite{Z26}. It can be shown in \cite{Z24-2} that this will also lead to the isomorphism of the geometric quantum loop group of the affine type $A$ and the quantum toroidal algebras.

\subsection{Quantum difference equation for $M(\mbf{v},\mbf{w})$}
In this paper we study the algebraic quantum difference equation of the equivariant $K$-theory of the affine type $A$ quiver varieties $M(\mbf{v},\mbf{w})$ from the algebraic construction given in \cite{Z23}. The first main result is that the algebraic quantum difference operator $\mbf{B}_{\mc{L}}^s(z)$ can have a good representation:
\begin{thm}{(See Theorem \ref{good-representation})}
Let $[s-\mc{L},s)$ be a generic path, with $s$ being generic. Then the quantum difference operator factorises as:
\begin{align*}
\mbf{B}_{\mc{L}}^s(z)=\mc{L}\prod^{\leftarrow}_{\mbf{m}\in[s,s-\mc{L})}\mbf{B}_{\mbf{m}}(z)
\end{align*}
and each $\mbf{B}_{\mbf{m}}$ can be written in one of the following two forms:
\begin{itemize}
	\item $U_{q}(\mf{sl}_2)$ type:
	\begin{equation*}
	\begin{aligned}
	&\prod_{k=0}^{\substack{\rightarrow\\\infty}}\exp_{q^{2}}(-(q-q^{-1})z^{-k(\mbf{v}_{\gamma})}p^{-k\mbf{m}\cdot(\mbf{v}_{\gamma})}q^{k(\mbf{v}_{\gamma})^T((\frac{n^2r-1}{2})\bm{\theta}+\mbf{e}_1)-2-2\delta_{1\gamma}}f_{\gamma}e_{\gamma}'))\\
	=&\sum_{n=0}^{\infty}\frac{(q-q^{-1})^n}{[n]_{q^2}!}\frac{(-1)^n}{\prod_{\nu=1}^{n}(1-z^{-\mbf{v}_{\gamma}}p^{-\mbf{m}\cdot\mbf{v}_{\gamma}}q^{\nu(\mbf{v}_{\gamma})^T((\frac{n^2r-1}{2})\bm{\theta}+\mbf{e}_1)-2-2\delta_{1\gamma}})}f_{\gamma}^{n}e_{\gamma}'^n.
	\end{aligned}
	\end{equation*}
	\item $U_{q}(\hat{\mf{gl}}_1)$ type.
	\begin{equation*}
	\begin{aligned}
	\mu(\prod_{h=1}^{g}(\exp(-\sum_{k=1}^{\infty}\frac{n_kq^{-\frac{k\lvert\bm{\delta_h}\lvert}{2}}}{1-z^{-k\lvert\bm{\delta_h}\lvert}p^{k\mbf{m}\cdot\bm{\delta}_{h}}q^{-\frac{k\lvert\bm{\delta_h}\lvert}{2}}}\alpha^{\mbf{m},h}_{-k}\otimes\alpha^{\mbf{m},h}_{k}).
	\end{aligned}
	\end{equation*}
	with $\mu:\mc{B}_{\mbf{m}}\otimes\mc{B}_{\mbf{m}}\rightarrow\mc{B}_{\mbf{m}}$ the multiplication map.
\end{itemize}
\end{thm}

Roughly speaking, the theorem implies that for the generic slope points $\mbf{m}\in\mbb{Q}^n$, the corresponding monodromy operator $\mbf{B}_{\mbf{m}}(z)$ can be totally written in either $U_{q}(\mf{sl}_2)$ type or $U_{q}(\hat{\mf{gl}}_1)$ type. This explains the philosophy that the slope subalgebra $\mc{B}_{\mbf{m}}$ can be thought of as being generated by the wall subalgebra $U_{q}(\mf{g}_w)$ such that the wall $w\subset\mbb{R}^n$ contains the slope point $\mbf{m}$.

\subsection{Connection matrix of the quantum difference equation}
One important aspect of the quantum difference equation is the connection matrix between the solution expanding from $z=0^{\bm{\theta}_1}$ to $z=0^{\bm{\theta}_2}$. It is defined as:
\begin{align*}
\textbf{Mon}_{\bm{\theta}_1,\bm{\theta}_2}(z):=\Psi_{\bm{\theta}_1}(z)^{-1}\Psi_{\bm{\theta}_2}(z).
\end{align*}

Here $\bm{\theta}_i\in(\mbb{Z}-\{0\})^n$ such that $0^{\theta}=0$ if $\theta>0$, $0^{\theta}=\infty$ if $\theta<0$. $\Psi_{\bm{\theta}}(z)$ is the fundamental solution around $z=0^{\bm{\theta}}$.

In the settings of the Fuchsian $q$-difference equation, the connection matrix plays the role of the monodromy operator in the settings of the differential equations with regular singularities.

For simplicity, in this paper we only give the explicit calculation of the connection matrix for $\bm{\theta}_1=\bm{\theta}=(1,1,\cdots,1)$, $\bm{\theta}_2=-\bm{\theta}=(-1,-1,\cdots,-1)$. Part of the reason is that this connection matrix $\textbf{Mon}_{\bm{\theta}_1,\bm{\theta}_2}(z)$ can be thought of as the "largest" connection matrix with the factorization:
\begin{align*}
\textbf{Mon}_{\bm{\theta},-\bm{\theta}}(z)=\textbf{Mon}_{\bm{\theta},\bm{\theta}_1}(z)\textbf{Mon}_{\bm{\theta}_1,\bm{\theta}_2}(z)\cdots\textbf{Mon}_{\bm{\theta}_n,-\bm{\theta}}(z)
\end{align*}

If we know the fundamental solutions, the connection matrix can be computed out via the terms in the fundamental solutions. The main result for the computation is the $p\rightarrow0$ limit of the regular part of the connection matrix $\textbf{Mon}^{reg}(zp^s)$:

\begin{thm}{(Theorem \ref{p0-limit-connection})}
For generic $s\in\mbb{R}^n$ with either $s_i\geq0$ for all $i$ or $s_i<0$ for all $i$, the connection matrix has asymptotics as $p\rightarrow0$:
\begin{align*}
\lim_{p\rightarrow0}\textbf{Mon}^{reg}(p^sz)=
\begin{cases}
\prod^{\leftarrow}_{0\leq\mbf{m}<s}(\mbf{B}_{\mbf{m}}^*)^{-1}\cdot\mbf{T},\qquad s\geq0\\
\prod_{s<\mbf{m}<0}^{\leftarrow}\mbf{B}_{\mbf{m}}^*\cdot\mbf{T},\qquad s<0.
\end{cases}
\end{align*}

Here $0\leq\mbf{m}<s$ means the slope points $\mbf{m}$ in one generic path from $0$ to the point $s$ without intersecting $s$.
\end{thm}

Here $\mbf{B}_{\mbf{m}}=\mu((1\otimes S_{\mbf{m}})(R_{\mbf{m}}^{-})^{-1})$. In the case of finite $A$ type, the operator $\mbf{B}_{\mbf{m}}$ was identified as the generators of the monodromy representation of the trigonometric Casimir connection \cite{GL11}.

\subsection{Degeneration of the quantum difference equation}

One important feature about the quiver variety $M(\mbf{v},\mbf{w})$ is the Dubrovin connection over its quantum cohomology $H_{T}^*(M(\mbf{v},\mbf{w}))[[q^d]]_{d\in H^2(M(\mbf{v},\mbf{w}))_{eff}}$. 

For the Dubrovin connection for the affine type $A$ quiver varieties $M(\mbf{v},\mbf{w})$, it can be written as:
\begin{align*}
\nabla_{\lambda}=d_{\lambda}-Q(\lambda),\qquad\lambda\in H^2(M(\mbf{v},\mbf{w}),\mbb{Z}).
\end{align*}

Here $Q(\lambda)=c_{1}(\lambda)\cup(-)+\cdots$ is the quantum multiplication operator. It has been proved by Maulik-Okounkov \cite{MO12} that the quantum multiplication operator can be written as:
\begin{align*}
Q(\lambda)=c_1(\lambda)+\sum_{\alpha>0}\frac{(\alpha,\lambda)}{1-z^{-\alpha}}e_{\alpha}e_{-\alpha}+\text{constant}.
\end{align*}
where $e_{\pm\alpha}\in\mf{g}^{MO}$ are the generators of the Maulik-Okounkov Lie algebra of the affine type $A$, and constant stands for some scaling operators. It has been proved in \cite{BD23} that the Lie algebra $\mf{g}^{MO}_{Q}$ is isomorphic to the affine Lie algebra $\hat{\mf{gl}}_n$, which means that the quantum multiplication operator can be written as:
\begin{align*}
Q(\lambda)=c_1(\lambda)+\sum_{i<j}\frac{(\lambda\cdot[i,j))}{1-z^{-[i,j)}}E_{[i,j)}E_{-[i,j)}+\text{const}
\end{align*}
where $E_{\pm[i,j)}\in\hat{\mf{gl}}_n$ are the generators of the affine Lie algebra $\hat{\mf{gl}}_n$. Thus, one can identify the Dubrovin connection of the quiver varieties as the trigonometric Casimir connection.

The important feature of the Dubrovin connection is that the solution and the monodromy have played a vital role in the Gromov-Witten theory of the quiver varieties. On the side of the representation theory, the monodromy of the trigonometric Casimir connection is also an important problem in the representation theory of the Yangian, one can refer to \cite{ATL15}\cite{GL11} for the further reference on the computation for the monodromy. It is also conjectured that the monodromy representation of the Dubrovin connection for $M(\mbf{v},\mbf{w})$ should be generated by $\mbf{B}_{w}:=\mbf{m}((1\otimes S_{w})(R_{w}^{-})^{-1})$ with $w\subset\text{Pic}(M(\mbf{v},\mbf{w}))\otimes\mbb{Q}$ the wall hyperplane in the rational Picard group. $R_{w}^-$ is the lower triangular $R$-matrix for the wall subalgebra $U_{q}(\mf{g}_{w})$.

The second main result of the paper is that we can express the monodromy representation of the Dubrovin connection in terms of the monodromy operators:
\begin{thm}{(Theorem \ref{monodromy-rep-casimir})}
The monodromy representation:
\begin{align*}
\pi_{1}(\mbb{P}^n\backslash\textbf{Sing}, 0^+)\rightarrow\text{Aut}(\widehat{H_{T}}(M(\mbf{v},\mbf{w})))
\end{align*}
of the Dubrovin connection is generated by $\mbf{B}_{\mbf{m}}=\mu((1\otimes S_{\mbf{m}})(R_{\mbf{m}}^{-})^{-1})$.
\end{thm}

The proof relies on two facts: First, a suitable degeneration limit of the quantum difference equation for $M(\mbf{v},\mbf{w})$ recovers the Dubrovin connection on $H_{T}^*(M(\mbf{v},\mbf{w}))$. This can be stated as the following theorem:

\begin{thm}{(Theorem \ref{degeneration-qde-thm})}
The degeneration limit of the algebraic quantum difference operator $\mbf{B}_{\mc{L}}^s(z)$ coincides with the quantum multiplication operator $Q(\mc{L})$ up to a constant operator.
\end{thm}

Second, this theorem tells us that one could describe the monodromy representation of the Dubrovin connection in terms of the degeneration limit of the connection matrix in \eqref{trans}:
\begin{align*}
\text{Trans}(s):=\psi_{0}(e^{2\pi is})^{-1}\psi_{\infty}(e^{2\pi is})=\lim_{\tau\rightarrow0}\textbf{Mon}(z=e^{2\pi is},t_1=e^{2\pi i\hbar_1\tau},t_2=e^{2\pi i\hbar_2\tau},q=e^{-2\pi i\tau})
\end{align*}
and we call it the transition matrix. The application of the transition matrix is explained in detail in Section \ref{sub:_p_rightarrow1_limit_of_the_solution_of_the_quantum_difference_equation}. 

It turns out that we could identify $\text{Trans}(s)$ with the $p\rightarrow0$ limit of the regular part of the connection matrix:
\begin{prop}{(Proposition \ref{identi-p-0-connection})}
$\text{Trans}(s)=\lim_{p\rightarrow0}\textbf{Mon}^{reg}(p^s,e^{2\pi i\hbar_1},e^{2\pi i\hbar_2},p)$ for $s\in\mbb{R}^n\backslash\text{Walls}$
\end{prop}

The monodromy representation is generated by $\text{Trans}(s')^{-1}\text{Trans}(s)$ for different generic $s',s$. For $s',s$ close enough, $\text{Trans}(s')^{-1}\text{Trans}(s)$ is equal to $\mbf{B}_{\mbf{m}}$. Thus in this way we have finished the proof of Theorem \ref{monodromy-rep-casimir}.

The proof of the degeneration of the connection matrix to the monodromy representation heavily uses the fact on the modular duality for the Riemann theta function. For completeness we will review the Riemann theta function on the abelian variety in the appendix.

It is worth noting that the construction of the quantum difference equations in \cite{Z23} can be generalized to arbitrary quivers and arbitrary weighted modules of the corresponding shuffle algebras $\mc{S}$. It is given in \cite{Z26}.    For the quiver of finite $ADE$ type, the construction is the analog of the construction given by Etingof and Varchenko in \cite{EV02}, and the degeneration to the corresponding trigonometric Casimir connection is given in \cite{BM15}.

\subsection{Structure of the paper}
The structure of the paper goes in the following way. In section \ref{sec:_textbf_quantum_toroidal_algebra_and_quiver_varieties} and section \ref{sec:_textbf_affine_yangian_and_kac_moody_algebra}, we will introduce some basic ingredients about the quantum toroidal algebra $U_{q,t}(\hat{\hat{\mf{sl}}}_n)$, affine Yangians $Y_{\hbar_1,\hbar_2}(\hat{\mf{sl}}_n)$ and their action on the equivariant $K$-theory $K_T(M(\mbf{w}))$ and on the equivariant cohomology $H_{T}^*(M(\mbf{w}))$ of the affine type $A$ quiver varieties. We will connect them by the degeneration limit in section \ref{sec:_textbf_degeneration_to_the_dubrovin_connection}. Also we will introduce the slope factorization of the quantum toroidal algebra in section \ref{sec:_textbf_quantum_toroidal_algebra_and_quiver_varieties}.

In section \ref{sec:_textbf_review_of_the_quantum_difference_equation} we will review the construction of the quantum difference equation for the affine type $A$ quiver varieties.

In section \ref{sec:_textbf_analysis_of_the_quantum_difference_equations} we will give the detailed analysis of the quantum difference equation, we will show that the quantum difference operator admits good representation (Theorem \ref{good-representation}), which means that one can choose the generic path such that each monodromy operator $\mbf{B}_{\mbf{m}}(z)$ on the path can be written in a self-contained formula that is well-defined for arbitrary $z$ outside of some hyperplanes. In section \ref{sec:_textbf_degeneration_to_the_dubrovin_connection} we use this result to prove that the degeneration limit of the quantum difference equation is equivalent to the Dubrovin connection (Theorem \ref{degeneration-qde-thm}) and the monodromy representation is generated by the monodromy operators (Theorem \ref{monodromy-rep-casimir}).

In section $8$, we will show some computation for the examples of the equivariant Hilbert scheme $\text{Hilb}_{n}([\mbb{C}^2/\mbb{Z}_r])$.

\subsection*{Acknowledgments.}The author would like to thank Andrei Negu\c{t}, Andrei Okounkov, Andrey Smirnov and Hunter Dinkins for their helpful discussions on quantum groups and stable envelopes. The author is supported by the international collaboration grant BMSTC and ACZSP (Grant no. Z221100002722017). Part of this work was done while the author was visiting Department of Mathematics at Columbia University and Mathematics Department of
University of North Carolina at Chapel Hill. The author thanks for their hospitality and provisio of excellent working environment.

\section{\textbf{Quantum toroidal algebra and quiver varieties}}\label{sec:_textbf_quantum_toroidal_algebra_and_quiver_varieties}

\subsection{Quantum toroidal algebra $U_{q,t}(\hat{\hat{\mf{sl}}}_{n})$}
In this subsection we introduce the basic notions for the quantum toroidal algebras. For details one can refer to \cite{N15}.

The quantum toroidal algebra $U_{q,t}(\hat{\hat{\mf{sl}}}_{n})$ is a $\mbb{Q}(q,t)$-algebra defined as:
\begin{align*}
U_{q,t}(\hat{\hat{\mf{sl}}}_{n})=\mbb{Q}(q,t)\langle\{e_{i,d}^{\pm}\}_{1\leq i\leq n}^{d\in\mbb{Z}},\{\varphi_{i,d}^{\pm}\}_{1\leq i\leq n}^{d\in\mbb{N}_{0}}\rangle/(~).
\end{align*}
The relation between the generators can be described in the generating functions:
\begin{equation*}
e_i^{ \pm}(z)=\sum_{d \in \mathbb{Z}} e_{i, d}^{ \pm} z^{-d} \quad \varphi_i^{ \pm}(z)=\sum_{d=0}^{\infty} \varphi_{i, d}^{ \pm} z^{\mp d}.
\end{equation*}

with $\varphi_{i,d}^{\pm}$ commute among themselves and:
\begin{equation*}
\begin{gathered}
e_i^{ \pm}(z) \varphi_j^{ \pm^{\prime}}(w) \cdot \zeta_{ij(ji)}\left(\frac{w^{ \pm 1}}{z^{ \pm 1}}\right)=\varphi_j^{ \pm^{\prime}}(w) e_i^{ \pm}(z) \cdot \zeta_{ji(ij)}\left(\frac{z^{ \pm 1}}{w^{ \pm 1}}\right) \\
e_i^{ \pm}(z) e_j^{ \pm}(w) \cdot \zeta_{ij(ji)}\left(\frac{w^{ \pm 1}}{z^{ \pm 1}}\right)=e_j^{ \pm}(w) e_i^{ \pm}(z) \cdot \zeta_{ji(ij)}\left(\frac{z^{ \pm 1}}{w^{ \pm 1}}\right) \\
{\left[e_i^{+}(z), e_j^{-}(w)\right]=\delta_i^j \delta\left(\frac{z}{w}\right) \cdot \frac{\varphi_i^{+}(z)-\varphi_i^{-}(w)}{q-q^{-1}}}.
\end{gathered}
\end{equation*}

Here $i,j\in\{1,\cdots,n\}$ and $z$,$w$ are variables of color $i$ and $j$:
\begin{equation*}
\zeta_{ij}\left(\frac{x_i}{x_j}\right)=\frac{\left[\frac{x_j}{q t x_i}\right]^{\delta_{j-1}^i}\left[\frac{t x_j}{q x_i}\right]^{\delta_{j+1}^i}}{\left[\frac{x_j}{x_i}\right]^{\delta_j^i}\left[\frac{x_j}{q^2 x_i}\right]^{\delta_j^i}},\qquad[x]=x^{1/2}-x^{-1/2}.
\end{equation*}

with the Serre relation:
\begin{equation*}
\begin{aligned}
& e_i^{ \pm}\left(z_1\right) e_i^{ \pm}\left(z_2\right) e_{i \pm^{\prime} 1}^{ \pm}(w)+\left(q+q^{-1}\right) e_i^{ \pm}\left(z_1\right) e_{i \pm^{\prime} 1}^{ \pm}(w) e_i^{ \pm}\left(z_2\right)+e_{i \pm^{\prime} 1}^{ \pm}(w) e_i^{ \pm}\left(z_1\right) e_i^{ \pm}\left(z_2\right)+\\
& +e_i^{ \pm}\left(z_2\right) e_i^{ \pm}\left(z_1\right) e_{i \pm^{\prime} 1}^{ \pm}(w)+\left(q+q^{-1}\right) e_i^{ \pm}\left(z_2\right) e_{i \pm^{\prime} 1}^{ \pm}(w) e_i^{ \pm}\left(z_1\right)+e_{i \pm^{\prime} 1}^{ \pm}(w) e_i^{ \pm}\left(z_2\right) e_i^{ \pm}\left(z_1\right)=0.
\end{aligned}
\end{equation*}

The standard coproduct structure is imposed as:
\begin{equation*}
\Delta: U_{q,t}(\hat{\hat{\mf{sl}}}_{n})\longrightarrow U_{q,t}(\hat{\hat{\mf{sl}}}_{n})\widehat{\otimes} U_{q,t}(\hat{\hat{\mf{sl}}}_{n}).
\end{equation*}

\begin{equation*}
\begin{array}{ll}
\Delta\left(e_i^{+}(z)\right)=\varphi_i^{+}(z) \otimes e_i^{+}(z)+e_i^{+}(z) \otimes 1 & \Delta\left(\varphi_i^{+}(z)\right)=\varphi_i^{+}(z) \otimes \varphi_i^{+}(z) \\
\Delta\left(e_i^{-}(z)\right)=1 \otimes e_i^{-}(z)+e_i^{-}(z) \otimes \varphi_i^{-}(z) & \Delta\left(\varphi_i^{-}(z)\right)=\varphi_i^{-}(z) \otimes \varphi_i^{-}(z).
\end{array}
\end{equation*}
\subsection{Quantum groups $U_{q}(\hat{\mf{sl}}_{n})$}
The quantum group $U_{q}(\hat{\mf{sl}}_{n})$ is a $\mbb{Q}(q)$-Hopf algebra:
\begin{align*}
U_{q}(\hat{\mf{sl}}_{n})=\mbb{Q}(q)\langle x_{i}^{\pm1},\psi^{\pm1}_{s},c^{\pm1}\rangle^{i\in\mbb{Z}/n\mbb{Z}}_{s\in\{1,\cdots,n\}}
\end{align*}
modulo the fact that $c$ is central, as well as the following relations:
\begin{align*}
&\psi_{s}\psi_{s'}=\psi_{s'}\psi_{s}\\
&\psi_{s}x_{i}^{\pm}=q^{\pm(\delta^{i+1}_{s}-\delta^{i}_{s})}x_{i}^{\pm}\psi_{s}\\
&[x_{i}^{\pm},x_{j}^{\pm}]=0,\qquad\text{if }j\notin\{i-1,i+1\}\\
&[x_{i}^{\pm},[x_{i}^{\pm},x_{j}^{\pm}]_{q}]_{q^{-1}}=0\qquad\text{if }j\in\{i-1,i+1\}\\
&[x_{i}^{+},x_{j}^{-}]=\frac{\delta^{j}_{i}}{q-q^{-1}}(\frac{\psi_{i+1}}{\psi_{i}}-\frac{\psi_{i}}{\psi_{i+1}})
\end{align*}
for all $i,j\in\mbb{Z}/n\mbb{Z}$ and $s,s'\in\{1,\cdots,n\}$. Here $[-,-]_{q}$ is the $q$-bracket:
\begin{align*}
[a,b]_{q}=ab-qba.
\end{align*}

We further require the following relation:
\begin{align*}
\psi_{s+n}=c\psi_{s},\forall s\in\mbb{Z}.
\end{align*}

The counit is given by $\epsilon(x_{i}^{\pm})=0$ and $\epsilon(\psi_{s})=1$, while the coproduct is given by
\begin{align*}
&\Delta(c)=c\otimes c,\Delta(\psi_{s})=\psi_{s}\otimes\psi_{s}\\
&\Delta(x_{i}^{+})=\frac{\psi_{i+1}}{\psi_{i}}\otimes x_{i}^{+}+x_{i}^{+}\otimes1\\
&\Delta(x_{i}^{-})=1\otimes x_{i}^{-}+x_{i}^{-}\otimes\frac{\psi_{i}}{\psi_{i+1}}.
\end{align*}

We also introduce the two half-subalgebras of $U_{q}(\hat{\mf{sl}}_{n})$:
\begin{align*}
&U_{q}^{\geq}(\hat{\mf{sl}}_{n})=\mbb{Q}(q)\langle x_{i}^{+},\psi_{s}^{\pm1},c^{\pm1}\rangle^{i\in\mbb{Z}/n\mbb{Z}}_{s\in\{1,\cdots,n\}}\subset U_{q}(\hat{\mf{sl}}_{n})\\
&U_{q}^{\leq}(\hat{\mf{sl}}_{n})=\mbb{Q}
(q)\langle x_{i}^{-},\psi_{s}^{\pm1},c^{\pm1}\rangle^{i\in\mbb{Z}/n\mbb{Z}}_{s\in\{1,\cdots,n\}}\subset U_{q}(\hat{\mf{sl}}_{n}).
\end{align*}

Both are bialgebras, and they are equipped with the following bialgebra pairing, called the Drinfeld pairing:
\begin{align*}
\langle-,-\rangle:U_{q}^{\geq}(\hat{\mf{sl}}_{n})\otimes U_{q}^{\leq}(\hat{\mf{sl}}_{n})\rightarrow\mbb{Q}(q)
\end{align*}
such that:
\begin{align*}
\langle x_{i}^{+},x_{j}^{-}\rangle=\frac{\delta^{i}_{j}}{q^{-1}-q},\qquad\langle\psi_{s},\psi_{s'}\rangle=q^{-\delta^{s}_{s'}}.
\end{align*}

\subsection{Quantum affine algebra $U_{q}(\hat{\mf{gl}}_{n})$} We use Section 5.1 of \cite{N15} as the detailed reference for the quantum affine algebra $U_{q}(\hat{\mf{gl}}_n)$.

The quantum affine algebra $U_{q}(\hat{\mf{gl}}_{n})$ is a $\mbb{Q}(q)$-algebra generated by:
\begin{align}\label{RTT-relation-1}
\mbb{Q}(q)\langle e_{\pm[i;j)},\psi_{s}^{\pm1},c^{\pm1}\rangle^{s\in\{1,\cdots,n\}}_{(i<j)\in\mbb{Z}^2/(n,n)\mbb{Z}}
\end{align}
The generators $e_{\pm[i;j)}$ satisfy the well-known RTT relations:
\begin{align}\label{RTT-relation-2}
R(\frac{z}{w})T_{1}^{+}(z)T_{2}^{+}(w)=T_{2}^{+}(w)T_{1}^{+}(z)R(\frac{z}{w})
\end{align}
\begin{align}\label{RTT-relation-3}
R(\frac{z}{w})T_{1}^{-}(z)T_{2}^{-}(w)=T_{2}^{-}(w)T_{1}^{-}(z)R(\frac{z}{w})
\end{align}
\begin{align}\label{RTT-relation-4}
R(\frac{z}{wc})T_{2}^{-}(w)T_{1}^{+}(z)=T_{1}^{+}(z)T_{2}^{-}(w)R(\frac{zc}{w}).
\end{align}

Here $T^{\pm}(z)$ are the generating functions:
\begin{equation*}
\begin{gathered}
T^{+}(z)=\sum_{1 \leq i \leq n}^{i \leq j} e_{[i ; j)} \psi_i \cdot E_{j \bmod n, i} z^{\left\lceil\frac{j}{n}\right\rceil-1} \\
T^{-}(z)=\sum_{1 \leq i \leq n}^{i \leq j} e_{-[i ; j)} \psi_i^{-1} \cdot E_{i, j \bmod n} z^{-\left\lceil\frac{j}{n}\right\rceil+1}
\end{gathered}
\end{equation*}
and $R(z/w)$ is the standard $R$-matrix for $\mf{gl}_{n}$:
\begin{equation*}
R\left(\frac{z}{w}\right)=\sum_{1 \leq i, j \leq n} E_{i, i} \otimes E_{j, j}\left(\frac{z q-w q^{-1}}{w-z}\right)^{\delta_j^i}+\left(q-q^{-1}\right) \sum_{1 \leq i \neq j \leq n} E_{i,j} \otimes E_{j,i} \frac{w^{\delta_{i>j}} z^{\delta_{i<j}}}{w-z}
\end{equation*}
where $[i,j)$ is a vector in $\mbb{R}^n$ defined as:
\begin{align*}
[i,j):=e_i+\cdots+e_{j-1}
\end{align*}
and here $e_k$ is the unit vector in $\mbb{R}^n$ at $k$-th direction with $0\leq k\leq n-1$, and $e_k$ is understood as $e_{k\text{ mod }n}$.

In addition, we have the relations:
\begin{equation*}
\psi_k \cdot e_{ \pm[i ; j)}=q^{ \pm\left(\delta_k^i-\delta_k^j\right)} e_{ \pm[i ; j)} \cdot \psi_k \quad \forall \operatorname{arcs}[i ; j) \text { and } k \in \mathbb{Z}.
\end{equation*}

The algebra $U_{q}(\hat{\mf{gl}}_{n})$ has the coproduct structure given by:
\begin{align*}
\Delta(T^{+}(z))=T^{+}(z)\otimes T^{+}(zc_1),\qquad\Delta(T^{-}(z))=T^{-}(zc_2)\otimes T^{-}(z)
\end{align*}
with $c_1:=c\otimes1$ and $c_2=1\otimes c$. Its antipode map is defined as:
\begin{align*}
S(T^{\pm}(z))=(T^{\pm}(z))^{-1}
\end{align*}
which may equivalently be written as:
\begin{equation*}
\begin{gathered}
S(e_{[i,j)})=\sum_{1 \leq i \leq n}^{i \leq j}(-1)^{j-i} f_{[i ; j)} \psi_j \cdot E_{j \bmod n, i} z^{\left\lceil\frac{j}{n}\right\rceil-1} \\
S(e_{-[i,j)})=\sum_{1 \leq i \leq n}^{i \leq j}(-1)^{j-i} f_{-[i ; j)} \psi_j^{-1} \cdot E_{i, j \bmod n} z^{-\left\lceil\frac{j}{n}\right]+1}.
\end{gathered}
\end{equation*}

Unwinding the relation \eqref{RTT-relation-1}\eqref{RTT-relation-2}\eqref{RTT-relation-3}, we obtain
\begin{align*}
\frac{e_{\pm[a,c)}e_{\pm[b,d)}}{q^{\delta^b_a-\delta^d_b+\delta^d_a}}-\frac{e_{\pm[b,d)}e_{\pm[a,c)}}{q^{\delta^b_c-\delta^d_b+\delta^d_c}}=(q-q^{-1})[\sum_{a\leq x<c}^{x\equiv d}e_{\pm[b,c+d-x)}e_{\pm[a,x)}-\sum_{a<x\leq c}^{x\equiv b}e_{\pm[x,c)}e_{\pm[a+b-x,d)}]
\end{align*}

\begin{align}\label{gln-relation}
[e_{[a,c)},e_{-[b,d)}]=(q-q^{-1})[\sum^{x\equiv b}_{a\leq x<c}\frac{e_{[c+b-x,d)}e_{[a,x)}}{q^{\delta^b_c+\delta^a_c-\delta^a_b}}\frac{\psi_x}{\psi_c}-\sum_{a<x\leq c}^{x\equiv d}\frac{e_{[x,c)}e_{-[b,a+d-x)}}{q^{-\delta^b_a+\delta^d_b-\delta^d_a}}\frac{\psi_x}{\psi_a}].
\end{align}

In \cite{N15} it was shown that the isomorphism $U_{q}(\hat{\mf{gl}}_{n})\cong U_{q}(\hat{\mf{sl}}_{n})\otimes U_{q}(\hat{\mf{gl}}_{1})/(c\otimes1-1\otimes c)$ is given by the following:
\begin{align*}
e_{[i;i+1)}=x_{i}^{+}(q-q^{-1})\\
e_{-[i;i+1)}=x_{i}^{-}(q^{-2}-1).
\end{align*}

For the generators $\{p_{k}\}_{k\in\mbb{Z}}$ of the Heisenberg algebra $U_{q}(\hat{\mf{gl}}_{1})\subset U_{q}(\hat{\mf{gl}}_n)$, they satisfy the following commutation relations:
\begin{align*}
[p_k,p_{-k}]=\frac{c^{nk}-c^{-nk}}{n_k},\qquad n_k=\frac{q^k-q^{-k}}{k}.
\end{align*}

It can be packaged into the following proposition proved by Negu\c{t} in \cite{N19}:
\begin{prop}
There exist constants $\alpha_1,\alpha_2,\cdots\in\mbb{Q}(q)$ such that:
\begin{equation*}
f_{ \pm[i ; j)}=\sum_{k=0}^{\left\lfloor\frac{j-i}{n}\right\rfloor} e_{ \pm[i ; j-n k)} g_{ \pm k}.
\end{equation*}

Here $\sum_{k=0}^{\infty}g_{\pm k}x^{k}=\exp(\sum_{k=1}^{\infty}\alpha_{k}p_{\pm k}x^k)$
\end{prop}

This formula allows us to obtain the expression for $p_{\pm k}$ iteratively.

\subsection{Shuffle algebra realization of $U_{q,t}(\hat{\hat{\mf{sl}}}_n)$}
Here we review the reconstruction of the quantum toroidal algebra via the shuffle algebra and the induced slope factorization of the quantum toroidal algebra. For details see \cite{N15}.

Fix the field of fractions $\mbb{F}=\mbb{Q}(q,t)$ and consider the space of symmetric rational functions:
\begin{align*}
\widehat{\text{Sym}}(V):=\bigoplus_{\mbf{k}=(k_1,\cdots,k_{n})\in\mbb{N}^n}\mbb{F}(\cdots,z_{i1},\cdots,z_{ik_i},\cdots)^{\text{Sym}}_{1\leq i\leq n}
\end{align*}
Here "Sym" means symmetrisation of the rational function for each $z_{i1},\cdots,z_{ik_i}$. We endow the vector space with the shuffle product:
\begin{align*}
F*G=\text{Sym}[\frac{F(\cdots,z_{ia},\cdots)G(\cdots,z_{jb},\cdots)}{\mbf{k}!\mbf{l}!}\prod_{1\leq a\leq k_i}^{1\leq i\leq n}\prod_{k_{j}+1\leq b\leq k_{j}+l_{j}}^{1\leq j\leq n}\zeta_{ij}(\frac{z_{ia}}{z_{jb}})].
\end{align*}

We define the subspace $\mc{S}^{\pm}_{\mbf{k}}\subset\widehat{\text{Sym}}(V)$ by:
\begin{align*}
\mc{S}^{+}:=\{F(\cdots,z_{ia},\cdots)=\frac{r(\cdots,z_{ia},\cdots)}{\prod_{1\leq a\neq b\leq k_i}^{1\leq i\leq n}(qz_{ia}-q^{-1}z_{ib})}\}
\end{align*}
where $r(\cdots,z_{i1},\cdots,z_{ik_i},\cdots)^{1\leq i\leq n}_{1\leq a\leq k_i}$ is any symmetric Laurent polynomial that satisfies the wheel conditions:
\begin{align*}
r(\cdots,z_{i1},\cdots,z_{ik_i},\cdots)|_{z_{ia}=q^{-1},z_{ib}=q,z_{i\pm1,c}=t^{\pm1}}=0
\end{align*}
for any three variables of colors $i,i\pm1,i$ respectively.

The shuffle algebra $\mc{S}^+$ has a natural bigrading given by the number of the variables in each color $\mbf{k}\in\mbb{N}^n$ and the homogeneous degree of the rational functions $d\in\mbb{Z}$:
\begin{align*}
\mc{S}^+=\bigoplus_{(\mbf{k},d)\in\mbb{N}^n\times\mbb{Z}}\mc{S}^+_{\mbf{k},d}.
\end{align*}

Similarly we can define the negative shuffle algebra $\mc{S}^{-}:=(\mc{S}^{+})^{op}$ which is the same as $\mc{S}^{+}$ with the opposite shuffle product. Now we slightly enlarge the positive and negative shuffle algebras by the generators $\{\varphi_{i,d}^{\pm}\}_{1\leq i\leq n}^{d\geq0}$:
\begin{align*}
\mc{S}^{\geq}=\langle\mc{S}^+,\{(\varphi^+_{i,d})^{d\geq0}_{1\leq i\leq n}\}\rangle/(\sim),\qquad \mc{S}^{\leq}=\langle\mc{S}^-,\{(\varphi^-_{i,d})^{d\geq0}_{1\leq i\leq n}\}\rangle/(\sim)
\end{align*}
and here $\varphi^{\pm}_{i,d}$ commute with themselves and have the relation with $\mc{S}^{\pm}$ as follows:
\begin{align*}
\varphi^{+}_i(w)F=F\varphi^+_i(w)\prod_{1\leq a\leq k_j}^{1\leq j\leq n}\frac{\zeta_{ij}(w/z_{ja})}{\zeta_{ji}(z_{ja}/w)}
\end{align*}
\begin{align*}
\varphi^{-}_i(w)G=G\varphi^-_i(w)\prod_{1\leq a\leq k_j}^{1\leq j\leq n}\frac{\zeta_{ji}(z_{ja}/w)}{\zeta_{ij}(w/z_{ja})}.
\end{align*}

The shuffle algebra $\mc{S}:=\mc{S}^{\leq}\hat{\otimes}\mc{S}^{\geq}$ is defined with respect to the bialgebra pairing as explained in \cite{N15}. Also the coproduct $\Delta:\mc{S}\rightarrow\mc{S}\hat{\otimes}\mc{S}$ is given by:
\begin{align*}
&\Delta(\varphi_{i}^{\pm}(w))=\varphi^{\pm}_{i}(w)\otimes\varphi^{\pm}_{i}(w)\\
&\Delta(F)=\sum_{0\leq\mbf{l}\leq\mbf{k}}\frac{[\prod_{1\leq j\leq n}^{b>l_{j}}\varphi_{j}^+(z_{jb})]F(\cdots,z_{i1},\cdots,z_{il_i}\otimes z_{i,l_i+1},\cdots,z_{ik_i},\cdots)}{\prod_{1\leq i\leq n}^{a\leq l_i}\prod_{1\leq j\leq n}^{b>l_j}\zeta_{ji}(z_{jb}/z_{ia})},\qquad F\in\mc{S}^+\\
&\Delta(G)=\sum_{0\leq\mbf{l}\leq\mbf{k}}\frac{G(\cdots,z_{i1},\cdots,z_{il_i}\otimes z_{i,l_i+1},\cdots,z_{ik_i},\cdots)[\prod_{1\leq j\leq n}^{b>l_{j}}\varphi_{j}^-(z_{jb})]}{\prod_{1\leq i\leq n}^{a\leq l_i}\prod_{1\leq j\leq n}^{b>l_j}\zeta_{ij}(z_{ia}/z_{jb})},\qquad G\in\mc{S}^-.
\end{align*}

In human language, the above coproduct formula means the following: Expand the right hand side in the region $\lvert z_{ia}\lvert<<\lvert z_{jb}\lvert$ for all $a\leq l_i$ and $b>l_{j}$, and then place all monomials in $\{z_{ia}\}_{a\leq l_i}$ to the left of the $\otimes$ symbol and all monomials in $\{z_{jb}\}_{b>l_j}$ to the right of the $\otimes$ symbol. Also we have the antipode map $S:\mc{S}\rightarrow\mc{S}$ which is an anti-homomorphism of both algebras and coalgebras:
\begin{align*}
&S(\varphi^+_i(z))=(\varphi^+_i(z))^{-1},\qquad S(F)=[\prod_{1\leq a\leq k_i}^{1\leq i\leq n}(-\varphi_i^+(z_{ia}))^{-1}]*F\\
&S(\varphi^-_i(z))=(\varphi^-_i(z))^{-1},\qquad S(G)=G*[\prod_{1\leq a\leq k_i}^{1\leq i\leq n}(-\varphi_i^-(z_{ia}))^{-1}].
\end{align*}

The following theorem has been proved by \cite{N15}:
\begin{thm}
There is a bigraded isomorphism of bialgebras
\begin{align*}
Y:U_{q,t}(\hat{\hat{\mf{sl}}}_{n})\rightarrow\mc{S}
\end{align*}
given by:
\begin{align*}
Y(\varphi_{i,d}^{\pm})=\varphi^{\pm}_{i,d},\qquad Y(e_{i,d}^{\pm})=\frac{z_{i1}^d}{[q^{-2}]}.
\end{align*}
\end{thm}

\subsection{Slope subalgebra of the shuffle algebra}
A remarkable property of the shuffle algebra is that it admits a slope factorization.

For the definition of the slope subalgebra see section $5.2$ in \cite{N15}.

It is known that the slope subalgebra $\mc{B}_{\mbf{m}}$ for the quantum toroidal algebra $U_{q,t}(\hat{\hat{\mf{sl}}}_{n})$ is isomorphic to:
\begin{align}\label{rootquantum}
\mc{B}_{\mbf{m}}\cong\bigotimes^{g}_{h=1}U_{q}(\widehat{\mf{gl}}_{l_h}).
\end{align}

For the proof of the above isomorphism \eqref{rootquantum} see \cite{N15}. For the rules of computing $l_h$, one can refer to Section 5.3 in \cite{N15}. The isomorphism is constructed in the following way. For $\mbf{m}\in\mbb{Q}^{n}$, it is proved in \cite{N15} that as the shuffle algebra, $\mc{B}_{\mbf{m}}$ is the Drinfeld double of the following:
\begin{align*}
\mc{B}_{\mbf{m}}:=\mc{B}_{\mbf{m}}^{+}\otimes\mbb{C}[\varphi_{i,0}^{\pm}]\otimes\mc{B}_{\mbf{m}}^{-}/(\text{relations}).
\end{align*}

Here $\mc{B}_{\mbf{m}}^{+}$ and $\mc{B}_{\mbf{m}}^{-}$ are defined as $(5.29)$ and $(5.30)$ in \cite{N15}, and we omit the details of the definition here. In this section we use the generators of the slope subalgebra $\mc{B}_{\mbf{m}}$ as the definition of the slope subalgebra.

Let us define the generators of the slope algebra $\mc{B}_{\mbf{m}}$. For $\mbf{m}\cdot[i;j)\in\mbb{Z}$, we denote the following elements:
\begin{equation}\label{positive-generators}
P_{\pm[i ; j)}^{\pm \mathbf{m}}=\operatorname{Sym}\left[\frac{\prod_{a=i}^{j-1} z_a^{\left\lfloor m_i+\ldots+m_a\right\rfloor-\left\lfloor m_i+\ldots+m_{a-1}\right\rfloor}}{t^{\text {ind } d_{[i, j)}^m} q^{i-j} \prod_{a=i+1}^{j-1}\left(1-\frac{q_2 z_a}{z_{a-1}}\right)} \prod_{i \leq a<b<j} \zeta_{ba}(\left(\frac{z_b}{z_a}\right)\right]\in\mc{S}^{\pm}
\end{equation}
\begin{equation}\label{negative generators}
Q_{\mp[i ; j)}^{\pm \mathbf{m}}=\operatorname{Sym}\left[\frac{\prod_{a=i}^{j-1} z_a^{\left\lfloor m_i+\ldots+m_{a-1}\right\rfloor-\left\lfloor m_i+\ldots+m_a\right\rfloor}}{t^{-\mathrm{ind}_{[i, j)}^m} \prod_{a=i+1}^{j-1}\left(1-\frac{q_1 z_{a-1}}{z_a}\right)} \prod_{i \leq a<b<j} \zeta_{ab}\left(\frac{z_a}{z_b}\right)\right]\in\mc{S}^{\mp}.
\end{equation}

Here $\text{ind}^{\mbf{m}}_{[i;j)}$ is defined as:
\begin{align*}
\text{ind}^{\mbf{m}}_{[i;j)}=\sum_{a=i}^{j-1}(m_i+\cdots+m_a-\lfloor m_i+\cdots+m_{a-1}\rfloor).
\end{align*}

The positive and negative part of the slope subalgebra $\mc{B}_{\mbf{m}}^{\pm}$ for the shuffle algebra $\mc{S}^{\pm}$ can be defined as the algebra generated by $\{P^{\mbf{m}}_{\pm[i;j)}\}_{i\leq j}$ and $\{Q^{\mbf{m}}_{\pm[i;j)}\}_{i\leq j}$. The slope subalgebra $\mc{B}_{\mbf{m}}$ is the similarly the Drinfeld double of $\mc{B}_{\mbf{m}}^{\pm}$ with the neutral elements $\{\varphi_{\pm[i;j)}\}$.

We can check that the antipode map $S_{\mbf{m}}:\mc{B}_{\mbf{m}}\rightarrow\mc{B}_{\mbf{m}}$ has the following relation:
\begin{align*}
S_{\mbf{m}}(P^{\mbf{m}}_{[i;j)})=Q^{\mbf{m}}_{[i;j)},\qquad S_{\mbf{m}}(Q^{\mbf{m}}_{-[i;j)})=P^{\mbf{m}}_{-[i;j)}
\end{align*}

\begin{equation*}
\Delta_{\mathbf{m}}\left(P_{[i ; j)}^{\mathbf{m}}\right)=\sum_{a=i}^j P_{[a ; j)}^{\mathbf{m}} \varphi_{[i ; a)} \otimes P_{[i ; a)}^{\mathbf{m}} \quad \Delta_{\mathbf{m}}\left(Q_{[i ; j)}^{\mathrm{m}}\right)=\sum_{a=i}^j Q_{[i ; a)}^{\mathbf{m}} \varphi_{[a ; j)} \otimes Q_{[a ; j)}^{\mathbf{m}}
\end{equation*}

\begin{equation*}
\Delta_{\mathbf{m}}\left(P_{-[i ; j)}^{\mathbf{m}}\right)=\sum_{a=i}^j P_{-[a ; j)}^{\mathbf{m}} \otimes P_{-[i ; a)}^{\mathbf{m}} \varphi_{-[a ; j)} \quad \Delta_{\mathbf{m}}\left(Q_{-[i ; j)}^{\mathbf{m}}\right)=\sum_{a=i}^j Q_{-[i ; a)}^{\mathbf{m}} \otimes Q_{-[a ; j)}^{\mathbf{m}} \varphi_{-[i ; a)}.
\end{equation*}

Here $\Delta_{\mbf{m}}$ is defined as $(5.27)$ and $(5.28)$ in \cite{N15}. The isomorphism \eqref{rootquantum} is given by:
\begin{align*}
e_{[i;j)}=P^{\mbf{m}}_{[i;j)_{h}},\qquad e_{-[i;j)}=Q^{\mbf{m}}_{-[i;j)_h},\qquad\varphi_{k}=\varphi_{[k;v_{\mbf{m}}(k))}.
\end{align*}
\subsection{Geometric action on $K_{T}(M(\mbf{w}))$}
We first fix the notation for the Nakajima quiver varieties.

Given a quiver $Q=(I,E)$, consider the following quiver representation space
\begin{align*}
\text{Rep}(\mbf{v},\mbf{w}):=\bigoplus_{h\in E}\text{Hom}(V_{i(h)},V_{o(h)})\oplus\bigoplus_{i\in I}\text{Hom}(V_i,W_i)
\end{align*}
and here $\mbf{v}=(\text{dim}(V_1),\cdots,\text{dim}(V_{|I|}))$ is the dimension vector for the vector spaces at the vertex, $\mbf{w}=(\text{dim}(W_1),\cdots,\text{dim}(W_I))$ is the dimension vector for the framing vector spaces.

Denote $G_{\mbf{v}}:=\prod_{i\in I}GL(V_i)$ and $G_{\mbf{w}}:=\prod_{i\in I}GL(W_i)$. There is a natural action of $G_{\mbf{v}}$ and $G_{\mbf{w}}$ on $\text{Rep}(\mbf{v},\mbf{w})$, and thus a natural Hamiltonian action on $T^*\text{Rep}(\mbf{v},\mbf{w})$ with respect to the standard symplectic form $\omega$. The corresponding moment map is:
\begin{align*}
\mu:T^*\text{Rep}(\mbf{v},\mbf{w})\rightarrow\mf{g}_{\mbf{v}}^*.
\end{align*}

One can define \textbf{Nakajima variety}:
\begin{align}\label{Quiver-variety}
M_{\theta,0}(\mbf{v},\mbf{w}):=\mu^{-1}(0)//_{\theta}G_{\mbf{v}}
\end{align}
with all the $\theta$-stable $G_{\mbf{v}}$-orbits in $\mu^{-1}(0)$. Here we choose the stability condition $\theta=(1,\cdots,1)$.

We will abbreviate $M(\mbf{v},\mbf{w})$ as the Nakajima quiver variety defined in \eqref{Quiver-variety}.

In this paper we only consider the affine type $A$ quiver, i.e. the edges are given by $e=(i,i+1)$ with $i=i\text{ mod }(n)$. i.e. the quiver is:

\begin{tikzcd}
	&&& {W_k} \\
	{W_{k-1}} &&& {V_{k}} &&& {W_{k+1}} \\
	& {V_{k-1}} &&&& {V_{k+1}} \\
	& \cdots &&&& \cdots \\
	& {V_2} &&&& {V_n} \\
	{W_2} &&& {V_1} &&& {W_n} \\
	&&& {W_1} \\
	&&& {}
	\arrow[curve={height=6pt}, from=1-4, to=2-4]
	\arrow[curve={height=6pt}, from=2-1, to=3-2]
	\arrow[curve={height=6pt}, from=2-4, to=1-4]
	\arrow[curve={height=6pt}, from=2-4, to=3-2]
	\arrow[curve={height=6pt}, from=2-4, to=3-6]
	\arrow[curve={height=6pt}, from=2-7, to=3-6]
	\arrow[curve={height=6pt}, from=3-2, to=2-1]
	\arrow[curve={height=6pt}, from=3-2, to=2-4]
	\arrow[curve={height=6pt}, from=3-2, to=4-2]
	\arrow[curve={height=6pt}, from=3-6, to=2-4]
	\arrow[curve={height=6pt}, from=3-6, to=2-7]
	\arrow[curve={height=6pt}, from=3-6, to=4-6]
	\arrow[curve={height=6pt}, from=4-2, to=3-2]
	\arrow[curve={height=6pt}, from=4-2, to=5-2]
	\arrow[curve={height=6pt}, from=4-6, to=3-6]
	\arrow[curve={height=6pt}, from=4-6, to=5-6]
	\arrow[curve={height=6pt}, from=5-2, to=4-2]
	\arrow[curve={height=6pt}, from=5-2, to=6-1]
	\arrow[curve={height=6pt}, from=5-2, to=6-4]
	\arrow[curve={height=6pt}, from=5-6, to=4-6]
	\arrow[curve={height=6pt}, from=5-6, to=6-4]
	\arrow[curve={height=6pt}, from=5-6, to=6-7]
	\arrow[curve={height=6pt}, from=6-1, to=5-2]
	\arrow[curve={height=6pt}, from=6-4, to=5-2]
	\arrow[curve={height=6pt}, from=6-4, to=5-6]
	\arrow[curve={height=6pt}, from=6-4, to=7-4]
	\arrow[curve={height=6pt}, from=6-7, to=5-6]
	\arrow[curve={height=6pt}, from=7-4, to=6-4]
\end{tikzcd}

The action of $\mbb{C}^*_{q}\times\mbb{C}^*_t\times G_{\mbf{w}}$ on $M(\mbf{v},\mbf{w})$ is given by:
\begin{align*}
(q,t,U_i)\cdot(X_e,Y_e,A_i,B_i)_{e\in E,i\in I}=(qtX_e,qt^{-1}Y_e,qA_iU_i,qU_i^{-1}B_i).
\end{align*}

Choose $T_{\mbf{w}}\subset G_{\mbf{w}}$ the maximal torus of $G_{\mbf{w}}$. The fixed-point set of the torus action $\mbb{C}^*_q\times\mbb{C}^*_t\times T_{\mbf{w}}$ on $M_{\theta,0}(\mbf{v},\mbf{w})$ is indexed by the $|\mbf{w}|$-partitions:
\begin{align*}
\bm{\lambda}=(\lambda_1,\cdots,\lambda_{|\mbf{w}|})
\end{align*}
such that $|\bm{\lambda}|=|\mbf{v}|$, and such that the set of the boxes
\begin{align*}
\{\square\in\bm{\lambda}|c_{\square}\equiv i(\text{ mod } n)\}
\end{align*}
has size $v_{i}$. For each box $\square\in(\lambda_1,\cdots.\lambda_{|\mbf{w}|})$, we define the following two character functions:
\begin{align*}
\chi_{\square}=u_iq^{x+y+1}t^{x-y},\qquad\square\in\lambda_i
\end{align*}

\begin{align*}
\chi^{ch}_{\square}=v_i+\hbar_1(x+1)+\hbar_2(y+1).
\end{align*}

The first one is often used in the equivariant $K$-theory. The second one is often used in the equivariant cohomology theory.

 If we choose the cocharacter $\sigma:\mbb{C}^*\rightarrow G_{\mbf{w}}$ such that $\mbf{w}=u_1\mbf{w}_1+\cdots+u_k\mbf{w}_k$, we have that:
\begin{align*}
M(\mbf{v},\mbf{w})^{\sigma}=\bigsqcup_{\mbf{v}_1+\cdots+\mbf{v}_k=\mbf{v}}M(\mbf{v}_1,\mbf{w}_1)\times\cdots\times M(\mbf{v}_k,\mbf{w}_k).
\end{align*}

Using the notation $T:=\mbb{C}_q^*\times\mbb{C}_t^*\times T_{\mbf{w}}$, we define:
\begin{align*}
K_{T}(M(\mbf{w})):=\bigoplus_{\mbf{v}}K_{T}(M(\mbf{v},\mbf{w}))_{loc}.
\end{align*}

We denote $I_{\bm{\lambda}}$ as the skyscraper sheaf corresponding to the fixed point $\bm{\lambda}$. We define the fixed-point basis by:
\begin{align*}
|\bm{\lambda}\rangle=\frac{I_{\bm{\lambda}}}{[T_{\bm{\lambda}}M(\mbf{v},\mbf{w})]}.
\end{align*}

Thus we can choose the fixed point basis $|\bm{\lambda}\rangle$ of the corresponding partition $\bm{\lambda}$ to span the vector space $K_{T}(M(\mbf{w}))$.

Since $K_T(M(\mbf{v},\mbf{w}))$ is generated by tautological classes and $K_{T}(pt)$, we may denote $x_{ij}$ with $1\leq i\leq n$ and $1\leq j\leq v_{i}$ as the Chern roots. Given $f\in K_{T}(M(\mbf{v},\mbf{w}))$ the colored-symmetric Laurent polynomial in Chern roots (or Laurent polynomial of tautological classes), we define $\overline{f}$ as:
\begin{align}\label{fixed-point-expression}
\overline{f}:=\sum_{\bm{\lambda}}f|_{\bm{\lambda}}|\bm{\lambda}\rangle.
\end{align}

Here we briefly review the geometric action of $U_{q,t}(\hat{\hat{\mf{sl}}}_{n})$ on $K_{T}(M(\mbf{w}))$. The construction is based on Nakajima's simple correspondence. 

Consider a pair of vectors $(\mbf{v}_+,\mbf{v}_-)$ such that $\mbf{v}_{+}=\mbf{v}_{-}+\mbf{e}_{i}$. There is a simple correspondence
\begin{align*}
Z(\mbf{v}_+,\mbf{v}_{-},\mbf{w})\hookrightarrow M(\mbf{v}_+,\mbf{w})\times M(\mbf{v}_-,\mbf{w})
\end{align*}
parametrises pairs of quadruples $(X^{\pm},Y^{\pm},A^{\pm},B^{\pm})$ that respect a fixed collection of quotients $(V^+\rightarrow V^-)$ of codimension $\delta^{i}_{j}$ consisting of the semistable points on which the moment map $\mu$ vanishes for each $M(\mbf{v}_+,\mbf{w})$. The variety $Z(\mbf{v}_+,\mbf{v}_{-},\mbf{w})$ is smooth with a tautological line bundle:
\begin{align*}
\mc{L}|_{V^+\rightarrow V^-}=\text{Ker}(V_{i}^+\rightarrow V_{i}^{-})
\end{align*}
and the natural projection maps:
\begin{equation*}
\begin{tikzcd}
&Z(\mbf{v}_+,\mbf{v}_{-},\mbf{w})\arrow[ld,"\pi_{+}"]\arrow[rd,"\pi_{-}"]&\\
M(\mbf{v}_+,\mbf{w})&&M(\mbf{v}_-,\mbf{w}).
\end{tikzcd}
\end{equation*}

Using this we could construct the operator:
\begin{align*}
e_{i,d}^{\pm}:K_{T}(M(\mbf{v}^{\mp},\mbf{w}))\rightarrow K_{T}(M(\mbf{v}^{\pm},\mbf{w})),\qquad e_{i,d}^{\pm}(\alpha)=\pi_{\pm*}(\text{Cone}(d\pi_{\pm})\mc{L}^d\pi_{\mp}^{\pm}(\alpha)).
\end{align*}

Here $\text{Cone}(d\pi_{\pm})$ is defined in \cite{N15}. Taking all $\mbf{v}$, we obtain operators $e_{i,d}^{\pm}:K_{T}(M(\mbf{w}))\rightarrow K_{T}(M(\mbf{w}))$. Also we have the action of $\varphi^{\pm}_{i,d}$ given by the multiplication of the tautological class, which means that:
\begin{align*}
\varphi_{i}^{\pm}(z)=\overline{\frac{\zeta{(\frac{z}{X})}}{\zeta(\frac{X}{z})}}\prod^{u_{j}\equiv i}_{1\leq j\leq\mbf{w}}\frac{[\frac{u_j}{\hbar z}]}{[\frac{z}{\hbar u_{j}}]},\qquad\overline{\frac{\zeta{(\frac{z}{X})}}{\zeta(\frac{X}{z})}}:=\prod_{j\in I}\prod_{1\leq a\leq v_j}\overline{\frac{\zeta_{ij}(\frac{z}{x_{ja}})}{\zeta_{ji}(\frac{x_{ja}}{z})}}.
\end{align*}

In particular, the element $\varphi_{i,0}$ acts on $K_{T}(M(\mbf{v},\mbf{w}))$ as $q^{\alpha_{i}^T(\mbf{w}-C\mbf{v})}$.

The above construction gives the following well-known result:
\begin{thm}
For all $\mbf{w}\in\mbb{N}^r$, the operators $e_{i,d}^{\pm}$ and $\varphi^{\pm}_{i,d}$ give rise to an action of $U_{q,t}(\hat{\hat{\mf{sl}}}_{n})$ on $K_{T}(M(\mbf{w}))$.
\end{thm}

In terms of the shuffle algebra, we can give the explicit formula of the action of the quantum toroidal algebra $U_{q,t}(\hat{\hat{\mf{sl}}}_{n})$ on $K_{T}(M(\mbf{w}))$:

Given $F\in\mc{S}_{\mbf{k}}^+$, we have that
\begin{align}\label{shuffle-formula-1}
\langle\bm{\lambda}|F|\bm{\mu}\rangle=F(\chi_{\bm{\lambda}\backslash\bm{\mu}})\prod_{\blacksquare\in\bm{\lambda}\backslash\bm{\mu}}[\prod_{\square\in\bm{\mu}}\zeta_{ij}(\frac{\chi_{\blacksquare}}{\chi_{\square}})\prod_{i=1}^{\mbf{w}}[\frac{u_i}{q\chi_{\blacksquare}}]],\qquad i=\text{color of }\blacksquare, j=\text{color of }\square.
\end{align}

Similarly, for $G\in\mc{S}_{-\mbf{k}}^{-}$, we have
\begin{align}\label{shuffle-formula-2}
\langle\bm{\mu}|G|\bm{\lambda}\rangle=G(\chi_{\bm{\lambda}\backslash\bm{\mu}})\prod_{\blacksquare\in\bm{\lambda}\backslash\bm{\mu}}[\prod_{\square\in\bm{\lambda}}\zeta(\frac{\chi_{\square}}{\chi_{\blacksquare}})\prod_{i=1}^{\mbf{w}}[\frac{\chi_{\blacksquare}}{qu_i}]]^{-1},\qquad i=\text{color of }\blacksquare, j=\text{color of }\square.
\end{align}

Here $F(\chi_{\bm{\lambda}\backslash\bm{\mu}})$ and $G(\chi_{\bm{\lambda}\backslash\bm{\mu}})$ are the rational symmetric function evaluation at the box $\chi_{\square}$ with $\square\in\bm{\lambda}\backslash\bm{\mu}$.

\section{\textbf{Affine Yangian and Kac-Moody algebra}}\label{sec:_textbf_affine_yangian_and_kac_moody_algebra}
In this section we introduce the presentation of the affine Yangian that is suitable in our settings.

We fix $\mbb{F}:=\mbb{Q}(\hbar_1,\hbar_2)$ with $\hbar:=\hbar_1+\hbar_2$. We introduce the following color-dependent rational function:
\begin{align*}
\omega_{ij}(x_i-x_j)=\frac{(x_j-x_i-\hbar_1)^{\delta^i_{j-1}}(x_j-x_i-\hbar_2)^{\delta^{i}_{j+1}}}{(x_j-x_i)^{\delta^i_j}(x_j-x_i-\hbar_1-\hbar_2)^{\delta^i_j}}.
\end{align*}

Consider the vector space
\begin{align*}
\mc{V}^+:=\bigoplus_{\mbf{n}\in\mbb{N}^I}\mc{V}_{\mbf{n}},\qquad\mc{V}_{\mbf{n}}:=\mbb{F}[z_{i1},\cdots,z_{in_i}]^{sym}_{i\in I}.
\end{align*}

we endow $\mc{V}^+$ with the following \textbf{shuffle product}:
\begin{equation*}
\begin{aligned}
R(z_{i1},\cdots, z_{in_i})_{i\in I}*R'(z_{i1},\cdots,z_{in_i'})_{i\in I}=&\text{Sym}[\frac{R(z_{i1},\cdots,z_{in_i})_{i\in I}R'(z_{i,n_i+1},\cdots, z_{i,n_i+n_{i}'})_{i\in I}}{\mbf{n}!\mbf{n}'!}\times\\
&\times\prod_{i,j\in I}\prod_{a=1}^{n_i}\prod_{b=n_{j}+1}^{n_{j}+n_{j}'}\omega_{ij}(z_{ia}-z_{jb})].
\end{aligned}
\end{equation*}

We consider the subalgebra $\mc{A}^+\subset\mc{V}^+$ generated by:
\begin{align*}
\{z_{i1}^d\}_{i\in I,d\geq0}.
\end{align*}

In this setting, we define the positive affine Yangian $Y_{\hbar_1,\hbar_2}^+(\hat{\mf{sl}}_n)$ as $\mc{A}^+$. The generators $z_{i1}^d$ is denoted by $E_{i,d}^+$.

\subsection{Affine Yangians}
In this subsection we define the affine Yangian $Y_{\hbar_1,\hbar_2}(\hat{\mf{sl}}_n)$ by means of the additive shuffle algebra. Most of the details can be found in \cite{BT19}.

The affine Yangian $Y_{\hbar_1,\hbar_2}(\hat{\mf{sl}}_n)$ is defined as a vector space:
\begin{align}\label{Yangian-formula-1}
Y_{\hbar_1,\hbar_2}(\hat{\mf{sl}}_n):=\mc{S}^+\otimes\mbb{Q}[\xi_{i,n}]_{n\geq0}\otimes\mc{S}^{-}
\end{align}
Here $\mc{S}^-:=(\mc{S}^+)^{op}$. The relations are defined as follows:
\begin{align}\label{Yangian-formula-2}
&p_{ij}(z,\sigma_j^{+})\xi_i(z)E_{j,s}^+=-p_{ji}(\sigma_j^{+},z)E_{j,s}^+\xi_{i}(z)\\
&p_{ji}(\sigma_j^-,z)\xi_{i}(z)E_{j,s}^-=-p_{ij}(z,\sigma_j^-)E_{j,s}^-\xi_i(z).
\end{align}

Here $p_{ij}(z,w)$ is defined as follows:
\begin{align*}
p_{ij}(z,w)=
\begin{cases}
z-w-\hbar_1&i=j-1\\
z-w-\hbar_2&i=j+1\\
(w-z)(w-z-\hbar_1-\hbar_2)&i=j.
\end{cases}
\end{align*}

$\sigma_{i}^{\pm}:Y_{\hbar_1,\hbar_2}(\hat{\mf{sl}}_n)^{\pm}\rightarrow Y_{\hbar_1,\hbar_2}(\hat{\mf{sl}}_n)^{\pm}$ is the homomorphism given by $\xi_{j,r}\mapsto\xi_{j,r}$, $E_{j,r}^{\pm}\mapsto E_{j,r+\delta_{ij}}^{\pm}$ 

\begin{align*}
[E_{i}^+(z),E_{j}^-(w)]=-\hbar_2\delta_{ij}\frac{\xi_i(z)-\xi_j(w)}{z-w}.
\end{align*}

As well as the relation in $\mc{S}^+$ and $\mc{S}^-$ given by the shuffle algebra structure.

\subsubsection{Lie algebra}
The affine Yangian $Y_{\hbar_1,\hbar_2}(\hat{\mf{sl}}_n)$ contains the Lie algebra $\hat{\mf{gl}}_n$, where:
\begin{align*}
\hat{\mf{gl}}_n=(n\times n\text{ matrices with value in }\mbb{C}[z^{\pm1}])\oplus\mbb{C}\cdot\gamma
\end{align*}
with $\gamma$ being central in $\hat{\mf{gl}}_n$ and the Lie bracket is written as:
\begin{align*}
[Xz^k,Yz^l]=[X,Y]z^{k+l}+\delta^{0}_{k+l}k\cdot\text{Tr}(XY)\gamma.
\end{align*}

We express the generators of the affine Lie algebra $\hat{\mf{gl}}_n$ as follows:
\begin{align*}
E_{[i,j)}=E_{i\text{ mod }n,j\text{ mod }n}\cdot z^{[\frac{j-1}{n}]-[\frac{i-1}{n}]}\\
E_{-[i,j)}=E_{j\text{ mod }n,i\text{ mod }n}\cdot z^{[\frac{j-1}{n}]-[\frac{i-1}{n}]}.
\end{align*}

The Lie algebra $\hat{\mf{gl}}_n$ contains the Heisenberg Lie algebra $\hat{\mf{gl}}_1$ with elements:
\begin{align*}
\alpha_{k}:=\sum_{i=1}^nE_{[i,i+nk)}=\text{Id}\cdot z^k
\end{align*}
which satisfies the relation:
\begin{align*}
[\alpha_k,\alpha_l]=\delta^0_{k+l}kn\gamma.
\end{align*}

\subsection{Action on the equivariant cohomology of quiver varieties}

The algebra action of $Y_{\hbar_1,\hbar_2}(\hat{\mf{sl}}_n)$ on the equivariant cohomology of quiver varieties $H_{T}^*(M(\mbf{w}))$ can be obtained as the degeneration limit of the action of the quantum toroidal algebra action $U_{q,t}(\hat{\hat{\mf{sl}}}_n)$ on the equivariant $K$-theory of quiver varieties $K_{T}^*(M(\mbf{w}))$.

\begin{equation*}
\begin{aligned}
\langle\bm{\lambda}|F|\bm{\mu}\rangle=F(\chi_{\bm{\lambda}\backslash\bm{\mu}}^{ch})\prod_{\blacksquare\in\bm{\lambda}\backslash\bm{\mu}}[\prod_{\square\in\bm{\mu}}\omega(\chi_{\blacksquare}-\chi_{\square})\prod_{i=1}^{\mbf{w}}(v_i-\hbar_1-\hbar_2-\chi_{\blacksquare})].
\end{aligned}
\end{equation*}

Here $|\bm{\lambda}\rangle$ is the fixed-point basis in the localised equivariant cohomology $H_{T}^*(M(\mbf{w})):=\bigoplus_{\mbf{v}}H_{T}^*(M(\mbf{v},\mbf{w}))_{loc}$ of the affine type $A$ quiver varieties. Similarly, the action of the Cartan current $\xi(z)$ is given by the multiplication of the tautological class:

\begin{align*}
\xi(z)=(\overline{\frac{\omega{(z-X)}}{\omega(X-z)}})^{+}\prod^{u_{j}\equiv i}_{1\leq j\leq\mbf{w}}\frac{a_j-\hbar-z}{z-\hbar-a_j}.
\end{align*}

Here $(\overline{\frac{\omega{(z-X)}}{\omega(X-z)}})^{+}$ means expanding the rational function of $z$ in terms of the non-negative power of $z$.

\section{\textbf{Degeneration of the quantum toroidal algebra to affine Yangians}}\label{sec:_textbf_degeneration_of_the_quantum_toroidal_algebra_to_affine_yangians}
In this section we describe the degeneration of the quantum toroidal algebra to the affine Yangians.

First let us explain what does it mean when we say the word "degeneration".

We first need to use the following fact:
\begin{prop}
There is an algebra embedding:
\begin{align*}
U_{q,t}(\hat{\hat{\mf{sl}}}_{n})\hookrightarrow\prod_{\mbf{w}}\text{End}(K(\mbf{w}))
\end{align*}
given by the matrix coefficients of the generators of $U_{q,t}(\hat{\hat{\mf{sl}}}_{n})$.
\end{prop}
\begin{proof}
See \cite{N23} for the complete proof of the Proposition.
\end{proof}

Note that $K(\mbf{w})$ is the $K_{T_{\mbf{w}}}(pt)_{loc}\cong\mbb{Q}(q,t,u_1,\cdots,u_{\mbf{w}})$-module. Using the Chern character map:
\begin{equation*}
\begin{aligned}
&ch: K_{T_{\mbf{w}}}(pt)\rightarrow\widehat{H_{T_{\mbf{w}}}(pt)},\qquad\widehat{H_{T_{\mbf{w}}}(pt)}:=\mbb{Q}[[\hbar_1,\hbar_2,z_1,\cdots,z_{\mbf{w}},\kappa]]\\
&q_i\mapsto e^{\kappa \hbar_i},\qquad u_{i}\mapsto e^{\kappa z_i}.
\end{aligned}
\end{equation*}

It induces the following map:
\begin{equation*}
\begin{aligned}
&ch: K_{T_{\mbf{w}}}(M(\mbf{w}))\rightarrow\widehat{H_{T_{\mbf{w}}}(M(\mbf{w}))}
\end{aligned}
\end{equation*}
such that the $\widehat{H_{T_{\mbf{w}}}(pt)}$ module $\widehat{H_{T_{\mbf{w}}}(M(\mbf{w}))}$ is spanned by the fixed point basis with coefficients in $\widehat{H_{T_{\mbf{w}}}(pt)}$. 

Now we turn to the definition of the asymptotic expansion of elements in $\widehat{H_{T_{\mbf{w}}}(M(\mbf{w}))}$. Given $f\in Im(ch)$, now we do the Laurent expansion of $f$ with respect to $\kappa$. We denote $f^{ch}$ as the lowest order component of $f$ with respect to $\kappa$ if it exists.

The following fact is easy to prove:
\begin{lem}
For arbitrary $f\in Im(ch)$, $f^{ch}$ exists.
\end{lem}
\begin{proof}
If $f\in\text{Im}(ch)$, then $f$ can be written as:
\begin{align*}
f(e^{\kappa z_i},e^{\kappa \hbar_i})=\frac{P(e^{\kappa z_i},e^{\kappa \hbar_i})}{Q(e^{\kappa z_i},e^{\kappa \hbar_i})}
\end{align*}
where $P(e^{\kappa z_i},e^{\kappa \hbar_i})$ and $Q(e^{\kappa z_i},e^{\kappa \hbar_i})$ are polynomials in the variables $e^{\kappa z_i}$ and $e^{\kappa \hbar_i}$. Hence the lowest-order term $f^{ch}$ exists.
\end{proof}

\subsection{Degeneration of the quantum toroidal algebra to the affine Yangian}

Now we turn to the case of the degeneration of the quantum toroidal algebra $U_{q,t}(\hat{\hat{\mf{sl}}}_n)$ to the affine Yangian $Y_{\hbar_1,\hbar_2}(\hat{\mf{sl}}_n)$.

\begin{prop}
The quantum toroidal algebra $U_{q,t}(\hat{\hat{\mf{sl}}}_n)$ is generated by
\begin{align*}
\{e_{i,0}^{+},e_{i,0}^{-},\varphi_{i,0}^{\pm},\varphi^{\pm}_{i,\pm1}\}
\end{align*}
with the central elements.

Similarly, the affine Yangian $Y_{\hbar_1,\hbar_2}(\hat{\mf{sl}}_n)$ is generated by $\{e_{i,0}^{+},e_{i,0}^{-},H_{i,0}^{\pm},H_{i,\pm1}\}$ with the central elements.
\end{prop}

The following proposition will be vital in the degeneration of the quantum difference equation to the quantum differential equation:
\begin{prop}\label{degeneration-of-algebra}
The degeneration limit of $U_{q,t}(\hat{\hat{\mf{sl}}}_n)$ lies in $Y_{\hbar_1,\hbar_2}(\hat{\mf{sl}}_n)$.
\end{prop}
\begin{proof}
This is proved by comparing the degeneration limit of the matrix coefficients of the elements of $U_{q,t}(\hat{\hat{\mf{sl}}}_n)$ with the matrix coefficients of the elements of $Y_{\hbar_1,\hbar_2}(\hat{\mf{sl}}_n)$.
\end{proof}

\section{\textbf{Review of the quantum difference equation}}\label{sec:_textbf_review_of_the_quantum_difference_equation}
In this section we give a review of the construction of the quantum difference equation in \cite{Z23} for the quantum toroidal algebra $U_{q,t}(\hat{\hat{\mf{sl}}}_{n})$ and quiver varieties of affine type $A$.

The quantum toroidal algebra $U_{q,t}(\hat{\hat{\mf{sl}}}_{n})$ admits the slope factorization:
\begin{align*}
U_{q,t}(\hat{\hat{\mf{sl}}}_{n})=\bigotimes^{\rightarrow}_{\mu\in\mbb{Q}}\mc{B}_{\mbf{m}+\mu\bm{\theta}}.
\end{align*}

For each slope subalgebra $\mc{B}_{\mbf{m}}$, one can associate an element $J_{\mbf{m}}^{\pm}(z)\in\mc{B}_{\mbf{m}}^{\pm}\hat{\otimes}\mc{B}_{\mbf{m}}^{\mp}$ such that they satisfy the ABRR equation:
\begin{align}\label{ABRR-equation}
J_{\mbf{m}}^{+}(z)z_{(1)}^{-1}q^{\Omega}R_{\mbf{m}}^{+}=z_{(1)}^{-1}q^{\Omega}J_{\mbf{m}}^{+}(z),\qquad q^{\Omega}R_{\mbf{m}}^{-}z_{(1)}^{-1}J_{\mbf{m}}^{-}(z)=J_{\mbf{m}}^{-}(z)q^{\Omega}z_{(1)}^{-1}.
\end{align}

\subsection{Asymptotics of $J_{\mbf{m}}(z)$}
with $R_{\mbf{m}}^{-}:=(R_{\mbf{m}})_{21}$. The operator $J_{\mbf{m}}^{\pm}$ satisfies the following relation:
\begin{align*}
S_{\mbf{m}}\otimes S_{\mbf{m}}((J_{\mbf{m}}^{+}(\lambda))_{21})=J_{\mbf{m}}^{-}(\lambda).
\end{align*}

In \cite{Z23} we have shown the following lemma:
\begin{lem}\label{asymptotic-lemma}
\begin{align*}
\lim_{z\rightarrow\infty}J_{\mbf{m}}^+(z^{\bm{\theta}})=1,\qquad\lim_{z\rightarrow0}J_{\mbf{m}}^+(z^{\bm{\theta}})=R_{\mbf{m}}^{+}
\end{align*}
here $\bm{\theta}$ is the parameter above in the slope factorization of the quantum toroidal algebra. Here we require that $\bm{\theta}>0$
\end{lem}

The following proposition reveals that for all $\mbf{m}\in\mbb{R}^n$, $\lim_{p\rightarrow0}J_{\mbf{m}}(zp^s)$ is a locally constant function on $s$:
\begin{prop}
$\lim_{p\rightarrow0}J_{\mbf{m}}(zp^{s})$ is locally constant on the variable $s\in\mbb{R}^n$, and it jumps as $s$ crosses the following hyperplanes:
\begin{equation*}
\begin{aligned}
&\mbf{n}\cdot s=0,\qquad \mbf{n}\in(\mbb{Z}_{\geq0})^n-\{0\}
\end{aligned}
\end{equation*}
\end{prop}
\begin{proof}
By the ABRR equation \eqref{ABRR-equation}, if we write down $R_{\mbf{m}}^+=\text{Id}+\sum_{\mbf{n}\in\mbb{N}^I}R_{\mbf{m}|\mbf{n}}$, we have the solution $J_{\mbf{m}}^+=\text{Id}+\sum_{\mbb{n}\in\mbb{N}^I}J_{\mbf{m}|\mbf{n}}(z)$ written as:
\begin{equation*}
\begin{aligned}
J_{\mbf{m}|\mbf{n}}(zp^s)=&\frac{1}{z^{\mbf{n}}p^{\mbf{n}\cdot s}q^{k}-1}\sum_{\substack{\mbf{n}_1+\mbf{n}_2=\mbf{n}\\\mbf{n}_1<\mbf{n}}}J_{\mbf{m}|\mbf{n}_1}(zp^s)R_{\mbf{m}|\mbf{n}_2}\\
=&\sum_{k=1}^{\infty}\sum_{\substack{\mbf{n}_1+\cdots+\mbf{n}_{k+1}=\mbf{n}\\\mbf{n}_{k}<\mbf{n}_{k-1}<\cdots<\mbf{n}_1<\mbf{n}}}\frac{R_{\mbf{m}|\mbf{n}_2+\cdots+\mbf{n}_{k+1}}R_{\mbf{m}|\mbf{n}_3+\cdots+\mbf{n}_{k+1}}\cdots R_{\mbf{m}|\mbf{n}_{k+1}}}{(z^{\mbf{n}}p^{\mbf{n}\cdot s}q^{k}-1)(z^{\mbf{n}_1}p^{\mbf{n}_1\cdot s}q^{k}-1)\cdots(z^{\mbf{n}_k}p^{\mbf{n}_k\cdot s}q^k-1)}.
\end{aligned}
\end{equation*}

So as $p\rightarrow0$, if $\mbf{n}_i\cdot s>0$ for every $1\leq i\leq k$, the term $\lim_{p\rightarrow0}\frac{R_{\mbf{m}|\mbf{n}_2+\cdots+\mbf{n}_{k+1}}R_{\mbf{m}|\mbf{n}_3+\cdots+\mbf{n}_{k+1}}\cdots R_{\mbf{m}|\mbf{n}_{k+1}}}{(z^{\mbf{n}}p^{\mbf{n}\cdot s}q^{k}-1)(z^{\mbf{n}_1}p^{\mbf{n}_1\cdot s}q^{k}-1)\cdots(z^{\mbf{n}_k}p^{\mbf{n}_k\cdot s}q^k-1)}$ is nonzero and equals $R_{\mbf{m}|\mbf{n}_2+\cdots+\mbf{n}_{k+1}}R_{\mbf{m}|\mbf{n}_3+\cdots+\mbf{n}_{k+1}}\cdots R_{\mbf{m}|\mbf{n}_{k+1}}$. Otherwise the term will be zero. Thus combining every term, one can see that $\lim_{p\rightarrow0}J_{\mbf{m}}(zp^s)$ will only change every time when $s$ goes across the hyperplane $\mbf{n}\cdot s=0$.
\end{proof}

The monodromy operator is defined as:
\begin{align}\label{defn-of-quantum-difference-operator}
\mbf{B}_{\mbf{m}}(z)=m((1\otimes S_{\mbf{m}})(J_{\mbf{m}}^{-}(zp^{-\mbf{m}})^{-1}))|_{z\mapsto zq^{\kappa}}.
\end{align}

Here $\kappa=\frac{C\mbf{v}-\mbf{w}}{2}$ if $\mbf{B}_{\mbf{m}}(z)$ acts on $K_{T}(M(\mbf{v},\mbf{w}))$.

Let $\mc{L}\in Pic(X)$ be a line bundle. Now we fix a slope $s\in H^2(X,\mbb{R})$ and choose a path in $H^2(X,\mbb{R})$ from $s$ to $s-\mc{L}$. This path crosses finitely many slope points in the order $\{\mbf{m}_1,\mbf{m}_2,\cdots,\mbf{m}_m\}$.  And for this choice of a slope, line bundle and a path we associate the following operator:
\begin{align}\label{defnofqdeoperator}
\mbf{B}_{\mc{L}}^{s}(\lambda)=\mc{L}\mbf{B}_{\mbf{m}_m}(\lambda)\cdots\mbf{B}_{\mbf{m}_1}(\lambda)
\end{align}

We define the $q$-difference operators:
\begin{align*}
\mc{A}^{s}_{\mc{L}}=T_{\mc{L}}^{-1}\mbf{B}^s_{\mc{L}}(\lambda).
\end{align*}

It has been proved in \cite{Z23} that the $q$-difference operator $\mc{A}^{s}_{\mc{L}}$ is independent of the choice of the path from $s$ to $s-\mc{L}$ if we choose the generic path.

The well-definedness and the rigidity of the definition of the $q$-difference operators $\mc{A}^{s}_{\mc{L}}$ lies in the following lemma:
\begin{lem}
For generic $s$ in $\text{Pic}(M(\mbf{v},\mbf{w}))\otimes\mbb{R}$, there exists a small open neighborhood $V_s$ of $s$ such that for every point $s'\in V_s$:
\begin{align*}
\mc{A}^{s'}_{\mc{L}}=\mc{A}^{s}_{\mc{L}}.
\end{align*}
\end{lem}
\begin{proof}
Note that for arbitrary path from $s$ to $s-\mc{L}$, there are only finitely many slope points $\{\mbf{m}_1,\cdots,\mbf{m}_n\}$ on the path giving the nontrivial action on $K_{T}(M(\mbf{v},\mbf{w}))$. Hence as long as $s$ is not on the wall in $\text{Pic}(M(\mbf{v},\mbf{w}))$, one can choose a small neighborhood $V_s$ of $s$ which does not intersect with any walls in $\text{Pic}(M(\mbf{v},\mbf{w}))$. In this way for arbitrary point $s'\in V_s$, $\mc{A}^{s'}_{\mc{L}}=\mc{A}^{s}_{\mc{L}}$. 
\end{proof}

In our settings, we always choose the slope points $s$ to be generic, i.e. $s$ does not lie on the wall. Note that this means that the $q$-difference operator $\mc{A}^{s}_{\mc{L}}$ is locally constant on $s$, thus we always fix $s$ in a generic position, i.e. $\mbf{B}_{s}(z)=1$ on $K_{T}(M(\mbf{v},\mbf{w}))$. We will see that the generic choice of the starting point $s$ would give us a convenient way to compute the quantum difference operator $\mbf{B}^s_{\mc{L}}(z)$.

Recall the wall set of a quiver variety $M(\mbf{v},\mbf{w})$ is defined as:
\begin{align}\label{wall-set-defn}
\text{Walls}(M(\mbf{v},\mbf{w}))=\{\mbf{m}\in\mbb{R}^{n}|R_{\mbf{m}}^+\text{ acts on } K_{T}(M(\mbf{v},\mbf{w}))\otimes K_{T}(M(\mbf{v},\mbf{w}))\text{ non trivially}\}.
\end{align}

\section{\textbf{Analysis of the quantum difference equations}}\label{sec:_textbf_analysis_of_the_quantum_difference_equations}
Fix the affine type $A$ quiver variety $M(\mbf{v},\mbf{w})$. In this section we will fix the slope $s$  of the quantum difference operator $\mbf{B}^s_{\mc{L}}(z)$ as defined in \eqref{defnofqdeoperator}. We assume throughout that $s$ is generic, i.e. $s$ will not lie on the wall of $M(\mbf{v},\mbf{w})$.

\subsection{Good representation for the quantum difference operators}

In \cite{Z23} we have shown that each monodromy operator $\mbf{B}_{\mbf{m}}(z)$ admits an explicit formula when expanding around $z=0$.
In this subsection we shall show that for any algebraic quantum difference operators $\mbf{B}_{\mc{L}}^s(z)$ of the affine type $A$ quiver varieties, it can be expressed as the ordered product of the monodromy operators $\mbf{B}_{\mbf{m}}(z)$ such that $\mbf{B}_{\mbf{m}}(z)$ can be represented via either $U_{q}(\mf{sl}_2)$-type or $U_{q}(\hat{\mf{gl}}_1)$-type, which is a closed formula in terms of the rational function over $z$.

\begin{thm}\label{good-representation}
Let $[s-\mc{L},s)$ be a generic path, with $s$ being generic. Then the quantum difference operator factorises as:
\begin{align*}
\mbf{B}_{\mc{L}}^s(z)=\mc{L}\prod^{\leftarrow}_{\mbf{m}\in[s-\mc{L},s)}\mbf{B}_{\mbf{m}}(z)
\end{align*}
and each $\mbf{B}_{\mbf{m}}(z)$ can be written either in one of the following two forms:
\begin{itemize}
	\item $U_{q}(\mf{sl}_2)$ type:
	\begin{equation}\label{goodpresentationwall-1}
	\begin{aligned}
	&\prod_{k=0}^{\substack{\rightarrow\\\infty}}\exp_{q^{2}}(-(q-q^{-1})z^{-k(\mbf{v}_{\gamma})}p^{-k\mbf{m}\cdot(\mbf{v}_{\gamma})}q^{k(\mbf{v}_{\gamma})^T((\frac{n^2r-1}{2})\bm{\theta}+\mbf{e}_1)-2-2\delta_{1\gamma}}f_{\gamma}e_{\gamma}'))\\
	=&\sum_{n=0}^{\infty}\frac{(q-q^{-1})^n}{[n]_{q^2}!}\frac{(-1)^n}{\prod_{\nu=1}^{n}(1-z^{-\mbf{v}_{\gamma}}p^{-\mbf{m}\cdot\mbf{v}_{\gamma}}q^{\nu(\mbf{v}_{\gamma})^T((\frac{n^2r-1}{2})\bm{\theta}+\mbf{e}_1)-2-2\delta_{1\gamma}})}f_{\gamma}^{n}e_{\gamma}'^n
	\end{aligned}
	\end{equation}
	with 
	\begin{align*}
	f_{\gamma}=P^{\mbf{m}}_{[i,j)},\qquad e_{\gamma}'=S_{\mbf{m}}(Q^{\mbf{m}}_{-[i,j)}).
	\end{align*}

	\item $U_{q}(\hat{\mf{gl}}_1)$ type.
	\begin{equation}\label{goodpresentationwall-2}
	\begin{aligned}
	\mu(\prod_{h=1}^{g}(\exp(-\sum_{k=1}^{\infty}\frac{n_kq^{-\frac{k\lvert\bm{\delta_h}\lvert}{2}}}{1-z^{-k\lvert\bm{\delta_h}\lvert}p^{k\mbf{m}\cdot\bm{\delta}_{h}}q^{-\frac{k\lvert\bm{\delta_h}\lvert}{2}}}\alpha^{\mbf{m},h}_{-k}\otimes\alpha^{\mbf{m},h}_{k})
	\end{aligned}
	\end{equation}
	with $\mu:\mc{B}_{\mbf{m}}\otimes\mc{B}_{\mbf{m}}\rightarrow\mc{B}_{\mbf{m}}$ the multiplication map, and
	\begin{align*}
	\alpha^{\mbf{m},h}_{-k}=P^{\mbf{m}}_{[h,h+k\delta_{\mbf{m}})},\qquad \alpha^{\mbf{m},h}_{k}=S_{\mbf{m}}(Q^{\mbf{m}}_{-[h,h+k\delta_{\mbf{m}})}).
	\end{align*}
\end{itemize}
\end{thm}

\textbf{Remark. }The difference of Theorem \ref{good-representation} and the monodromy operator formula given in \cite{Z23} is the following: The theorem says that one can find the generic path for the quantum difference operator $\mbf{B}^s_{\mc{L}}(z)$ such that with the fixed $M(\mbf{v},\mbf{w})$, the monodromy operator $\mbf{B}_{\mbf{m}}(z)$ can be written as a rational function over $z$ without doing the Taylor series expansion around $z=0$.

\begin{proof}
The construction of the path can be given as follows. Fix the quiver variety $M(\mbf{v},\mbf{w})$ and the slope subalgebra $\mc{B}_{\mbf{m}}$. Recall that the generators of $\mc{B}_{\mbf{m}}$ are given by $P^{\mbf{m}}_{\pm[i;j)_h}$, and $P^{\mbf{m}}_{\pm[i;i+1)_h}$ gives the generators of $\mc{B}_{\mbf{m}}$ of the minimal degree. All we need to do is to find a slope $\mbf{m}\in\mbb{R}^n$ such that $[i;i+1)_h$ exceeds the quiver dimension vector $\mbf{v}$ except for a single generator $P^{\mbf{m}}_{[j;j+1)_h}$ for some $j$.

To carry out the procedure, recall that $\mbf{m}\cdot[i;i+1)_h$ satisfying the above conditions can be written as:
\begin{align*}
(n-1)\mbf{m}\cdot\bm{\theta}+m_{i+1}+\cdots+m_{k}\in\mbb{Z}.
\end{align*}

Thus we define the following set of hyperplanes in $\mbb{R}^n$:
\begin{align*}
\mbf{v}\in\mbb{R}^n,\qquad[i;j)\cdot\mbf{v}\in\mbb{Z}.
\end{align*}

We require that for the integer vector $[i;j)=k_1e_1+k_2e_2+\cdots+k_ne_n$, we have $k_i\leq v_i$. Here $\mbf{v}=(v_1,\cdots,v_n)$ is the dimension vector for the quiver variety $M(\mbf{v},\mbf{w})$. The following lemma is easy to prove:
\begin{lem}
The monodromy operator $\mbf{B}_{\mbf{m}}(z)$ acts on $K_{T}(M(\mbf{v},\mbf{w}))$ nontrivially if $[i,j)\cdot\mbf{m}\in\mbb{Z}$ for some $[i,j)$ satisfies the condition above. Moreover, if $[i,j)\cdot\mbf{m}\in\mbb{Z}$, $P^{\mbf{m}}_{[i,j)}$ acts on $K_{T}(M(\mbf{v},\mbf{w}))$ non-trivially.
\end{lem}

For these hyperplanes $H_{[i,j),k}:=\{\mbf{m}\in\mbb{R}^n|\mbf{m}\cdot[i,j)+k=0\}\subset\mbb{R}^n$, choose $\mbf{m}\in\mbb{R}^n$ such that $\mbf{m}$ intersects with only one hyperplane $H_{[i,j),k}$. It is easy to see that the set of points $\mbf{m}\in\mbb{R}^n$ intersecting with only one or zero of the hyperplanes $H_{[i,j),k}$ is dense. Therefore for each hyperplane, one can choose a generic point $\mbf{m}$ for each hyperplane and connect them together. 

Then we determine the type of the monodromy operator $\mbf{B}_{\mbf{m}}(z)$ at each slope points $\mbf{m}$ with the corresponding intersecting hyperplane $H_{[i,j),k}$. If $[i,j)=n\bm{\theta}$, $\mc{B}_{\mbf{m}}\cong U_{q}(\hat{\mf{gl}}_1)^{\otimes n}$. In this case the monodromy operator $\mbf{B}_{\mbf{m}}$ can be written as:
\begin{align*}
\mbf{B}_{\mbf{m}}(z)=(\textbf{Heisenberg type})
\end{align*}
as in \eqref{goodpresentationwall-2}. Here $\delta_{\mbf{m}}$ is the minimal vector $[i,j)$ such that $\delta_{\mbf{m}}\cdot\mbf{m}\in\mbb{Z}$.

If $[i,j)\neq n\bm{\theta}$ for any $n\in\mbb{N}$, by \eqref{rootquantum}, $\mc{B}_{\mbf{m}}$ is a tensor product of the quantum affine algebras $U_{q}(\hat{\mf{gl}}_{l_h})$ such that there is at least one $l_{h}$ satisfying $l_h\geq2$. The corresponding element $P^{\mbf{m}}_{\pm[i,j)}$ and $\varphi_{[i,j)}$ gives a $U_q(\mf{sl}_2)$ subalgebra in $\mc{B}_{\mbf{m}}$. Therefore we have that in this case:
\begin{align*}
\mbf{B}_{\mbf{m}}(z)=(U_{q}(\mf{sl}_2)\textbf{ type})
\end{align*}
as in \eqref{goodpresentationwall-1}. Here $[i,j)$ is chosen such that $\mbf{m}\cdot[i,j)\in\mbb{Z}$. 
\end{proof}

\textbf{Remark}. The theorem implies that for the generic points $\mbf{m}\in\mbb{R}^n$ on the wall $[i,j)\cdot\mbf{m}\in\mbb{Z}$, the corresponding monodromy operator $\mbf{B}_{\mbf{m}}(z)$ is independent of the choice of the generic points $\mbf{m}$. This agrees with the monodromy operator $\mbf{B}_{w}(z)$ defined by Okounkov and Smirnov in \cite{OS22} where the monodromy operator $\mbf{B}_{w}(z)$ is defined only on the wall $w$ for the $K$-theoretic stable envelope. 

\subsection{General solution for the quantum difference equations}
In this section we use the good representation to give the fundamental solution of the quantum difference equation of the affine type $A$ quiver varieties $M(\mbf{v},\mbf{w})$.

For simplicity, we denote $X=M(\mbf{v},\mbf{w})$. It is well-known that $\text{Pic}(X)$ is generated by its tautological line bundles $\mc{L}_i$, and we denote $\mc{O}(1)=\mc{L}_1\cdots\mc{L}_{n}$. In this way the quantum difference equation at direction $\mc{O}(1)$ can be written as:
\begin{align*}
\Psi(p^{\mc{O}(1)}z)=\Psi(pz_1,\cdots,pz_n)=\mbf{B}_{\mc{O}(1)}^s(z)\Psi(z),\qquad z\in\text{Pic}(X)\otimes\mbb{C}^{*}.
\end{align*}

For simplicity we only concentrate on the quantum difference equation for $\mc{O}(1)$. For the quantum difference equation at other directions, i.e. $\mbf{B}_{\mc{L}}^s(z)$ with $\mc{L}$ other than $\mc{O}(1)$, the analysis is similar.

\subsubsection{Fundamental solution of $QDE$ near $z=0$}

First note that at $z=0\in\text{Pic}(X)\otimes\mbb{C}^{\times}$, i.e. the corresponding chamber in $\text{Pic}(X)\otimes\mbb{R}$ corresponds to the chamber such that all the parameters $a_1,\cdots,a_n$ go to $-\infty$ with $z_i:=e^{a_i}$. Let $\mbf{M}(0)=\mbf{B}_{\mc{L}_n}^s(0)\mbf{B}_{\mc{L}_{n-1}}^s(0)\cdots\mbf{B}_{\mc{L}_1}^s(0)=\mc{L}_{n}\otimes\cdots\otimes\mc{L}_{1}$ and this matrix is diagonal in $K_{T}(X)$ with respect to the fixed point basis. Let $P$ be the matrix with columns given by fixed point eigenvectors $[\lambda]$. We denote by $E_{0}$ the diagonal matrix of eigenvalues, so that:
\begin{align*}
\mbf{M}(0)P=PE_0
\end{align*}
and the eigenvalue diagonal matrix $E_{0}$ is $\mc{L}_1|_{\lambda}\cdots\mc{L}_n|_{\lambda}\in K_{T}(\text{pt})$.  There are two choices of the solution, the first is the elliptic type solutions: 
\begin{align*}
\Psi_0^{ell}(z)=P\Psi_{0}^{reg}(z)\prod_{i=1}^n(\Theta(E_{0}^{(i)},z_{i})).
\end{align*}

Here $E_{0}^{(i)}$ is the diagonal matrix with the value $\mc{L}_i|_{\lambda}$. $\Theta(E_{0}^{(i)},z_{i})$ is defined as:
\begin{align*}
\Theta(E_{0}^{(i)},z_{i})=\text{diag}(\frac{\theta(z_i,p)}{\theta(z_{i}\mc{L}_i^{-1}|_{\lambda},p)})_{\lambda\in M(\mbf{v},\mbf{w})^T}=\text{diag}(\frac{\theta(z_i,p)}{\theta(z_i\prod_{\square\in\lambda}\chi_{\square}^{-1},p)})_{\lambda\in M(\mbf{v},\mbf{w})^T}
\end{align*}
and this is the single-valued solution in $z$. For the definition of theta functions the reader can refer to Appendix \ref{sec:_textbf_appendix_i_basics_of_the_abelian_functions}.

The second is the multiplicative type:
\begin{align*}
\Psi_{0}(z)=P\Psi^{reg}_0(z)\exp(\sum_{i}\frac{\ln(E_0^{(i)})\ln(z_i)}{\ln(p)})
\end{align*}
which is the multi-valued solution in $z$ since $\exp(\sum_{i}\frac{\ln(E_0^{(i)})\ln(z_i)}{\ln(p)})$ is a multi-valued function. The regular part $\Psi_{0}^{reg}(z)$ satisfies the difference equation:
\begin{align}\label{regular-difference-equation}
\Psi_0^{reg}(zp^{\mc{L}})\mc{L}=\mbf{B}^s_{\mc{L}}(z)\Psi_0^{reg}(z).
\end{align}

\subsubsection{\textbf{Computation of $\Psi^{reg}_{0}(z)$ via} $\mc{O}(1)$}

For the second procedure, we can just write down the regular part of the solution as:
\begin{align}\label{formula-for-0}
\Psi_0^{reg}(z)=\prod^{\substack{\leftarrow\\\infty}}_{k=0}\mbf{B}^{s,*}(p^kz)^{-1}=\prod^{\substack{\leftarrow\\\infty}}_{k=0}\mbf{B}^{s,*}_{k}(z)^{-1}
\end{align}
with
\begin{align*}
\mbf{B}^{s,*}(z)=E_{0}^{-1}P^{-1}\mbf{B}_{\mc{O}(1)}^s(z)P,\qquad\mbf{B}^{s,*}_{k}(z):=E^{-k-1}_{0}P^{-1}\mbf{B}_{\mc{O}(1)}^s(p^kz)PE_{0}^k.
\end{align*}

Recall that by Theorem \ref{good-representation} the quantum difference operator $\mbf{B}_{\mc{O}(1)}^s(z)$ can be written as the product of monodromy operators $\mbf{B}_{\mbf{m}}(z)$ such that $\mbf{B}_{\mbf{m}}(z)$ are either of $U_{q}(\mf{sl}_2)$-type or $U_{q}(\hat{\mf{gl}}_1)$-type.

For example, if we take the path from $s$ to $s-\mc{O}(1)$ as the straight interval from $0$ to $-\mc{O}(1)$, the regular part of the fundamental solution around $z=0$ can also be written as:
\begin{align*}
\Psi^{reg}_0(z)=\prod^{\leftarrow}_{\mbf{m}\in(-\infty,s)}\mbf{B}_{\mbf{m}}^{*}(z)^{-1},\qquad\mbf{B}_{\mbf{m}}^*(z):=P^{-1}\mbf{B}_{\mbf{m}}(z)P
\end{align*}
Here $\mbf{m}\in(-\infty,0)$ means that we fix the initial point $s$ close to $0$ and the end point $s-n\mc{O}(1)$, and we require that $n\rightarrow\infty$.

\subsubsection{Fundamental solution near $z=\infty$}
For the solution around $z=\infty$, by the Lemma 3.2 in \cite{Z23}, the operator $\mbf{B}_{\mc{O}(1)}^s(\infty)$ consists of the product of the form $\mbf{m}((1\otimes S_{\mbf{m}})(R_{\mbf{m}}^{-})^{-1})$. It is a diagonalizable matrix over $\mbb{Q}((q,t))$ with eigenvalues given by the formal power series of $q$ and $t$, i.e. it has distinct eigenvalues for the generic value of $q$ and $t$. In fact, $\mbf{B}_{\mc{O}(1)}^s(z)$ has distinct eigenvalues for generic $z$. We will provide the proof of this fact in the Appendix \ref{appendix-2}. We denote by $H$ the matrix of eigenvectors and by $\mbf{E}_{\infty}$ the diagonal matrix of eigenvalues.
\begin{align*}
\mbf{B}_{\mc{O}(1)}^s(\infty)H=H\mbf{E}_{\infty}.
\end{align*}

Now we do the decomposition of $\mbf{E}_{\infty}$ in each direction $z_i$. Note that for each direction $z_i$ we have the difference equation:
\begin{align*}
\Psi(p^{\mc{L}_i}z)=\mbf{B}_{\mc{L}_i}^s(z)\Psi(z)
\end{align*}
So equivalently, the difference equation on $\mc{O}(1)=\mc{L}_1\cdots\mc{L}_n$ can be written as:
\begin{align*}
\Psi(p^{\mc{O}(1)}z)=\mbf{B}_{\mc{L}_1}^s(p^{\mc{O}(1)-\mc{L}_1}z)\mbf{B}_{\mc{L}_2}^s(p^{\mc{O}(1)-\mc{L}_1-\mc{L}_2}z)\cdots\mbf{B}_{\mc{L}_n}^s(z)\Psi(z).
\end{align*}

In this way we can see that:
\begin{align*}
\mbf{B}_{\mc{O}(1)}^s(\infty)=\mbf{B}_{\mc{L}_1}^s(\infty)\cdots\mbf{B}_{\mc{L}_n}^s(\infty).
\end{align*}

Furthermore $\mbf{B}_{\mc{L}_i}^s(\infty)\mbf{B}_{\mc{L}_j}^s(\infty)=\mbf{B}_{\mc{L}_i\mc{L}_j}^s(\infty)=\mbf{B}_{\mc{L}_j}^s(\infty)\mbf{B}_{\mc{L}_i}^s(\infty)$, and this means that the operators $\mbf{B}_{\mc{L}_i}^s(\infty)$ share the same eigenvectors $H$ with different eigenvalues $E_{\infty}^{(i)}$.

In this way we can first construct the elliptic solution as:
\begin{align*}
\Psi_{\infty}^{ell}=H\Psi_{\infty}^{reg}(z)\prod_{i}\Theta(E_{\infty}^{(i)},z_i).
\end{align*}

the solution around $z=\infty$ can be written as:
\begin{align*}
\Psi_{\infty}(z)=H\Psi_{\infty}^{reg}(z)\exp(\sum_{i}\frac{\ln(E_{\infty}^{(i)})\ln(z_i)}{\ln(p)})
\end{align*}
such that
\begin{align*}
\Psi^{reg}_{\infty}(z)\mbf{E}_{\infty}=H^{-1}\mbf{B}_{\mc{O}(1)}^s(z)H\Psi^{reg}_{\infty}(z).
\end{align*}

We use the similar procedure to compute the regular part of $\Psi^{reg}_{\infty}(z)$, and it is easy to see that the solution is of the form:
\begin{align*}
\Psi^{reg}_{\infty}(z)=\prod^{\substack{\rightarrow\\\infty}}_{k=0}\mbf{B}^{s,\heartsuit}_{-k}(z)
\end{align*}
with
\begin{align*}
\mbf{B}^{s,\heartsuit}_{-k}(z)=\mbf{E}_{\infty}^{k-1}\mbf{B}^{s,*}
(zp^{-k})\mbf{E}_{\infty}^{-k},\qquad\mbf{B}^{s,\heartsuit}_{\mc{O}(1)}(z):=H^{-1}\mbf{B}^s_{\mc{O}(1)}(z)H.
\end{align*}

Also we denote 
\begin{align*}
\mbf{B}_{\mbf{m}}^{\heartsuit}(z):=H^{-1}\mbf{B}_{\mbf{m}}(z)H,\qquad\mbf{B}_{\mbf{m}}^{\heartsuit}=H^{-1}\mbf{B}_{\mbf{m}}H.
\end{align*}

Now the fundamental solution around $z=\infty$ can be written as:
\begin{align*}
\Psi^{reg}_{\infty}(z)=\lim_{N\rightarrow\infty}\Psi^{reg}_{N}(z)=\lim_{N\rightarrow\infty}(\prod^{\rightarrow}_{\mbf{m}\in[0,N)}\mbf{B}_{\mbf{m}}^{\heartsuit}(z))(\mbf{B}_{\mc{O}(1)}^{s,\heartsuit}(0))^N\mbf{E}_{\infty}^{-N}.
\end{align*}

\subsection{Regular fundamental solution of different slopes}

Given two distinct generic slopes $s,s'\in\mbb{R}^n$. We can connect the regular parts of the fundamental solutions $\Psi_{0}^{reg,s}(z)$, $\Psi^{reg,s'}_{0}(z)$ at slopes $s,s'$ by the monodromy operators.

\begin{lem}
Suppose that $\Psi_{0}^{reg,s}(z)$, $\Psi^{reg,s'}_{0}(z)$ are the regular fundamental solutions at $z=0$ for the quantum difference operators $\mc{A}^{s}_{\mc{L}}$, $\mc{A}^{s'}_{\mc{L}}$. Then,
\begin{align*}
\Psi^{reg,s'}_{0}(z)=\prod_{\mbf{m}\in(s',s)}^{\leftarrow}\mbf{B}_{\mbf{m}}(z)\Psi^{reg,s}_{0}(z)
\end{align*}
and here $(s,s')$ means that we choose a path from $s'$ to $s$ such that $\mbf{m}\in(s',s)$ are the slope points in the path, which is a finite product.
\end{lem}
\begin{proof}
The proof is a direct computation:
\begin{equation*}
\begin{aligned}
\prod_{\mbf{m}\in(s',s)}^{\leftarrow}\mbf{B}_{\mbf{m}}(zp^{\mc{L}})\Psi^{reg,s}_{0}(zp^{\mc{L}})\mc{L}=&\prod_{\mbf{m}\in(s,s')}^{\leftarrow}\mbf{B}_{\mbf{m}}(zp^{\mc{L}})\mbf{B}_{\mc{L}}^{s}(z)\Psi^{reg,s}_{0}(z)\\
=&\prod_{\mbf{m}\in(s,s')}^{\leftarrow}\mbf{B}_{\mbf{m}}(zp^{\mc{L}})\mc{L}\prod_{\mbf{m}\in(s-\mc{L},s)}^{\leftarrow}\mbf{B}_{\mbf{m}}(z)\Psi^{reg,s}_{0}(z)\\
=&\mc{L}\prod_{\mbf{m}\in(s'-\mc{L},s-\mc{L})}^{\leftarrow}\mbf{B}_{\mbf{m}}(z)\prod_{\mbf{m}\in(s-\mc{L},s)}^{\leftarrow}\mbf{B}_{\mbf{m}}(z)\Psi^{reg,s}_{0}(z)\\
=&\mc{L}\prod_{\mbf{m}\in(s'-\mc{L},s)}^{\leftarrow}\mbf{B}_{\mbf{m}}(z)\Psi^{reg,s}_{0}(z)\\
=&\mc{L}\prod_{\mbf{m}\in(s'-\mc{L},s')}^{\leftarrow}\mbf{B}_{\mbf{m}}(z)\prod_{\mbf{m}\in(s',s)}^{\leftarrow}\mbf{B}_{\mbf{m}}(z)\Psi^{reg,s}_{0}(z)
\end{aligned}
\end{equation*}
where the first equality comes from the difference equation \eqref{regular-difference-equation}.

By the uniqueness of the solution, the proof is finished.
\end{proof}

\subsection{Connection matrix for the two solutions}
The transition matrix between two fundamental solutions $\Psi_{0}(z)$ and $\Psi_{\infty}(z)$ is defined as:
\begin{align*}
\textbf{Mon}(z):=\Psi_{0}(z)^{-1}\Psi_{\infty}(z)
\end{align*}

It is clear that $\textbf{Mon}(p^{\mc{L}}z)=\textbf{Mon}(z)$, and we expect that $\textbf{Mon}(z)$ is an abelian function over the abelian variety $\text{Pic}(M(\mbf{v},\mbf{w}))\otimes\mbb{C}^{\times}/(z\sim zp^{\mc{L}})$. However, note that $\textbf{Mon}(z)$ might be a multivalued function over $z$, which means that $\textbf{Mon}(ze^{2\pi i})\neq\textbf{Mon}(z)$.

Also for the sake of the regular part of the fundamental solution, we can also define the regular part of the transition matrix as:
\begin{align*}
\textbf{Mon}^{reg}(z):=\exp(\sum_{i}\frac{\ln(E_{0}^{(i)})\ln(z_i)}{\ln(p)})\textbf{Mon}(z)\exp(-\sum_{i}\frac{\ln(E_{\infty}^{(i)})\ln(z_i)}{\ln(p)}).
\end{align*}

It is easy to check that:
\begin{align*}
\textbf{Mon}^{reg}(pz)=E_{0}\textbf{Mon}^{reg}(z)E_{\infty}^{-1}.
\end{align*}

We also define the connection matrix for the elliptic type solution $\Psi_{0}^{ell}(z)$, $\Psi_{\infty}^{ell}(z)$ as:
\begin{align*}
\textbf{Mon}^{ell}(z):=\Psi_{0}^{ell}(z)^{-1}\Psi_{\infty}^{ell}(z)
\end{align*}
which also satisfies $\textbf{Mon}^{ell}(zp^{\mc{L}})=\textbf{Mon}^{ell}(z)$. Moreover, $\textbf{Mon}^{ell}(z)$ is a single-valued meromorphic function over $z$. Therefore $\mbf{Mon}^{ell}(z)$ lives in $GL(M(\mbf{E}_{p}^r))$ the space of invertible matrices valued in the meromorphic function over the abelian variety $M(\mbf{E}_{p}^r)$. Using the fact that for the variety $\mbf{E}_{p}^r$, the meromorphic function over $\mbf{E}_{p}^r$ is generated by the linear combinations of the ratio of the theta function of the form (For the proof, see page 91 in \cite{Sie89}:
\begin{align}\label{theta-function-form-of-ell-connection}
\prod_{i_1}\frac{\theta(z_1a_{i_1})}{\theta(z_1b_{i_1})}\cdots\prod_{i_r}\frac{\theta(z_ra_{i_r})}{\theta(z_{r}b_{i_r})}\times\cdots\times\prod_{\mbf{k}}\frac{\theta(z_1^{k_1}z_2^{k_2}\cdots z_{r}^{k_r}a_{\mbf{k}})}{\theta(z_1^{k_1}z_2^{k_2}\cdots z_{r}^{k_r}b_{\mbf{k}})}
\end{align}
such that $\prod_{i_1}a_{i_1}=\prod_{i_1}b_{i_1}$,$\cdots$, $\prod_{i_r}a_{i_r}=\prod_{i_r}b_{i_r}$, $\prod_{\mbf{k}}a_{\mbf{k}}=\prod_{\mbf{k}}b_{\mbf{k}}$. Now via the relation between $\mbf{Mon}^{ell}(z)$ and $\mbf{Mon}(z)$ we have that:
\begin{align}\label{elliptic-to-multivalued}
\exp(-\sum_{i}\frac{\ln(E_{0}^{(i)})\ln(z_i)}{\ln(p)})\textbf{Mon}(z)\exp(\sum_{i}\frac{\ln(E_{\infty}^{(i)})\ln(z_i)}{\ln(p)})=\prod_{i}(\Theta(E_{0}^{(i)},z_i))^{-1}\textbf{Mon}^{ell}(z)\prod_{i}(\Theta(E_{\infty}^{(i)},z_i)).
\end{align}

By construction, we can see that $\textbf{Mon}^{reg}(z)$ can be written as the linear combination of the theta function as the following form:
\begin{align}\label{theta-function-form-of-reg-conn}
\prod_{i=1}^r(\Theta(E_{0}^{(i)},z_{i}))^{-1}\prod_{i_1}\frac{\theta(z_1a_{i_1})}{\theta(z_1b_{i_1})}\cdots\prod_{i_r}\frac{\theta(z_ra_{i_r})}{\theta(z_{r}b_{i_r})}\times\prod_{\mbf{k}}\frac{\theta(z_1^{k_1}z_2^{k_2}\cdots z_{r}^{k_r}a_{\mbf{k}})}{\theta(z_1^{k_1}z_2^{k_2}\cdots z_{r}^{k_r}b_{\mbf{k}})}\prod_{i=1}^r(\Theta(E_{\infty}^{(i)},z_{i})).
\end{align}

\subsection{$p\rightarrow0$ limit}
In this subsection we turn to analyze the limit $p\rightarrow0$ for the function of the type $f(p^sz)$ with $s\in\mbb{R}^{r}$.

Now we choose our path such that $\mbf{B}_{\mbf{m}}(z)$ are all of the generic point $\mbf{m}$ in the sense of Theorem \ref{good-representation}. It remains to analyze $\mbf{B}_{\mbf{m}}(p^sz)$, by the formula in Theorem \ref{good-representation}, one can directly compute that:
\begin{align}\label{limit-of-monodromy-operator}
\lim_{p\rightarrow0}\mbf{B}_{\mbf{m}}(p^sz)=
\begin{cases}
\mu((1\otimes S_{\mbf{m}})(R_{\mbf{m}}^{+})_{21}^{-1})&s<\mbf{m},\mbf{m}\in U_{s}\\
\mbf{B}_{\mbf{m}}(z)|_{p=1}&s=\mbf{m}\\
\mbf{B}_{\mbf{m}}(z_{I},\theta_{J})&s_{I}=\mbf{m}_{I},s_{I^c}\neq \mbf{m}_{I^c}\\
1&\text{Otherwise}
\end{cases}.
\end{align} 

Here $\mu:\mc{B}_{\mbf{m}}\otimes\mc{B}_{\mbf{m}}\rightarrow\mc{B}_{\mbf{m}}$ is the multiplication map. Here we use the notation $s<\mbf{m}$ to mean that for $s=(s_1,\cdots,s_n),\mbf{m}=(m_1,\cdots,m_n)$ $s<\mbf{m}$ if there is at least one $s_i,m_i$ such that $s_i<m_i$. The symbol $\mbf{B}_{\mbf{m}}(z_{I},\theta_{J})$ is defined as  follows: We first write down $\mbf{B}_{\mbf{m}}(z)=\mbf{B}_{\mbf{m}}(z_1,\cdots,z_n)$, and we split the set $\{1,\cdots,n\}$ into disjoint union of two subsets $I$ and $J$. Now $z_{I}$ is defined as the variables $(z_i)_{i\in I}$ with index $i$ in the set $I$. $\theta_{J}=(\theta_j)_{j\in J}$ such that $\theta_{j}=0$ if $s_{j}<m_{j}$ and $\infty$ if $s_{j}>m_{j}$. $U_{s}$ is a suitable small neighborhood of $s$.

Denote $\mbf{B}_{\mbf{m}}:=\mu((1\otimes S_{\mbf{m}})(R_{\mbf{m}}^{+})_{21}^{-1})$. Now we choose a good representation as in Theorem \ref{good-representation} of the quantum difference operator $\mbf{B}_{\mc{O}(1)}^s(z)=\mc{O}(1)\prod^{\leftarrow}_{\mbf{m}\in(s,s-\mc{O}(1)]}\mbf{B}_{\mbf{m}}(z)$ such that each monodromy operator $\mbf{B}_{\mbf{m}}(z)$ is of either the $U_{q}(\mf{sl}_2)$-type or of the $U_{q}(\hat{\mf{gl}}_1)$-type. The following theorem gives the expression for the $p\rightarrow0$ limit of the connection matrix:
\begin{thm}\label{p0-limit-connection}
For generic $s\in\mbb{Q}^n$ the connection matrix has the following asymptotic as $p\rightarrow0$:
\begin{align*}
\lim_{p\rightarrow0}\textbf{Mon}^{reg}(p^sz)=
\begin{cases}
\prod^{\leftarrow}_{0\leq\mbf{m}<s}(\mbf{B}_{\mbf{m}}^*)^{-1}\cdot\mbf{T},\qquad s\geq0\\
\prod_{s<\mbf{m}<0}^{\leftarrow}\mbf{B}_{\mbf{m}}^*\cdot\mbf{T},\qquad s<0
\end{cases}.
\end{align*}

Here $\mbf{T}:=P^{-1}H$, $0\leq\mbf{m}<s$ means the slope points $\mbf{m}$ in one generic path from $0$ to the point $s$ without intersecting $s$. Here $s\geq0$ and $s<0$ stands for $s=(s_1,\cdots,s_n)$ such that $s_{i}\leq 0$ or $s_{i}>0$ for every component $s_i$. Here $\mbf{B}_{\mbf{m}}^*$ is defined as $P^{-1}\mbf{B}_{\mbf{m}}P$.
\end{thm}

Roughly speaking, the theorem tells us that the nodal limit of the regular part of the connection matrix $\textbf{Mon}^{reg}(z)$ can be described as the ordered product of the monodromy operators. The choice of the monodromy operators on the path would not affect the expression of $\textbf{Mon}^{reg}(z)$.

For the generic $s$, we can in principle write down the formula for $\lim_{p\rightarrow0}\textbf{Mon}^{reg}(p^sz)$ in terms of the monodromy operators $\mbf{B}_{\mbf{m}}^*$, but in that case the resulting expression would be complicated and chaotic. We will see in the proof that there is an explicit formula for $\lim_{p\rightarrow0}\textbf{Mon}^{reg}(p^sz)$.

To prove the theorem, we first analyze the $p\rightarrow0$ limit of the fundamental solution $\psi^{reg}_{0}=\prod^{\substack{\leftarrow\\\infty}}_{k=0}\mbf{B}_{k}^{s,*}(z)^{-1}$ around $z=0$. We denote:
\begin{align*}
C(s)=\{\mbf{x}\in\mbb{R}^n|x_{i}\leq s_i\}\subset\mbb{R}^n
\end{align*}
and the following set: 
\begin{align*}
B_{k}:=\{\mbf{s}\in\mbb{R}^n|\exists w\in\text{Wall}_{0}\text{ such that }\exists I\subset\{1,\cdots,r\}, s_{I}+k\leq w_{I}\}\subset\mbb{R}^n.
\end{align*}

For convenience we consider the following stratification of $B_{k}$:
\begin{align*}
B_{k}=B_{k}^o\cup\partial B_{k}
\end{align*}
with $B_{k}^o$ denoting the interior points of $B_{k}$, and $\partial B_{k}$ denotes the points on the boundary of $B_{k}$.

It is easy to see that $B_{k}\subset B_{0}$ and moreover $B_{k+1}\subset B_{k}$. Using this fact we can compute the nodal limit of $\psi^{reg}_{0}(p^sz)$:
\begin{prop}\label{p0-limit-z-0}
The asymptotic behavior of $\psi^{reg}_{0}(z)$ for $r\in B_{l}-B_{l+1}$ with $l\geq0$ goes in the following way:
\begin{align*}
\lim_{p\rightarrow0}\psi^{reg}_{0}(p^rz)=\begin{cases}
\prod_{k=0}^{\substack{\leftarrow\\l-1}}(\mbf{B}_{k}^{s,*})^{-1}\cdot\mbf{B}_{l,asymp}^{s,*}(r)^{-1}&r\in B_{l}^o\\
\prod_{k=0}^{\substack{\leftarrow\\l-1}}(\mbf{B}_{k}^{s,*})^{-1}\cdot\mbf{B}_{l,asymp}^{s,*}(z,r)^{-1}&r\in\partial B_{l}
\end{cases}
\end{align*}
where 
\begin{equation*}
\begin{aligned}
&\mbf{B}_{l,asymp}^{s,*}(z,r):=E_{0}^{-l-1}P^{-1}(\mc{O}(1)\prod_{\substack{w\in\text{Walls}_{0}\cap C(r)^{c}\\w<r,r_{\sigma(1)}}}\mbf{B}_{w}\cdot \mbf{B}_{w_1,r_{\sigma(1)}}(z_{I},\theta_{J})|_{p=1}\cdot\prod_{\substack{w\in\text{Walls}_{0}\cap C(r)^{c}\\w<r,r_{\sigma(2)}}}\mbf{B}_{w}\cdot \mbf{B}_{w_1,r_{\sigma(2)}}(z_{I},\theta_{J})|_{p=1}\cdot\\
&\cdots\prod_{\substack{w\in\text{Walls}_{0}\cap C(r)^{c}\\w<r,r_{\sigma(l)}}}\mbf{B}_{w}\cdot \mbf{B}_{w_1,r_{\sigma(l)}}(z_{I},\theta_{J})|_{p=1}\cdot\prod_{\substack{w\in\text{Walls}_{0}\cap C(r)^{c}\\w<r,r_{\sigma(l+1)}}}\mbf{B}_{w})PE_{0}^k\\
&\mbf{B}_{l,asymp}^{s,*}(r):=E_{0}^{-l-1}P^{-1}(\mc{O}(1)\prod_{w\in\text{Walls}_{l}\cap C(r)^{c}}\mbf{B}_{w})PE_{0}^k.
\end{aligned}
\end{equation*}

Here $\text{Walls}_{l}$ stands for the set of wall slope points in a small neighborhood of $[l+s_i,l+1+s_i]_{1\leq i\leq n}$. $(w_1,r_{\sigma(l)})$ stands for the points $(w_1,\cdots,r_{\sigma(l)},\cdots,w_n)$ and here $r_{\sigma(i)}$ is ordered such that $r_{\sigma(1)}\leq r_{\sigma(2)}\leq\cdots\leq r_{\sigma(n)}$.
\end{prop}
\begin{proof}
This follows from the formula \eqref{formula-for-0} of the fundamental solution $\psi^{reg}_{0}(z)$. The term $\mbf{B}_{l,asymp}^{s,*}(r)$ and $\mbf{B}_{l,asymp}^{s,*}(z,r)$ comes from computing the precise limit of $\mbf{B}_{k}^{s,*}(zp^r)$ using the formula \eqref{limit-of-monodromy-operator}. We leave the computation of the limit as the exercise.
\end{proof}

We use the similar procedure to compute the regular part of $\Psi^{reg}_{\infty}(z)$, and it is easy to see that the solution is of the form:
\begin{align*}
\Psi^{reg}_{\infty}(z)=\prod^{\substack{\rightarrow\\\infty}}_{k=0}\mbf{B}^{s,\heartsuit}_{-k}(z)
\end{align*}
with
\begin{align*}
\mbf{M}^{\heartsuit}_{-k}(z)=E_{\infty}^{k-1}\mbf{B}^{s,\heartsuit}
(zp^{-k})\mbf{E}_{\infty}^{-k}.
\end{align*}

Also we denote 
\begin{align*}
\mbf{B}_{\mbf{m}}^{\heartsuit}(z):=H^{-1}\mbf{B}_{\mbf{m}}(z)H,\qquad\mbf{B}_{\mbf{m}}^{\heartsuit}=H^{-1}\mbf{B}_{\mbf{m}}H.
\end{align*}

For the general case, we set
\begin{align*}
D(s)=\{\mbf{x}\in\mbb{R}^n|x_{i}\geq s_i\}.
\end{align*}

Define $B_{l}$ for $l\leq0$ in a similar way:
\begin{align*}
B_{-k}:=\{\mbf{s}\in\mbb{R}^n|\exists \mbf{m}\in\text{Wall}_{0}\text{ such that }\exists I\subset\{1,\cdots,n\}, s_{I}-k\geq m_{I}\}\subset\mbb{R}^n.
\end{align*}

Under this situation, doing the similar calculation above, we have that:
\begin{prop}\label{p0-limit-z-infty}
The asymptotic behavior of $\psi^{reg}_{\infty}(z)$ for $r\in B_{-l}-B_{-l-1}$ with $l\geq0$ goes in the following way:
\begin{align*}
\lim_{p\rightarrow0}\psi^{reg}_{\infty}(p^rz)=\begin{cases}
\mbf{B}_{-l,asymp}^{\heartsuit}(r)\cdot\prod_{k=-l+1}^{\substack{\leftarrow\\-1}}(\mbf{B}_{k}^{s,\heartsuit})^{-1}&r\in B_{-l}^o\\
\mbf{B}_{-l,asymp}^{\heartsuit}(z,r)\cdot\prod_{k=-l+1}^{\substack{\leftarrow\\-1}}(\mbf{B}_{k}^{s,\heartsuit})^{-1}&r\in\partial B_{-l}.
\end{cases}
\end{align*}
where 
\begin{equation*}
\begin{aligned}
&\mbf{B}_{l,asymp}^{s,\heartsuit}(z,r):=E_{\infty}^{l-1}H^{-1}(\mc{O}(1)\prod_{\substack{w\in\text{Walls}_{0}\cap D(r)^{c}\\w<r,r_{\sigma(1)}}}\mbf{B}_{w}\cdot \mbf{B}_{w_1,r_{\sigma(1)}}(z_{I},\theta_{J})|_{p=1}\cdot\prod_{\substack{w\in\text{Walls}_{0}\cap D(r)^{c}\\w<r,r_{\sigma(2)}}}\mbf{B}_{w}\cdot\\ 
&\mbf{B}_{w_1,r_{\sigma(2)}}(z_{I},\theta_{J})|_{p=1}\cdot
\cdots\prod_{\substack{w\in\text{Walls}_{0}\cap D(r)^{c}\\w<r,r_{\sigma(l)}}}\mbf{B}_{w}\cdot \mbf{B}_{w_1,r_{\sigma(l)}}(z_{I},\theta_{J})|_{p=1}\cdot\prod_{\substack{w\in\text{Walls}_{0}\cap D(r)^{c}\\w<r,r_{\sigma(l+1)}}}\mbf{B}_{w})HE_{\infty}^{-l}\\
&\mbf{B}_{l,asymp}^{s,\heartsuit}(r):=E_{\infty}^{l-1}H^{-1}(\mc{O}(1)\prod_{w\in\text{Walls}_{l}\cap D(r)^{c}}\mbf{B}_{w})HE_{\infty}^{-l}
\end{aligned}
\end{equation*}
and here $\text{Walls}_{l}$ stands for the set of walls in a small neighborhood of $[l+s_i,l+1+s_i]_{1\leq i\leq n}$. $(w_1,r_{\sigma(l)})$ stands for the points $(w_1,\cdots,r_{\sigma(l)},\cdots,w_n)$ and here $r_{\sigma(i)}$ is ordered such that $r_{\sigma(1)}\leq r_{\sigma(2)}\leq\cdots\leq r_{\sigma(n)}$.

\end{prop}

\begin{proof}
The proof is similar to proposition $6$ in \cite{S21}. We consider a finite approximation of infinite product
\begin{align*}
\Psi_{N}^{reg}(z)=\mbf{B}_{-1}^{s,\heartsuit}(z)\mbf{B}_{-2}^{s,\heartsuit}(z)\cdots\mbf{B}_{-N}^{s,\heartsuit}(z)=\mbf{B}^{s,\heartsuit}(\frac{z}{p})\mbf{B}^{s,\heartsuit}(\frac{z}{p^2})\cdots\mbf{B}^{s,\heartsuit}(\frac{z}{p^N})E_{\infty}^{-N}
\end{align*}
via the translation formula, we can write down:
\begin{align*}
\Psi^{reg}_{N}(z)=(\prod^{\rightarrow}_{[0,N)\bm{\theta}\in\text{Walls}_0+\mbb{N}}\mbf{B}^{\heartsuit}_{w}(z))(\mbf{B}_{\mc{O}(1)}^{s,\heartsuit}(0))^NE_{\infty}^{-N}.
\end{align*}

Now for $s\in B_{-l}^o$, we have that:
\begin{align*}
\lim_{p\rightarrow0}\Psi_{N}^{reg}(zp^s)=(\prod^{\rightarrow}_{(\text{Walls}_0+\mbb{N})\cap[0,N)\cap C(s)^c}\mbf{B}^{\heartsuit}_{w}(z))(\mbf{B}_{\mc{O}(1)}^{s,\heartsuit}(0))^NE_{\infty}^{-N}.
\end{align*}

Also note that by taking $z=\infty$, we have that:
\begin{align*}
E_{\infty}^N=(\prod^{\rightarrow}_{[0,N)\bm{\theta}\in\text{Walls}_0+\mbb{N}}\mbf{B}^{\heartsuit}_{w}(z))(\mbf{B}_{\mc{O}(1)}^{s,\heartsuit}(0))^N.
\end{align*}
which means that 
\begin{align*}
(\prod^{\rightarrow}_{[0,N\bm{\theta})\in\text{Walls}_0+\mbb{N}}\mbf{B}^{\heartsuit}_{w}(\mbf{B}_{\mc{O}(1)}^{s,\heartsuit}(0))^NE_{\infty}^{-N}=1.
\end{align*}

Thus comparing the result we have the statement for $s\in B_{l}^o$. For the statement that $s\in\partial B_{l}$, just mimic the computation as in the $\psi_{0}^{reg}(z)$, which is complicated but straightforward computation. 

\end{proof}

Though the formula above is really complicated, here we only consider the nodal limit when $s$ is the generic point. Using the proposition we can have the following result for the asymptotics of $\Psi^{reg}_{\infty}(zp^s)$ and $\Psi^{reg}_{0}(zp^s)$ as $p\rightarrow0$:
\begin{prop}
The solutions $\Psi^{reg}_{\infty}(z)$ and $\Psi^{reg}_{0}(zp^s)$ has the following limit for generic $s$:
\begin{align}\label{final-limit-simple-01}
\lim_{p\rightarrow0}\Psi^{reg}_{\infty}(zp^s)=
\begin{cases}
\prod^{\leftarrow}_{\mbf{m}\in[0,s)}(\mbf{B}^{\heartsuit}_{\mbf{m}})^{-1}&s\geq0\\
1&s<0
\end{cases},\qquad B_{\mbf{m}}^{\heartsuit}:=H^{-1}\mbf{B}_{\mbf{m}}H
\end{align}

\begin{align}\label{final-limit-simple-02}
\lim_{p\rightarrow0}\Psi^{reg}_{0}(zp^s)=
\begin{cases}
\prod^{\leftarrow}_{\mbf{m}\in(s,0]}(\mbf{B}^{\heartsuit}_{\mbf{m}})&s\leq0\\
1&s\geq0
\end{cases},\qquad B_{\mbf{m}}^*:=P^{-1}\mbf{B}_{\mbf{m}}P.
\end{align}

Here $[0,s)$ ($(s,0]$) stands for a path from $0$ ($s$) to $s$ ($0$) containing 0. Here $s\geq0$ and $s<0$ stands for $s=(s_1,\cdots,s_n)$ such that $s_{i}\leq 0$ or $s_{i}>0$ for every component $s_i$.
\end{prop}

Now combining the formula \eqref{final-limit-simple-01} and \eqref{final-limit-simple-02}, we have obtained the Theorem \ref{p0-limit-connection}.

\section{\textbf{Degeneration to the Dubrovin connection}}\label{sec:_textbf_degeneration_to_the_dubrovin_connection}
In this section we shall construct the degeneration limit of the quantum difference equations.
\subsection{Dubrovin connection for quiver varieties}

Fix a quiver variety $M_{\theta}(\mbf{v},\mbf{w})$. It is proved in \cite{MO12} that the Dubrovin connection for the quantum equivariant cohomology $H_{T}^*(M_{\theta}(\mbf{v},\mbf{w}))$ is given as:
\begin{align*}
\nabla_{\lambda}=d_{\lambda}-Q(\lambda),\qquad\lambda\in H^2(M_{\theta}(\mbf{v},\mbf{w}))\cong\text{Pic}(M_{\theta}(\mbf{v},\mbf{w}))
\end{align*}
Here:
\begin{align*}
Q(\lambda)=c_{1}(\lambda)+\hbar\sum_{\theta\cdot\alpha>0}\frac{\alpha(\lambda)}{1-z^{-\alpha}}e_{\alpha}e_{-\alpha}+\text{constant opeartor}.
\end{align*}

Here $\{e_{\alpha}\}$ are the generators of the Maulik-Okounkov Lie algebra $\mf{g}_{Q}^{MO}\subset Y_{\hbar_1,\hbar_2}(\mf{g}_{Q}^{MO})$ defined in \cite{MO12}. It has been proved \cite{BD23} that the positive half of the Maulik-Okounkov Lie algebra $\mf{n}_{Q}^{MO}$ is isomorphic to the BPS Lie algebra $\mf{n}_{Q}\subset Y_{\hbar_1,\hbar_2}(\mf{g}_{Q})$. It was also shown in \cite{Dav25} that $\mf{n}_{Q}$ is isomorphic to the positive half of $\hat{\mf{gl}}_n$. Therefore we have the isomorphism of Lie algebras
\begin{align}\label{isomorphism-of-Lie-algebras}
\mf{g}_{Q}^{MO}\cong\hat{\mf{gl}}_n,\qquad e_{\pm\alpha}\mapsto E_{\pm[i,j)}
\end{align}
and here the generators $e_{\pm\alpha}$ in $\mf{g}^{MO}_{Q}$ are sent to $E_{\pm[i,j)}$ in $\hat{\mf{gl}}_n$ if $\alpha=\alpha_i+\alpha_{i+1}+\cdots+\alpha_{j-1}$, and here the index $i$ means that $\alpha_i=\alpha_{i\text{ mod }n}$.

There are many interesting properties for the Dubrovin connection. The Dubrovin connection has regular singularities $z^{-\alpha}=1$ and $z=0,\infty$, and the corresponding monodromy representation:
\begin{align*}
\pi_{1}((\mbb{C}^*)^n-\text{Sing})\rightarrow\text{Aut}(H_{T}^*(M_{\theta}(\mbf{v},\mbf{w})))
\end{align*}
can be represented by the quantum Weyl group for the corresponding quantum group $U_{q}(\mf{g}_{Q})$.

In our cases, for the affine type $A$ quiver varieties $M(\mbf{v},\mbf{w})$ for $\bm{\theta}=(1,\cdots,1)$, the corresponding BPS Lie algebra $\mf{g}_{Q}$ is isomorphic to $\hat{\mf{gl}}_n$. The corresponding Dubrovin connection can be written as:
\begin{equation*}
\begin{aligned}
Q(\lambda)=&c_1(\lambda)+\hbar\sum_{\alpha>0}\frac{\alpha(\lambda)}{1-z^{-\alpha}}e_{\alpha}e_{-\alpha}+\cdots\\
=&c_1(\lambda)+\hbar\sum_{\alpha\in\text{real positive roots}+\{0\}}\sum_{k\geq0}\frac{(\alpha+k\delta)\lambda}{1-z^{-\alpha-k\delta}}e_{\alpha+k\delta}e_{-\alpha-k\delta}+\cdots
\end{aligned}
\end{equation*}
and here the dots stand for constant operators. If we take $\lambda=\mc{O}(1)$, we have that:
\begin{align*}
Q(\mc{O}(1))=c_{1}(\mc{O}(1))+\hbar\sum_{\alpha\in\text{real positive roots}+\{0\}}\sum_{k\geq0}\frac{kn+|\alpha|}{1-z^{-\alpha-k\delta}}e_{\alpha+k\delta}e_{-\alpha-k\delta}+\cdots
\end{align*}
with $e_{\pm\alpha}$ defined as $E_{\pm[i,j)}$ as mentioned above. In this way we can write down the quantum differential equation as:
\begin{align}\label{good-form-for-Dubrovin-connection}
Q(\mc{L})=c_1(\mc{L})+\hbar\sum_{[i,j)}\frac{\mc{L}\cdot[i,j)}{1-z^{-[i,j)}}E_{[i,j)}E_{-[i,j)}+\cdots.
\end{align}

\subsection{Degeneration limit of the slope subalgebra}

In this subsection we describe the degeneration limit of the slope subalgebra $\mc{B}_{\mbf{m}}\subset U_{q,t}(\hat{\hat{\mf{sl}}}_n)$ into the subalgebra of $Y_{\hbar_1,\hbar_2}(\hat{\mf{sl}}_n)$.

The description we shall use is as follows: Using the algebra embedding:
\begin{align*}
U_{q,t}(\hat{\hat{\mf{sl}}}_n)\hookrightarrow\prod_{\mbf{w}}\text{End}(K(\mbf{w})).
\end{align*}

We express the generators of $\mc{B}_{\mbf{m}}\subset U_{q,t}(\hat{\hat{\mf{sl}}}_n)$ into the matrix coefficients in $\prod_{\mbf{w}}\text{End}(K(\mbf{w}))$. Then we proceed the degeneration limit of $K(\mbf{w})$ via the Chern character map $K(\mbf{w})\rightarrow\widehat{H(\mbf{w})}$. 

Then via the algebra embedding $U_{q,t}(\hat{\hat{\mf{sl}}}_n)\hookrightarrow U_{q,t}(\hat{\mf{g}}_{Q}^{MO})$ as constructed in \cite{N23} and the fact that the degeneration limit of $U_{q,t}(\hat{\mf{g}}_{Q}^{MO})$ lies in $Y_{\hbar_1,\hbar_2}(\mf{g}_{Q}^{MO})\cong Y_{\hbar_1,\hbar_2}(\mf{g}_{Q})$ as stated in \cite{SV17}, we have that the degeneration limit of $\mc{B}_{\mbf{m}}\subset U_{q,t}(\hat{\hat{\mf{sl}}}_n)$ lies in $Y_{\hbar}(\mf{g}_{Q})$.

To get the formula for the matrix element, recall the formula:
\begin{align*}
&\langle\bm{\lambda}|P^{\mbf{m}}_{[i;j)}|\bm{\mu}\rangle=P^{\mbf{m}}_{[i;j)}(\bm{\lambda}\backslash\bm{\mu})\prod_{\blacksquare\in\bm{\lambda}\backslash\bm{\mu}}[\prod_{\square\in\bm{\mu}}\zeta_{ji}(\frac{\chi_{\blacksquare}}{\chi_{\square}})\prod_{k=1}^{\mbf{w}}[\frac{u_k}{q\chi_{\blacksquare}}]],\qquad\mbf{m}\cdot[i,j)\in\mbb{Z}\\
&\langle\bm{\mu}|Q^{\mbf{m}}_{-[i;j)}|\bm{\lambda}\rangle=Q^{\mbf{m}}_{-[i;j)}(\bm{\lambda}\backslash\bm{\mu})\prod_{\blacksquare\in\bm{\lambda}\backslash\bm{\mu}}[\prod_{\square\in\bm{\lambda}}\zeta_{ji}(\frac{\chi_{\square}}{\chi_{\blacksquare}})\prod_{k=1}^{\mbf{w}}[\frac{\chi_{\blacksquare}}{qu_k}]]^{-1},\qquad\mbf{m}\cdot[i,j)\in\mbb{Z}.
\end{align*}

Here:
\begin{equation}\label{Pformula}
P_{\pm[i ; j)}^{\pm \mbf{m}}=\operatorname{Sym}\left[\frac{\prod_{a=i}^{j-1} z_a^{\left\lfloor m_i+\ldots+m_a\right\rfloor-\left\lfloor m_i+\ldots+m_{a-1}\right\rfloor}}{t^{\text {ind } d_{[i, j)}^m} q^{i-j} \prod_{a=i+1}^{j-1}\left(1-\frac{q_2 z_a}{z_{a-1}}\right)} \prod_{i \leq a<b<j} \zeta_{ba}\left(\frac{z_b}{z_a}\right)\right]
\end{equation}
\begin{equation}\label{Qformula}
Q_{\mp[i ; j)}^{\pm \mbf{m}}=\operatorname{Sym}\left[\frac{\prod_{a=i}^{j-1} z_a^{\left\lfloor m_i+\ldots+m_{a-1}\right\rfloor-\left\lfloor m_i+\ldots+m_a\right\rfloor}}{t^{-\mathrm{ind}_{[i, j)}^m} \prod_{a=i+1}^{j-1}\left(1-\frac{q_1 z_{a-1}}{z_a}\right)} \prod_{i \leq a<b<j} \zeta_{ab}\left(\frac{z_a}{z_b}\right)\right].
\end{equation}

Now we shall give the asymptotic behavior of the matrix coefficients $\langle\bm{\lambda}|P^{\mbf{m}}_{[i;j)}|\bm{\mu}\rangle$ in terms of the order of $\kappa$. In the very first order of $\kappa$, we have that:
\begin{align*}
\langle\bm{\lambda}|P^{\mbf{m}}_{[i;j)}|\bm{\mu}\rangle^{coh}=\text{Sym}[\frac{1}{\prod_{a=i+1}^{j-1}(\hbar_2+\chi_{a}-\chi_{a-1})}\prod_{i\leq a<b<j}\omega_{ji}(\chi_b-\chi_a)]\prod_{\blacksquare\in\bm{\lambda}\backslash\bm{\mu}}[\prod_{\square\in\bm{\mu}}\zeta(\frac{\chi_{\blacksquare}}{\chi_{\square}})\prod_{k=1}^{\mbf{w}}[\frac{u_k}{q\chi_{\blacksquare}}]]^{coh}.
\end{align*}

In this way one can find that $\langle\bm{\lambda}|P^{\mbf{m}}_{[i;j)}|\bm{\mu}\rangle^{coh}$ is independent of the choice of $\mbf{m}$. In this case we choose $P^{\mbf{0}}_{[i;j)}$.

We want to prove the following thing:
\begin{prop}\label{Degenerations-for-generators:prop}
$\langle\bm{\lambda}|P^{\mbf{m}}_{[i;j)}|\bm{\mu}\rangle^{coh}$ coincides with $-\langle\bm{\lambda}|E_{[i;j)}|\bm{\mu}\rangle$, and $\langle\bm{\lambda}|Q^{\mbf{m}}_{-[i;j)}|\bm{\mu}\rangle^{coh}$ coincides with $-\langle\bm{\lambda}|E_{-[i;j)}|\bm{\mu}\rangle$. Here $E_{\pm[i;j)}\in\hat{\mf{gl}}_n$. 
\end{prop}
\begin{proof}

The logic of the proof can be given as follows: It is known that the quantum toroidal algebra $U_{q,t}(\hat{\hat{\mf{sl}}}_n)$ action on $K_{T}(M(\mbf{w}))$ will degenerate naturally to the affine Yangian $Y_{\hbar_1,\hbar_2}(\hat{\mf{sl}}_n)$ action on $H_{T}^*(M(\mbf{w}))$. Using the fact that there are algebra embeddings:
\begin{align*}
U_{q,t}(\hat{\hat{\mf{sl}}}_n)\hookrightarrow\prod_{\mbf{w}}\text{End}(K_{T}(M(\mbf{w}))_{loc}),\qquad Y_{\hbar_1,\hbar_2}(\hat{\mf{sl}}_n)\hookrightarrow\prod_{\mbf{w}}\text{End}(H_{T}^*(M(\mbf{w}))_{loc})
\end{align*}
we only need to show that the matrix coefficients of $\langle\bm{\lambda}|P^{\mbf{m}}_{[i;j)}|\bm{\mu}\rangle^{coh}$ and $\langle\bm{\lambda}|Q^{\mbf{m}}_{-[i;j)}|\bm{\mu}\rangle^{coh}$ coincides with $\langle\bm{\mu}|E_{[i,j)}|\bm{\lambda}\rangle$ and $\langle\bm{\mu}|E_{-[i,j)}|\bm{\lambda}\rangle$ respectively. 

Here we only show the proof for $\langle\bm{\lambda}|P^{\mbf{m}}_{[i;j)}|\bm{\mu}\rangle^{coh}$. The proof for $\langle\bm{\lambda}|Q^{\mbf{m}}_{-[i;j)}|\bm{\mu}\rangle^{coh}$ is similar.

Define the following rational functions:
\begin{align*}
P^{coh}_{[i,j)}:=\text{Sym}[\frac{1}{\prod_{a=i+1}^{j-1}(\hbar_2+z_{a}-z_{a-1})}\prod_{i\leq a<b<j}\omega_{ji}(z_b-z_a)]
\end{align*}
which is the degeneration of the element $P^{\mbf{}}_{[i,j)}$ in \eqref{Pformula} under the degeneration limit via sending $\kappa\rightarrow0$ under the variable change $q_i\mapsto e^{\kappa\hbar_i}, z_{i}\mapsto e^{\kappa z_i}$.

It is easy to see that $\langle\bm{\mu}|P^{\mbf{m}}_{[i;j)}|\bm{\lambda}\rangle^{coh}=\langle\bm{\mu}|P^{coh}_{[i,j)}|\bm{\lambda}\rangle$. Thus it remains to prove that 
\begin{equation}\label{key-identity}
P^{coh}_{[i;j)}=-E_{[i;j)}\in\hat{\mf{gl}}_n.
\end{equation}

We use the induction on $[i,j)$, now suppose that the equality \eqref{key-identity} is true for $[i,k)$ and $[k,j)$ for $i<k<j$. We need to prove the following equality:
\begin{align*}
P^{coh}_{[i,j)}=[P^{coh}_{[i,k)},P^{coh}_{[k,j)}]
\end{align*}
since $e_{[i,j)}=[e_{[i,k)},e_{[k,j)}]$.
First we have that $E_{[i,i+1)}=(\frac{1}{q-q^{-1}}P^{\mbf{0}}_{[i,i+1)})^{coh}$ by the direct computation, and we define $x_{i}^{\pm}:=\frac{1}{q-q^{-1}}e_{\pm[i,i+1)}$ so that $x_{i}^{\pm,coh}$ is well-defined as $E_{[i,i+1)}$. Inductively, we also define $x_{[i,j)}^{\pm}:=\frac{1}{q-q^{-1}}e_{\pm[i,j)}$, we will prove that they have the well-defined degeneration limit.

To prove this, note that the isomorphism of the slope subalgebra of slope $\mbf{0}$  $\mc{B}_{\mbf{0}}\cong U_{q}(\hat{\mf{gl}}_n)$ with the quantum affine algebra sends $e_{[i,j)}$ to $P^{\mbf{0}}_{[i,j)}$, $e_{-[i,j)}$ to $Q^{\mbf{0}}_{-[i,j)}$ satisfying the relation \eqref{gln-relation}. By \eqref{gln-relation} we can see that for $1\leq i\leq j\leq k\leq n$:
\begin{align*}
e_{[i,k)}=-\frac{1}{q-q^{-1}}(e_{[i,j)}e_{[j,k)}-q^{-1}e_{[j,k)}e_{[i,j)})
\end{align*}
which is equivalent to
\begin{align*}
x_{[i,k)}=-(x_{[i,j)}x_{[j,k)}-q^{-1}x_{[j,k)}x_{[i,j)}).
\end{align*}

Now we take $q\rightarrow1$ limit, by the induction we have that
\begin{align*}
e_{[i,k)}^{coh}=-E_{[i,k)}.
\end{align*}

Moreover:
\begin{equation*}
\begin{aligned}
&\frac{1}{(q-q^{-1})}[e_{\pm[i,j)}e_{\pm[j,k)}-q^{-1}e_{\pm[j,k)}e_{\pm[i,j)}]=\sum_{i\leq x<j}^{x\equiv k}e_{\pm[j,j+k-x)}e_{\pm[i,x)}-\sum_{i<x\leq j}^{x\equiv j}e_{\pm[x,j)}e_{\pm[i+j-x,k)}\\
=&\sum_{i\leq x<j}^{x\equiv k}e_{\pm[j,j+k-x)}e_{\pm[i,x)}-\sum_{i<x<j}^{x\equiv j}e_{\pm[x,j)}e_{\pm[i+j-x,k)}-e_{\pm[i,k)}.
\end{aligned}
\end{equation*}

Using the induction we have that:
\begin{equation*}
\begin{aligned}
&[x_{\pm[i,j)}x_{\pm[j,k)}-q^{-1}x_{\pm[j,k)}x_{\pm[i,j)}]=(q-q^{-1})\sum_{i\leq x<j}^{x\equiv k}x_{\pm[j,j+k-x)}x_{\pm[i,x)}-(q-q^{-1})\sum_{i<x\leq j}^{x\equiv j}x_{\pm[x,j)}x_{\pm[i+j-x,k)}\\
=&(q-q^{-1})\sum_{i\leq x<j}^{x\equiv k}x_{\pm[j,j+k-x)}x_{\pm[i,x)}-(q-q^{-1})\sum_{i<x<j}^{x\equiv j}x_{\pm[x,j)}x_{\pm[i+j-x,k)}-x_{\pm[i,k)}.
\end{aligned}
\end{equation*}

As $q\rightarrow1$, the quadratic terms on the right hand side vanish.

Hence we have proved that:
\begin{align*}
E_{[i,j)}=-e_{[i,j)}^{coh}
\end{align*}
for arbitrary $i,j$.

\end{proof}

\textbf{Remark.} The above proof is really well-known since it just reflects the fact that the $q\rightarrow1$ of the quantum affine algebra $U_{q}(\hat{\mf{gl}}_n)$ degenerates to the affine Kac-moody algebra $U(\hat{\mf{gl}}_n)$. Also it reveals the fact that the degeneration of the slope subalgebra $\mc{B}_{\mbf{m}}$ will be the universal enveloping algebra $U(\mf{g}_{\mbf{m}})$ of a Lie subalgebra $\mf{g}_{\mbf{m}}$ in $\hat{\mf{gl}}_n$. The point of the proof here is to show the matrix coefficients of $U(\hat{\mf{gl}}_n)$ acting on the equivariant cohomology $H_{T}^*(M(\mbf{w}))$ of quiver varieties, which is useful in the precise computation.

\subsection{Main result}
Having established the degeneration property between quantum affine toroidal algebra and the affine Yangians, here we show our main result:
\begin{thm}\label{degeneration-qde-thm}
The degeneration limit of the quantum difference operator $\mbf{B}_{\mc{L}}^s(z)$ coincides with the quantum multiplication operator $Q(\mc{L})$ up to a constant operator.
\end{thm}
\begin{proof}
To prove the theorem, let us first state the following useful result on the degeneration.
\begin{prop}\label{Degeneration theorem}
With the notation in Theorem 3.5 in \cite{Z23},  the degeneration limit of the monodromy operator $\mbf{B}_{\mbf{m}}(z)$ reads:
\begin{equation*}
\begin{aligned}
\mbf{B}_{\mbf{m}}^{coh}(z)=&(\hbar_1+\hbar_2)\sum_{h=1}^{g}(\sum_{k=1}^{\infty}\frac{1}{1-z^{-k\lvert\bm{\delta}_h\lvert}}\alpha^{\bm{m},h(0),coh}_{k}\alpha^{\bm{m},h(0),coh}_{-k}+\sum_{\substack{\gamma\in\Delta(A)\\m\geq0}}\frac{1}{1-z^{-\mbf{v}_{r}+(m+1)\bm{\delta}_h}}e_{(\delta-\gamma)+m\delta}^{coh}f_{(\delta-\gamma)+m\delta}^{coh}\\
&+\sum_{m>0}\sum_{i,j=1}^{l_h}\frac{1}{1-z^{m\bm{\delta}_h}}u_{m,ij}e_{m\delta,\alpha_i}^{coh}f_{m\delta,\alpha_j}^{coh}+\sum_{\substack{\gamma\in\Delta(A)\\m\geq0}}\frac{1}{1-z^{\mbf{v}_{\gamma}+m\bm{\delta}_h}}e_{\gamma+m\delta}^{coh}f_{\gamma+m\delta}^{coh})
\end{aligned}
\end{equation*}
and here $u_{m,ij}=\text{min}\{i,j\}-\frac{ij}{n}$, and the generators (e.g. $f_{\gamma+m\delta}^{coh}e_{\gamma+m\delta}^{coh})$ ) are the degeneration limits of the generators (e.g.$f_{\gamma+m\delta}e_{\gamma+m\delta})$)under the image of the following map:
\begin{equation*}
\begin{tikzcd}
U_{q,t}(\hat{\hat{\mf{sl}}}_n)\arrow[r,hook]&\prod_{\mbf{w}}\text{End}(K(\mbf{w}))\arrow[r,"ch"]&\prod_{\mbf{w}}\text{End}(\widehat{H}(\mbf{w}))
\end{tikzcd}
\end{equation*}
\end{prop}

\begin{proof}
This is from the definition of the degeneration limit:
\begin{align*}
\mbf{B}^{coh}_{\mbf{m}}(z):=\lim_{\kappa\rightarrow0}\frac{1}{\kappa}(\mbf{B}_{\mbf{m}}(z)-1).
\end{align*}

Then one can use Theorem 3.5 in \cite{Z23} via replacing $q=e^{\kappa(\hbar_1+\hbar_2)}$. One has that the first order of $\kappa$ is given by:
\begin{equation*}
\begin{aligned}
&(\hbar_1+\hbar_2)\sum_{h=1}^g(\sum_{k=1}^{\infty}\frac{1}{1-z^{-k\lvert\bm{\delta}_h\lvert}}\alpha^{\bm{m},h(0),coh}_{k}\alpha^{\bm{m},h(0),coh}_{-k}+\sum_{\substack{\gamma\in\Delta(A)\\m\geq0}}\sum_{l=0}^{\infty}z^{-l(\mbf{v}_{\gamma}+(m+1)\bm{\delta}_h)}e_{(\delta-\gamma)+m\delta}^{coh}f_{(\delta-\gamma)+m\delta}^{coh}\\
+&\sum_{m>0}\sum_{i,j=1}^{l_h}\sum_{l=0}^{\infty}z^{-lm\bm{\delta}_h}u_{m,ij}e_{m\delta,\alpha_i}^{coh}f_{m\delta,\alpha_j}^{coh}+\sum_{\substack{\gamma\in\Delta(A)\\m\geq0}}\sum_{l=0}^{\infty}z^{-l(\mbf{v}_{\gamma}+m\bm{\delta}_h)}e_{\gamma+m\delta}^{coh}f_{\gamma+m\delta}^{coh})
\end{aligned}
\end{equation*}
which can be written in a compact form:
\begin{equation*}
\begin{aligned}
&(\hbar_1+\hbar_2)\sum_{h=1}^g(\sum_{k=1}^{\infty}\frac{1}{1-z^{-k\lvert\bm{\delta}_h\lvert}}\alpha^{\bm{m},h(0),coh}_{k}\alpha^{\bm{m},h(0),coh}_{-k}+\sum_{\substack{\gamma\in\Delta(A)\\m\geq0}}\frac{1}{1-z^{-(\mbf{v}_{\gamma}+(m+1)\bm{\delta}_h)}}e_{(\delta-\gamma)+m\delta}^{coh}f_{(\delta-\gamma)+m\delta}^{coh}\\
+&\sum_{m>0}\sum_{i,j=1}^{l_h}\frac{1}{1-z^{-m\bm{\delta}_h}}u_{m,ij}e_{m\delta,\alpha_i}^{coh}f_{m\delta,\alpha_j}^{coh}+\sum_{\substack{\gamma\in\Delta(A)\\m\geq0}}\frac{1}{1-z^{-(\mbf{v}_{\gamma}+m\bm{\delta}_h)}}e_{\gamma+m\delta}^{coh}f_{\gamma+m\delta}^{coh}).
\end{aligned}
\end{equation*}

Then using the Proposition \ref{degeneration-of-algebra} we can show that all the generators $f_{\gamma+m\delta}^{coh}e_{\gamma+m\delta}^{coh}$ lie in the affine Yangian $Y_{\hbar_1,\hbar_2}(\hat{\mf{sl}}_n)$. Therefore we have obtained the formula in the Proposition.

\end{proof}

Now we choose a path such that each $\mbf{m}$ is at the generic points. In this case the degeneration limit of the monodromy operator can be written as one of the following form:
\begin{align*}
\sum_{h=1}^{n}(\sum_{k=1}^{\infty}\frac{1}{1-z^{-k\lvert\bm{\delta}_h\lvert}}\alpha^{\bm{m},h(0),coh}_{k}\alpha^{\bm{m},h(0),coh}_{-k}),\qquad\frac{1}{1-z^{-\alpha_{i}^{\mbf{m}}}}e_{\alpha_i}^{\mbf{m},coh}e_{-\alpha_i}^{\mbf{m},coh}.
\end{align*}

Suppose that the quantum difference operator $\textbf{M}_{\mc{L}}(z)$ being written as:
\begin{align*}
\textbf{M}_{\mc{L}}(z)=\mc{L}\textbf{B}_{w_m}(z)\cdots\textbf{B}_{w_0}(z)
\end{align*}
such that $w_m,\cdots,w_0$ are the slope points on the path from $s$ to $s-\mc{L}$. 
Under the degenerate limit in Proposition \ref{Degenerations-for-generators:prop}, the formula can be written as:
\begin{equation*}
\begin{aligned}
(\hbar_1+\hbar_2)\sum_{\mbf{m}\in[s,s-\mc{L})}\sum_{[i,j)}^{\mbf{m}\cdot[i,j)\in\mbb{Z}}\frac{1}{1-z^{-[i;j)}}E_{[i,j)}E_{-[i,j)}
\end{aligned}
\end{equation*}
with $\bm{\theta}=(1,\cdots,1)$.

Note that each $Q^{\mbf{m}}_{-[i;j)}P^{\mbf{m}}_{[i;j)}$ has the degeneration limit as $E_{-[i,j)}E_{[i,j)}$. With this fact it remains to match the coefficients of the Dubrovin connection in \eqref{good-form-for-Dubrovin-connection} on each term. i.e. The coefficients of the formula of the form $\frac{1}{1-z^{-[i,j)}}E_{[i,j)}E_{-[i,j)}$ in $Q(\mc{L})$ is given by $\mc{L}\cdot[i,j)$.

For a strictly increasing monotone path $\mu:[0,1]\rightarrow\mbb{R}^n$ from $s$ to $s-\mc{L}$ such that $\mu(0)=s-\mc{L}$ and $\mu(1)=s$. This will associate a function $f:[0,1]\rightarrow\mbb{R}$ defined as $f(t)=\mu(t)\cdot[i,j)$. It is easy to check that $f(t)$ is strictly increasing. We also have $f(1)=\mu(1)\cdot[i,j)=s\cdot[i,j)$ and $f(0)=s\cdot[i,j)-\mc{L}\cdot[i,j)$. Since we are choosing the interval $[0,1)$, hence $f([0,1))=[s\cdot[i,j),s\cdot[i,j)-\mc{L}\cdot[i,j))$, which contains exactly $\mc{L}\cdot[i,j)$ integers. Thus the coefficients of the form $\frac{1}{1-z^{-[i,j)}}E_{[i,j)}E_{-[i,j)}$ is given by $f([0,1))\cap\mbb{Z}$, which coincides with $\mc{L}\cdot[i,j)$.

\end{proof}

\subsection{$p\rightarrow1$ limit of the solution of the quantum difference equation}\label{sub:_p_rightarrow1_limit_of_the_solution_of_the_quantum_difference_equation}

By the result of the degeneration limit of the quantum difference equation in Theorem \ref{degeneration-qde-thm}, we can start to construct the $p\rightarrow1$ degeneration limit of the solution to the quantum difference equation.

Let $\Psi_{0,\infty}(q_1,q_2,p,z)$ be the solution to the quantum difference equation described above. Then we take $z=e^{2\pi i s}$, $q_i=e^{2\pi i\hbar_i\tau}$, $q=e^{-2\pi i\tau}$. By the Theorem \ref{degeneration-qde-thm}, for the degeneration limit of the solution $\psi_{0,\infty}(z)$, which is defined as follows:
\begin{align*}
\psi_{0,\infty}(z)=\lim_{\tau\rightarrow0}\Psi_{0,\infty}(e^{2\pi i\hbar_1\tau},e^{2\pi i\hbar_2\tau},e^{-2\pi i\tau},z)\in \widehat{H_{T}(M(\mbf{v},\mbf{w}))}_{loc}.
\end{align*}

It is the solution to the corresponding Dubrovin connection. The solution lies in the completion of $H_{T}(M(\mbf{v},\mbf{w}))_{loc}$. Using this, we can define the transport of the solution of the Dubrovin connection:
\begin{align*}
\text{Trans}(s):=\psi_{0}(e^{2\pi is})^{-1}\psi_{\infty}(e^{2\pi is})\in \widehat{H_{T}(M(\mbf{v},\mbf{w}))}_{loc}
\end{align*}
Similarly as a result of the Theorem \ref{degeneration-qde-thm}, the transport of the solution is a limit of the monodromy:
\begin{align}\label{trans}
\text{Trans}(s)=\lim_{\tau\rightarrow0}\textbf{Mon}(z=e^{2\pi is},t_1=e^{2\pi i\hbar_1\tau},t_2=e^{2\pi i\hbar_2\tau},q=e^{-2\pi i\tau})
\end{align} 
in $\widehat{H_{T}(M(\mbf{v},\mbf{w}))}_{loc}$.
Here $s=(s_1,\cdots,s_n)\in\mbb{R}^n$ is the generic point in $\mbb{R}^n$.

The description for the monodromy operator from the monodromy operator generally can only be written as the infinite product of matrices. The elliptic stable envelope of \cite{AO21} would give a more compact description, which we do not pursue here.

\begin{prop}\label{identi-p-0-connection}
$\text{Trans}(s)=\lim_{p\rightarrow0}\textbf{Mon}^{reg}(p^s,e^{2\pi i\hbar_1},e^{2\pi i\hbar_2},p)$ for $s\in\mbb{R}^n\backslash\text{Walls}$
\end{prop}
\begin{proof}
Note that $\textbf{Mon}(q,t_1,t_2,z)$ can be written as the following:
\begin{align*}
\exp(-\sum_{i}\frac{\ln(E_{i}^{(0)})\ln(z_i)}{\ln(q)})\sum_{\mbf{a},\mbf{b}}L_{\mbf{a},\mbf{b}}\prod_{j=1}^{r}\prod_{i_j,l_j}\frac{\theta(z_{i_j}^{l_j}a_{i_j})}{\theta(z_{i_j}^{l_j}b_{i_j})}\exp(\sum_{i}\frac{\ln(E_{i}^{(\infty)})\ln(z_i)}{\ln(q)}).
\end{align*}

Using the formula \eqref{elliptic-to-multivalued}:
\begin{align*}
\exp(-\sum_{i}\frac{\ln(E_{0}^{(i)})\ln(z_i)}{\ln(p)})\textbf{Mon}(z)\exp(\sum_{i}\frac{\ln(E_{\infty}^{(i)})\ln(z_i)}{\ln(p)})=\prod_{i}(\Theta(E_{0}^{(i)},z_i))^{-1}\textbf{Mon}^{ell}(z)\prod_{i}(\Theta(E_{\infty}^{(i)},z_i)).
\end{align*}

It is known that $\textbf{Mon}^{ell}(z)$ is the linear combination of Riemann theta function as stated in \eqref{theta-function-form-of-ell-connection}, and it is $p$-periodic. Since $\textbf{Mon}(p,t_1,t_2,z)$ is also $p$-periodic, by definition it is written as:
\begin{align*}
\textbf{Mon}(p,t_1,t_2,z)=\exp(-\sum_{i}\frac{\ln(E_{0}^{(i)})\ln(z_i)}{\ln(p)})\textbf{Mon}^{reg}(p,t_1,t_2,z)\exp(\sum_{i}\frac{\ln(E_{\infty}^{(i)})\ln(z_i)}{\ln(p)})
\end{align*}
with $\textbf{Mon}(p,t_1,t_2,z)$ written in the form \eqref{theta-function-form-of-reg-conn}. This means that $\textbf{Mon}(p,t_1,t_2,z)$ is written in the linear combination of the following form:
\begin{align*}
\exp(\sum_{i}\frac{\ln(E_{\infty}^{(i)}/E_{0}^{(i)})\ln(z_i)}{\ln(p)})\prod_{\mbf{k}}\frac{\theta(z_1^{k_1}\cdots z_n^{k_n}a_{\mbf{k}})}{\theta(z_1^{k_1}\cdots z_n^{k_n}b_{\mbf{k}})}
\end{align*}
with 
\begin{align}\label{periodicity-condition}
\prod_{\mbf{k}}(\frac{a_{\mbf{k}}}{b_{\mbf{k}}})^{k_i}=E_{(0)}^{(i)}/E_{\infty}^{(i)}
\end{align}
since $\textbf{Mon}(p,t_1,t_2,zp^{\mc{L}})=\textbf{Mon}(p,t_1,t_2,z)$.

We can use the following identity \eqref{modular-transform} proved in the appendix:
\begin{align*}
\lim_{\tau\rightarrow0}\exp(2\pi is(\epsilon_1-\epsilon_2))\frac{\theta(e^{2\pi i(s+\tau\epsilon_1)},e^{2\pi i\tau})}{\theta(e^{2\pi i(s+\tau\epsilon_2)},e^{2\pi i\tau})}=\lim_{q\rightarrow0}\frac{\theta(q^se^{2\pi i\epsilon_1},q)}{\theta(q^se^{2\pi i\epsilon_2},q)}.
\end{align*}

This identity can be generalized to:
\begin{align*}
\lim_{\tau\rightarrow0}\exp(2\pi i(k\cdot s)\lambda)\frac{\prod_{i}\theta(e^{2\pi i(k\cdot s+\tau a_i)},e^{2\pi i\tau})}{\prod_{j}\theta(e^{2\pi i(k\cdot s+\tau b_j)},e^{2\pi i\tau})}=\lim_{q\rightarrow0}\frac{\prod_{i}\theta(q^{k\cdot s}e^{2\pi i a_i},q)}{\prod_{j}\theta(q^{k\cdot s}e^{2\pi i b_j},q)},\lambda=\sum_{i}a_i-\sum_{j}b_j.
\end{align*}

Then using the formula to our setting, we have:
\begin{equation*}
\begin{aligned}
&\lim_{\tau\rightarrow0}\exp(\sum_{i}\frac{\ln(E_0^{(i)}/E_{\infty}^{(i)})\ln(z_i)}{\ln(p)})\prod_{\mbf{k}}\frac{\theta(z_1^{k_1}\cdots z_n^{k_n}a_{\mbf{k}})}{\theta(z_1^{k_1}\cdots z_n^{k_n}b_{\mbf{k}})}\\
=&\lim_{\tau\rightarrow0}\exp(\sum_{i}\frac{\ln(E_0^{(i)}/E_{\infty}^{(i)})\ln(z_i)}{\ln(p)})\exp(-\sum_{\mbf{k}}\frac{k_i\ln(z_i)\ln(a_{\mbf{k}}/b_{\mbf{k}})}{\ln(p)})\exp(\sum_{\mbf{k}}\frac{k_i\ln(z_i)\ln(a_{\mbf{k}}/b_{\mbf{k}})}{\ln(p)})\times\\
&\times\prod_{\mbf{k}}\frac{\theta(z_1^{k_1}\cdots z_n^{k_n}a_{\mbf{k}})}{\theta(z_1^{k_1}\cdots z_n^{k_n}b_{\mbf{k}})}=\lim_{q\rightarrow0}\prod_{\mbf{k}}\frac{\theta(q^{\mbf{k}\cdot\mbf{s}}a_{\mbf{k}})}{\theta(q^{\mbf{k}\cdot\mbf{s}}b_{\mbf{k}})}.
\end{aligned}
\end{equation*}

Here we use the fact that $E_{(0)}^{(i)}/E_{\infty}^{(i)}=\prod_{\mbf{k}}(\frac{a_{\mbf{k}}}{b_{\mbf{k}}})^{k_i}$ as stated in \eqref{periodicity-condition}. Now we have finished the proof of the theorem.
\end{proof}

\textbf{Remark. }Readers can also compare the proof with the proof of Proposition 12 in \cite{S21}. The proof over here uses the intrinsic fact that the connection matrix for the quantum difference equation is written in the balanced form of the theta functions as stated in Section 8.2 in \cite{S21}. This is also the key point in the proof of Proposition 12 in \cite{S21} where the balanced property of the elliptic stable envelope was used. Readers might also do the similar computation on the side of the elliptic stable envelopes using the formula given in \cite{S20}\cite{D21}.

\begin{prop}
Let $\gamma$ be a path from $z=0$ to $z=\infty$ meeting the torus $|z_i|=1$ at a nonsingular point $z=e^{2\pi is}$ with all $z_i$ equal, and suppose that $s$ does not lie on the wall. Then
\begin{equation*}
\text{Trans}(s)=
\begin{cases}
\prod^{\leftarrow}_{w\in(0,s)}(\mbf{B}_{w}^*)^{-1}\cdot\mbf{T}&s\geq0\\
\prod^{\rightarrow}_{w\in(s,0)}\mbf{B}_{w}^{*}\cdot\mbf{T}&s<0
\end{cases}.
\end{equation*}
\end{prop}

Moreover, if we relax the condition of $(z_1,\cdots,z_n)$ to the torus $\{|z_i|=1\}\subset\mbb{C}^n$, one could obtain the following result for the transition operator:
\begin{thm}
The transport of the quantum connection from $z=0$ to $z=\infty$ intersecting the polydisk at a nonsingular point $z=e^{2\pi is}$ with $|z_i|=1$ for each $i=1,\cdots,n$, and $s$ does not meet the wall set, then it equals:
\begin{align*}
\text{Trans}(s)=
\begin{cases}
\prod_{w_i\in\text{Walls}_i\subset\text{Walls}\cap C(s)^c}^{\rightarrow}(\mbf{B}_{w_i}^{*})\prod_{w_i\in\text{Walls}_i\subset\text{Walls}^-\cap D(s)^c}^{\leftarrow}(\mbf{B}_{w_i}^{*})^{-1}\\s\in B_{k}\cap B_{-l}\\
\prod_{w_i\in\text{Walls}_i\subset\text{Walls}^-\cap D(s)^c}^{\leftarrow}(\mbf{B}_{w_i}^{*})^{-1},\qquad s\in\bigcup_{w\in\text{Walls}_0}\{\mbf{s}\in\mbb{R}^n|s_i-l\leq w_i\}\\
\prod_{w_i\in\text{Walls}_i\subset\text{Walls}\cap C(s)^c}^{\rightarrow}(\mbf{B}_{w_i}^{*}),\qquad s\in\bigcup_{w\in\text{Walls}_0}\{\mbf{s}\in\mbb{R}^n|s_i+k\geq w_i\}
\end{cases}.
\end{align*}
where $\text{Walls}^-:=\text{Walls}\cap(\mbb{R}^{-})^n$.
\end{thm}
\begin{proof}
Combine the nodal limits of $\psi^{reg}_{0}(p^sz)$ and $\psi^{reg}_{\infty}(p^sz)$ from Proposition \ref{p0-limit-z-0} with Proposition \ref{p0-limit-z-infty}. Then combine them with the Proposition \ref{identi-p-0-connection}, we obtain the theorem.
\end{proof}

This formula has the following consequences.

Suppose we are given a solution $\psi_{0}(z)$ of the quantum difference equation and extend it to $\psi_{\infty}(z)$ via the curve $\gamma$ intersecting the unit polydisk at $e^{2\pi is}$, $s\in\mbb{R}^n$, we know that the transport formula is given by $\text{Trans}_{DT}(s)$. Now we move the curve $\gamma$ a little bit such that the new curve $\gamma'$ intersect the unit polydisk with $e^{2\pi is'}$. Now we require that the move of $s$ and $s'$ is close enough such that:
\begin{align*}
(\text{Walls}\cap C(s)^c\cap C(s')^c)^c_{\text{Walls}}\cap(\text{Walls}^-\cap D(s)^c\cap D(s')^c)^{c}_{\text{Walls}^-}\text{ has only one element}.
\end{align*}

In human language, it means that the corresponding wall elements of $s$ and $s'$ has only one element difference. In this way, inversing the curve $\gamma'$ and consider the loop $\gamma'^{-1}\gamma$, this corresponds to the element $\text{Trans}(s')^{-1}\text{Trans}(s)$, and if we denote the wall of the difference by $w$, we can see that:
\begin{align*}
\text{Trans}(s')^{-1}\text{Trans}(s)=\mbf{B}_{w}^*.
\end{align*}

Since the analytic continuation along the loop $\gamma'^{-1}\gamma$ depends only on its homotopy class, we have defined a map:
\begin{align*}
\pi_{1}((\mbb{C}^*)^n-\text{Sing})\rightarrow\text{Aut}(\widehat{H_{T}}(M(\mbf{v},\mbf{w})))
\end{align*}
with each loop $\gamma$ mapped to an ordered product of $\prod_{w}\mbf{B}_{w}^*$.

In conclusion, we have the following monodromy representation:
\begin{thm}\label{monodromy-rep-casimir}
The monodromy representation:
\begin{align*}
\pi_{1}(\mbb{P}^n\backslash\textbf{Sing}, 0^+)\rightarrow\text{Aut}(\widehat{H_{T}}(M(\mbf{v},\mbf{w})))
\end{align*}
of the Dubrovin connection is generated by $\mbf{B}_{\mbf{m}}^*:=P^{-1}\mbf{B}_{\mbf{m}}P$ with $q_1=e^{2\pi i\hbar_1},q_2=e^{2\pi i\hbar_2}$, i.e. the monodromy operators $\mbf{B}_{\mbf{m}}$ in the image of the $K$-theoretic fixed point basis under the Chern character map. $0^+$ is a point infinitesimally near $0$.
\end{thm}

\section{\textbf{Case study: the equivariant Hilbert scheme of $A_{r-1}$ singularity}}\label{sec:_textbf_case_study_the_equivariant_hilbert_scheme_of_a_r_singularity}
In this section we study the equivariant Hilbert scheme $\text{Hilb}_{n}([\mbb{C}^2/\mbb{Z}_r])$ of $[\mbb{C}^2/\mbb{Z}_r]$ in detail.

As a quiver variety, the corresponding quiver representation is written as:
\begin{align*}
\text{Rep}(\mbf{v},\mbf{w})=\bigoplus_{i\in\mbb{Z}/r\mbb{Z}}\text{Hom}(V_i,W_i)\oplus\text{Hom}(V_0,W_0),\qquad V_i=\mbb{C}^n,W_0=\mbb{C}.
\end{align*}

The fixed point of the equivariant Hilbert scheme $\text{Hilb}_{n}([\mbb{C}^2/\mbb{Z}_r])$ is denoted by the partitions $\lambda$ such that $|\lambda|=nr$ such that:
\begin{align*}
\#\{\square\in\lambda|c_{\square}\text{mod }r=i\}=n.
\end{align*}

\subsection{Matrix coefficients of the quantum difference operators for $\text{Hilb}_{n}([\mbb{C}^2/\mbb{Z}_r])$}
Now we use the above example to give calculation of the wall structure of the Hilbert scheme $\text{Hilb}_{n}([\mbb{C}^2/\mbb{Z}_r])$. Also we give the potential methods for the computation of the matrix coefficients of the algebraic quantum difference operator in terms of the stable basis and fixed point basis. For the introduction of the stable basis, one can refer to Section 6 in \cite{N15} for the further details. Here we use the notation $s_{\bm{\lambda}}^{\mbf{m}}\in K_{T}(M(\mbf{v},\mbf{w}))$ as the stable basis of slope $\mbf{m}\in\mbb{Q}^r$ with respect to the multipartition $\bm{\lambda}$.

From the definition it is easy to see that the wall hyperplane in $\text{Pic}(\text{Hilb}_{n}([\mbb{C}^2/\mbb{Z}_r]))\otimes\mbb{Q}\cong\mbb{Q}^r$ for the Hilbert scheme $\text{Hilb}_{n}([\mbb{C}^2/\mbb{Z}_r])$ is given by:
\begin{align*}
\mbf{m}\cdot[i,j)\in\mbb{Z},\qquad1\leq i\leq j\leq nr.
\end{align*}

The finite hyperplane through the origin is given by:
\begin{align*}
\mbf{m}\cdot[i,j)=0.
\end{align*}

In this case we choose $s\in\mbb{R}^r$ to be infinitesimally close to zero, and the quantum difference operator can be written as:
\begin{equation*}
\mbf{B}_{\mc{L}}^s(z)=\mc{L}\prod_{\substack{\mbf{m}\in\text{generic points of walls}\\ \mbf{z}\cdot[i,j)=n, 1\leq i\leq j\leq nr}}\mbf{B}_{\mbf{m}}(z)
\end{equation*}
such that:
\begin{align*}
\mbf{B}_{\mbf{m}}(z)=
\begin{cases}
\sum_{n=0}^{\infty}\frac{(q-q^{-1})^n}{[n]_{q^2}!}\frac{(-1)^n}{\prod_{l=1}^{n}(1-z^{n[i,j)}p^{-n\mbf{m}\cdot[i,j)}q^{-ln\bm{\theta}\cdot((\frac{n^2r-1}{2}\bm{\theta}+\mbf{e}_1)-2-2\delta_{1l})})}(Q^{\mbf{m}}_{-[i,j)})^n(P^{\mbf{m}}_{[i,j)})^n&(j-i)\text{mod }r\neq0\\
\prod_{h=1}^{r}:\exp(-\sum_{k=1}\frac{n_kq^{-\frac{k|\bm{\delta}_{h}|}{2}}}{1-z^{-k|\bm{\delta}_h|}p^{k\mbf{m}\cdot\bm{\delta}_h}q^{-\frac{k|\bm{\delta}_{h}|}{2}}}\alpha^{\mbf{m},h}_{-k}\alpha_{k}^{\mbf{m},h}):& (j-i)\text{mod }r=0
\end{cases}.
\end{align*}

Here we introduce two ways of computing the matrix coefficients of quantum difference operator in terms of stable basis and fixed point basis.

For the stable basis case, it has been shown in Section 6 in \cite{N15} that one could write down the matrix coefficients of $\mbf{B}_{\mbf{m}}(z)$ in terms of the normalized stable basis:
\begin{align*}
&P^{\mbf{m}}_{[i,j)}\cdot s_{\bm{\mu}}^{\mbf{m}}=\sum_{\text{cavalcade of }\mbf{m}\text{-ribbons}}^{\bm{\lambda}\backslash\bm{\mu}=C\text{ is a type }[i,j)}s_{\bm{\lambda}}^{\mbf{m}}(1-q^2)^{\#_C}q^{\text{ht C}+\text{ind}_{C}^{\mbf{m}}+N_C^+}\\
&Q^{\mbf{m}}_{-[i,j)}\cdot s_{\bm{\mu}}^{\mbf{m}}=\sum_{\text{stampede of }\mbf{m}\text{-ribbons}}^{\bm{\lambda}\backslash\bm{\mu}=S\text{ is a type }[i,j)}s_{\bm{\mu}}^{\mbf{m}}(1-q^2)^{\#_S}q^{\text{wd}(S)-\text{ind}_{S}^{\mbf{m}}-j+i+N_{S}^-}.
\end{align*}

Using the formula, we have that:
\begin{equation*}
\begin{aligned}
(Q^{\mbf{m}}_{-[i,j)})^n(P^{\mbf{m}}_{[i,j)})^n\cdot s^{\mbf{m}}_{\bm{\lambda}_0}=&\sum_{\text{cavalcade of }\mbf{m}\text{-ribbons}}^{\bm{\lambda}_i\backslash\bm{\lambda}_{i-1}=C_i\text{ is a type }[i,j)}(Q^{\mbf{m}}_{-[i,j)})^ns^{\mbf{m}}_{\bm{\lambda}_n}(1-q^2)^{\sum_{i=1}^n\#_{C_i}}q^{\sum_{i=1}^n(\text{ht }C_i+\text{ind}_{C_i}^{\mbf{m}}+N_{C_i}^+)}\\
=&\sum_{\text{stampede of }\mbf{m}\text{-ribbons}}^{\bm{\lambda}_{n+i-1}\backslash\bm{\lambda}_{n+i}=S_i\text{ is a type }[i,j)}\sum_{\text{cavalcade of }\mbf{m}\text{-ribbons}}^{\bm{\lambda}_i\backslash\bm{\lambda}_{i-1}=C_i\text{ is a type }[i,j)}s^{\mbf{m}}_{\lambda_{2n}}(1-q^2)^{\sum_{i=1}^n\#_{C_i}+\#_{S_i}}\times\\
&q^{\sum_{i=1}^n(\text{ht }C_i+\text{ind}_{C_i}^{\mbf{m}}+N_{C_i}^+)}
\times q^{\sum_{i=1}^n(\text{wd}(S_i)-\text{ind}^{\mbf{m}}_{S_i}-j+i+N_{S_i}^-)}.
\end{aligned}
\end{equation*}

In conclusion, the matrix coefficients of the quantum difference operators $\mbf{B}_{\mc{L}}^s(z)$ can be written as:
\begin{equation*}
\begin{aligned}
(s^{\mbf{m}}_{\bm{\mu}},\mbf{B}_{\mc{L}}^s(z)s^{\mbf{m}}_{\bm{\lambda}})=&(s^{\mbf{m}}_{\bm{\mu}},\mc{L}\prod_{\substack{\mbf{m}\in\text{generic points of walls}\\ \mbf{z}\cdot[i,j)=n, 1\leq i\leq j\leq nr}}\mbf{B}_{\mbf{m}}(z)s^{\mbf{m}}_{\bm{\lambda}})\\
=&(s^{\mbf{m}}_{\bm{\mu}},\mc{L}s^{\mbf{m}}_{\bm{\lambda}_{o}})\prod_{\mbf{m},\bm{\lambda}_i}(s^{\mbf{m}}_{\bm{\lambda}_i},\mbf{B}_{\mbf{m}}(z)s^{\mbf{m}}_{\bm{\lambda}_{i-1}})R^{\mbf{m},\mbf{m}'}_{\bm{\lambda}_{i+1},\bm{\lambda}_{i}}.
\end{aligned}
\end{equation*}

Here $\mbf{m}'$ is the slope immediately following $\mbf{m}$ with respect to the order of the product. $R^{\mbf{m},\mbf{m}'}_{\bm{\lambda}_{i+1},\bm{\lambda}_{i}}$ is the $K$-theoretic wall $R$-matrix of slopes $\mbf{m}$ and $\mbf{m}'$ at the fixed point basis $\bm{\lambda}_{i+1}$ and $\bm{\lambda}_{i}$. 

In the fixed point basis, recall the formula:
\begin{align*}
&\langle\bm{\lambda}|P^{\mbf{m}}_{[i,j)}|\bm{\mu}\rangle=P^{\mbf{m}}_{[i,j)}(\bm{\lambda}\backslash\bm{\mu})\prod_{\blacksquare\in\bm{\lambda}\backslash\bm{\mu}}[\prod_{\square\in\bm{\mu}}\zeta_{ij}(\frac{\chi_{\blacksquare}}{\chi_{\square}})[\frac{u_1}{q\chi_{\blacksquare}}]],\qquad i=\text{color of }\blacksquare, j=\text{color of }\square\\
&\langle\bm{\mu}|Q^{\mbf{m}}_{-[i,j)}|\bm{\lambda}\rangle=Q^{\mbf{m}}_{-[i,j)}(\bm{\lambda}\backslash\bm{\mu})\prod_{\blacksquare\in\bm{\lambda}\backslash\bm{\mu}}[\prod_{\square\in\bm{\lambda}}\zeta_{ji}(\frac{\chi_{\square}}{\chi_{\blacksquare}})[\frac{\chi_{\blacksquare}}{qu_1}]]^{-1},\qquad i=\text{color of }\blacksquare, j=\text{color of }\square.
\end{align*}

So we can write that:
\begin{equation*}
\begin{aligned}
&(Q^{\mbf{m}}_{-[i,j)})^n(P^{\mbf{m}}_{[i,j)})^n|\bm{\lambda}_0\rangle=(Q^{\mbf{m}}_{-[i,j)})^n|\bm{\lambda}_n\rangle\prod_{l=1}^{n}P^{\mbf{m}}_{[i,j)}(\bm{\lambda}_{l}\backslash\bm{\lambda}_{l-1})\prod_{\blacksquare_l\in\bm{\lambda}_{l}\backslash\bm{\lambda}_{l-1}}[\prod_{\square_{l-1}\in\bm{\lambda}_{l-1}}\zeta(\frac{\chi_{\blacksquare_l}}{\chi_{\square_{l-1}}})[\frac{u_1}{q\chi_{\blacksquare_l}}]]\\
=&\prod_{l=1}^{n}Q^{\mbf{m}}_{-[i,j)}(\bm{\lambda}_{n+l-1}\backslash\bm{\lambda}_{n+l})P^{\mbf{m}}_{[i,j)}(\bm{\lambda}_{l}\backslash\bm{\lambda}_{l-1})\prod_{\substack{\blacksquare_{l}\in\bm{\lambda}_{l}\backslash\bm{\lambda}_{l-1}\\\blacksquare_{n+l}\in\bm{\lambda}_{n+l}\backslash\bm{\lambda}_{n+l-1}}}[\prod_{\substack{\square_{l-1}\in\bm{\lambda}_{l-1}\\\square_{n+l-1}\in\bm{\lambda}_{n+l-1}}}\frac{\zeta(\frac{\chi_{\blacksquare_l}}{\chi_{\square_{l-1}}})}{\zeta(\frac{\chi_{\square_{n+l-1}}}{\chi_{\blacksquare_{n+l}}})}\frac{[\frac{u_1}{q\chi_{\blacksquare_l}}]}{[\frac{\chi_{\blacksquare_{n+l}}}{qu_1}]}]]|\bm{\lambda}_{2n}\rangle.
\end{aligned}
\end{equation*}

\subsection{Example: $\text{Hilb}_{3}([\mbb{C}^2/\mbb{Z}_2])$}
Now we give the example of the equivariant Hilbert scheme $\text{Hilb}_{3}([\mbb{C}^2/\mbb{Z}_2])$. As the quiver variety, the corresponding quiver representation is
\begin{align*}
\text{Hom}(\mbb{C}^3,\mbb{C}^3)^{\oplus 2}\oplus\text{Hom}(\mbb{C}^3,\mbb{C}).
\end{align*}

The $T$-fixed point of the equivariant Hilbert scheme $\text{Hilb}_{3}([\mbb{C}^2/\mbb{Z}_2])$ is given by the single partition $\lambda$ such that $|\lambda|=6$ and such that for the boxes $\square\in\lambda$ such that $\#\{\square\in\lambda|c_{\square}\text{ mod }2=0\}=\#\{\square\in\lambda|c_{\square}\text{ mod }2=1\}=3$. Here $c_{\square}=i-j$ with $(i,j)$ the coordinate of the box. Written in the Young diagram, the fixed point set corresponds to the following Young diagrams:

\begin{equation*}
\ytableausetup{textmode, boxframe=0.1em, boxsize=1em}
\begin{ytableau}
0 \\
1 \\
0 \\
1 \\
0 \\
1 \\
\end{ytableau}
\qquad
\begin{ytableau}
0&1 \\
1 \\
0 \\
1 \\
0 \\
\end{ytableau}
\qquad
\begin{ytableau}
0&1\\
1&0\\
0\\
1\\
\end{ytableau}
\qquad
\begin{ytableau}
0&1\\ 1&0\\ 0&1
\end{ytableau}
\qquad
\begin{ytableau}
0&1&0\\ 1\\ 0\\ 1
\end{ytableau}
\qquad
\begin{ytableau}
0&1&0\\1&0&1
\end{ytableau}
\end{equation*}
\begin{equation*}
\begin{ytableau}
0&1&0&1\\1&0
\end{ytableau}
\qquad
\begin{ytableau}
0&1&0&1&0\\1
\end{ytableau}
\qquad
\begin{ytableau}
0&1&0&1&0&1
\end{ytableau}
\end{equation*}

The wall structure on $\text{Hilb}_{3}([\mbb{C}^2/\mbb{Z}_2])$ is given by the following hyperplanes in $\mbb{R}^2$:
\begin{equation*}
\begin{aligned}
&x=n_1,\qquad y=n_2,\qquad x+y=n_3,\qquad x+2y=n_4\\
&2x+y=n_5,\qquad 2x+2y=n_6,\qquad 3x+2y=n_7,\qquad 2x+3y=n_8,\qquad 3x+3y=n_9.
\end{aligned}
\end{equation*}

The following is a graph of the wall structure for $\text{Hilb}_{3}([\mbb{C}^2/\mbb{Z}_2])$ in the neighborhood of $[0,1]\times[0,1]$:
\begin{center}
\includegraphics[width=0.5\textwidth]{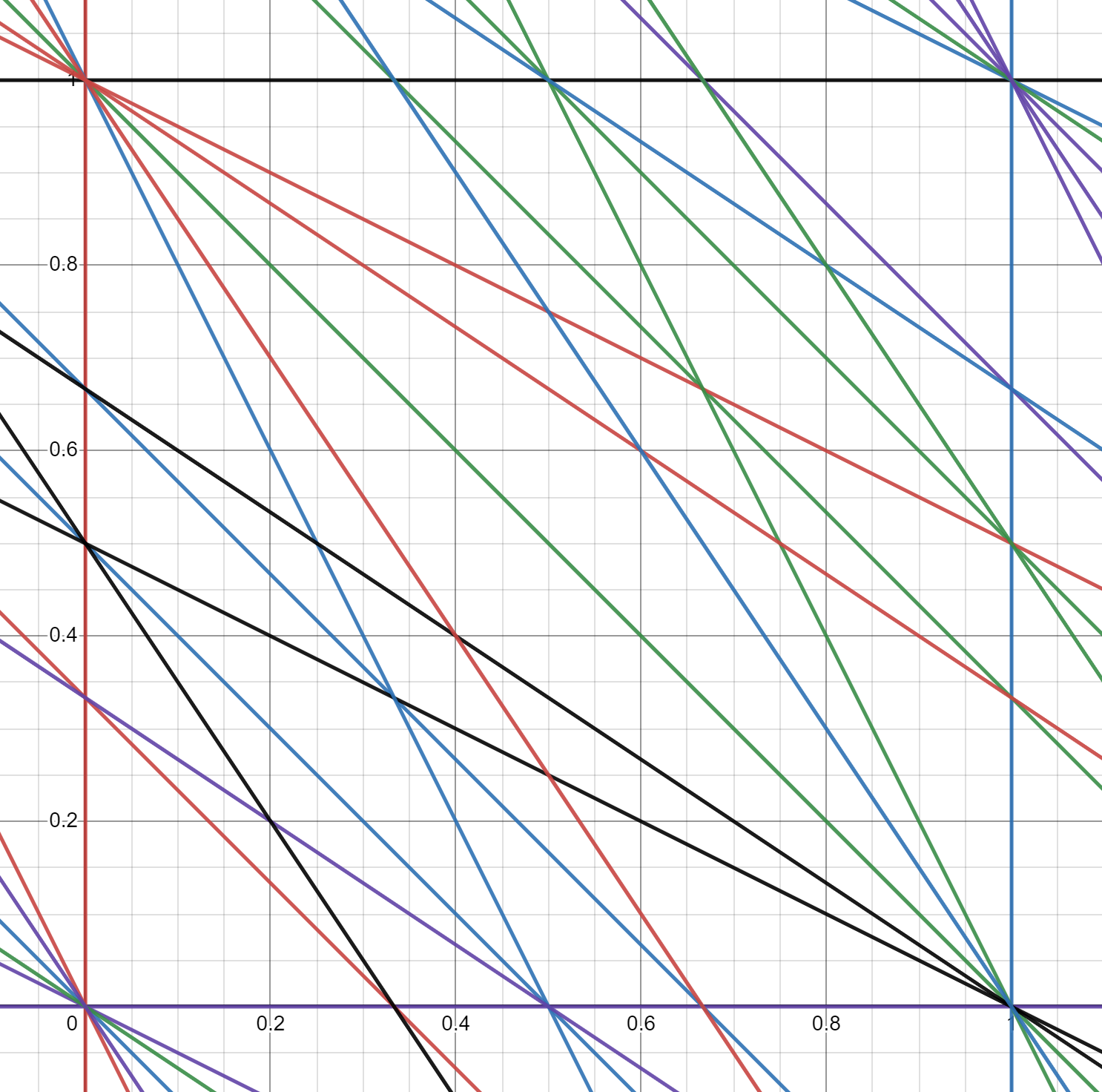}
\end{center}

A direct computation shows that for the generic slope point $\mbf{m}$ on the wall $x=n_1$, $y=n_2$, $x+2y=n_4$, $2x+y=n_5$, $3x+2y=n_7$, $2x+3y=n_8$, the monodromy operator $\mbf{B}_{\mbf{m}}$ is of the $U_{q}(\mf{sl}_2)$ type. For the generic slope points $\mbf{m}$ on the wall $x+y=n_3$, $2x+2y=n_6$, $3x+3y=n_9$, the monodromy operator $\mbf{B}_{\mbf{m}}$ is of the $U_{q}(\hat{\mf{gl}}_1)$-type.

\subsubsection{Matrix Coefficients of $\mbf{B}_{\mbf{m}}$}

We first write down the monodromy operators $\mbf{B}_{\mbf{m}}$ explicitly for $K_{T}(\text{Hilb}_{3}([\mbb{C}^2/\mbb{Z}_2])$:
\begin{align*}
\mbf{B}_{\mbf{m}}(z)=
\begin{cases}
\sum_{n=0}^{\infty}\frac{(q-q^{-1})^n}{[n]_{q^2}!}\frac{(-1)^n}{\prod_{l=1}^{n}(1-z^{n[i,j)}p^{|\mbf{m}|})}(Q^{\mbf{m}}_{-[i,j)})^{n}(P^{\mbf{m}}_{[i,j)})^n& (j-i)\text{mod }2\neq0\\
\prod_{h=1}^2:\exp(-\sum_{k=1}\frac{n_k}{1-z^{-k|\bm{\delta}_h|}p^{k\mbf{m}\cdot\bm{\delta}_h}}\alpha_{-k}^{\mbf{m},h}\alpha_{k}^{\mbf{m},h})&(j-i)\text{mod }2=0
\end{cases}.
\end{align*}

This means that
\begin{align*}
\mbf{B}_{\mbf{m}}=
\begin{cases}
\sum_{n=0}^{\infty}\frac{(q-q^{-1})^n}{[n]_{q^2}!}(Q^{\mbf{m}}_{-[i,j)})^{n}(P^{\mbf{m}}_{[i,j)})^n& (j-i)\text{mod }2\neq0\\
\prod_{h=1}^2:\exp(-\sum_{k=1}n_k\alpha_{-k}^{\mbf{m},h}\alpha_{k}^{\mbf{m},h})&(j-i)\text{mod }2=0
\end{cases}.
\end{align*}

Now we do the classification of the monodromy operators on each wall:
\begin{equation*}
\begin{aligned}
\mbf{B}_{\mbf{m}}(z)=
\begin{cases}
\sum_{n=0}^{1}\frac{(q-q^{-1})^n}{[n]_{q^2}!}\frac{(-1)^n}{\prod_{l=1}^{n}(1-z^{n}_2p^{|\mbf{m}|})}(Q^{\mbf{m}}_{-[2,3)})^{n}(P^{\mbf{m}}_{[2,3)})^n& x=1\\
\sum_{n=0}^{1}\frac{(q-q^{-1})^n}{[n]_{q^2}!}\frac{(-1)^n}{\prod_{l=1}^{n}(1-z^{n}_1p^{|\mbf{m}|})}(Q^{\mbf{m}}_{-[1,2)})^{n}(P^{\mbf{m}}_{[1,2)})^n& y=1\\
\sum_{n=0}^{1}\frac{(q-q^{-1})^n}{[n]_{q^2}!}\frac{(-1)^n}{\prod_{l=1}^{n}(1-z^{2n}_1z_2^np^{|\mbf{m}|})}(Q^{\mbf{m}}_{-[1,4)})^{n}(P^{\mbf{m}}_{[1,4)})^n& 2x+y=1,2,3\\
\sum_{n=0}^{1}\frac{(q-q^{-1})^n}{[n]_{q^2}!}\frac{(-1)^n}{\prod_{l=1}^{n}(1-z^{n}_1z_2^{2n}p^{|\mbf{m}|})}(Q^{\mbf{m}}_{-[2,5)})^{n}(P^{\mbf{m}}_{[2,5)})^n& x+2y=1,2,3\\
\sum_{n=0}^{1}\frac{(q-q^{-1})^n}{[n]_{q^2}!}\frac{(-1)^n}{\prod_{l=1}^{n}(1-z^{3n}_1z^{2n}_2p^{|\mbf{m}|})}(Q^{\mbf{m}}_{-[1,7)})^{n}(P^{\mbf{m}}_{[1,7)})^n& 3x+2y=1,2,3\\
1& 2x+3y=1,2,3\\
\prod_{h=1}^{2}(1-\sum_{k=1}^{3}\frac{n_k}{1-z^{-k}_1z^{-k}_2p^{k}}\alpha_{-k}^{\mbf{m},h}\alpha_{k}^{\mbf{m},h}+\frac{n_1n_2}{(1-z_1^{-1}z_2^{-1}p)(1-z_1^{-2}z_2^{-2}p^2)}\times\\
\times(\alpha^{\mbf{m},h}_{-1}\alpha^{\mbf{m},h}_{-2}\alpha^{\mbf{m},h}_{1}\alpha^{\mbf{m},h}_{2})
-\frac{n_1^3}{6(1-z_1^{-1}z_2^{-1}p)^3}(\alpha^{\mbf{m},h}_{-1})^3(\alpha^{\mbf{m},h}_{1})^3)&x+y=1,2\\
\prod_{h=1}^2(1-\frac{n_1}{1-z_1^{-2}z_2^{-2}p^2}\alpha^{\mbf{m},h}_{-2}\alpha^{\mbf{m},h}_{2})&2x+2y=1,3\\
\prod_{h=1}^2(1-\frac{n_1}{1-z_1^{-3}z_2^{-3}p^3}\alpha^{\mbf{m},h}_{-1}\alpha^{\mbf{m},h}_{1})&3x+3y=1,2,4,5
\end{cases}.
\end{aligned}
\end{equation*}

Using the formula in the fixed point basis as in \eqref{shuffle-formula-1} and \eqref{shuffle-formula-2}, we can write down the explicit matrix coefficients of the monodromy operators. Otherwise we can use the corresponding stable bases $s^{\mbf{m}}_{\bm{\lambda}}$ to compute the matrix coefficients of the monodromy operators. The transition between two stable basis $s^{\mbf{m}_1}_{\bm{\lambda}}$, $s^{\mbf{m}_2}_{\bm{\lambda}}$ of different slopes $\mbf{m}_1$ and $\mbf{m}_2$ is determined by the wall $R$-matrix $R^{\mbf{m}_1,\mbf{m}_2}_{\mc{C}}$ of $K$-theoretic stable envelope, which has been explained in \cite{Z23}.

\subsection{Connection to the Dubrovin connections of Hilbert scheme of $A_{r-1}$-singularities}

By Theorem \ref{degeneration-qde-thm}, the quantum difference equation on $K_{T}(\text{Hilb}_{n}([\mbb{C}^2/\mbb{Z}_r]))$ degenerates to the Dubrovin connection on the quantum cohomology $H_T^*(\text{Hilb}_{n}([\mbb{C}^2/\mbb{Z}_r]))$:
\begin{equation*}
\begin{aligned}
&Q_{\mc{O}(1)}(z)=c_{1}(\mc{O}(1))\cup(-)+(\hbar_1+\hbar_2)\sum_{[i,j)}\frac{|[i,j)|}{1-z^{-[i,j)}}E_{[i,j)}E_{-[i,j)}+\cdots\\
&Q_{\mc{L}_l}(z)=c_1(\mc{L}_l)\cup(-)+(\hbar_1+\hbar_2)\sum_{[i,j)}\frac{|[i,j)|_{l}}{1-z^{-[i,j)}}E_{[i,j)}E_{-[i,j)}+\cdots.
\end{aligned}
\end{equation*}

Here $|[i,j)|$ is the sum of the coefficients of $[i,j)$, $|[i,j)|_{l}$ is the sum of the coefficients of $[i,j)$ in the direction $\mbf{e}_{l}$.

In \cite{MO09}, Maulik and Oblomkov have computed the Dubrovin connection on the quantum cohomology of the Hilbert scheme of $A_{r-1}$-singularities $H_{T}^*(\text{Hilb}_{n}(\widehat{\mbb{C}^2/\mbb{Z}_r}))$, the Dubrovin connection can be written as:
\begin{equation*}
\begin{aligned}
M_{\mc{O}(1)}(z)=&c_1(\mc{O}(1))\cup(-)+(\hbar_1+\hbar_2)\sum_{1\leq i\leq j\leq r}\sum_{k\in\mbb{Z}}:e_{ji}(k)e_{ij}(-k):\frac{k(z_1\cdots z_r)^kz_i\cdots z_{j-1}}{1-(z_1\cdots z_r)^kz_{i}\cdots z_{j-1}}\\
&+\sum_{k\geq1}[rt_1t_2p_{-k}(1)p_{k}(1)+\sum_{i=1}^{r-1}p_{-k}(E_{i})p_{k}(\omega_i)](\frac{k(z_1\cdots z_r)^k}{1-(z_1\cdots z_r)^k}-\frac{z_1\cdots z_r}{1-z_1\cdots z_r}).
\end{aligned}
\end{equation*}
\begin{equation*}
\begin{aligned}
M_{\mc{L}_i}(z)=c_1(\mc{L}_i)\cup(-)-(\hbar_1+\hbar_2)\sum_{1\leq i\leq l\leq j\leq r}\sum_{k\in\mbb{Z}}:e_{ji}(k)e_{ij}(-k):\frac{(z_1\cdots z_r)^kz_i\cdots z_{j-1}}{1-(z_1\cdots z_r)^kz_{i}\cdots z_{j-1}}.
\end{aligned}
\end{equation*}

Here we use the notations in \cite{MO09}.  The equivariant Hilbert scheme $\text{Hilb}_{n}([\mbb{C}^2/\mbb{Z}_r])$ and the Hilbert scheme $\text{Hilb}_{n}(\widehat{\mbb{C}^2/\mbb{Z}_r})$ are the cyclic quiver varieties of the same quiver representation with different stability conditions. By the flop isomorphism $F:H_{T}^*(\text{Hilb}_{n}(\widehat{\mbb{C}^2/\mbb{Z}_r}))\cong H_{T}^*(\text{Hilb}_{n}([\mbb{C}^2/\mbb{Z}_r]))$ defined in \cite{MO12}, together with commutative diagram $(4.38)$ in \cite{MO12} and the Theorem $7.2.1$ in \cite{MO12}, we have that the quantum multiplication by divisors for $\text{Hilb}_{n}([\mbb{C}^2/\mbb{Z}_r])$ and $\text{Hilb}_{n}(\widehat{\mbb{C}^2/\mbb{Z}_r})$ are intertwined by the flop isomorphism up to a scaling operator:
\begin{equation*}
\begin{tikzcd}
H_{T}^*(\text{Hilb}_{n}([\mbb{C}^2/\mbb{Z}_r]))\arrow[r,"F"]\arrow[d,"Q_{\mc{L}_i}(z)"]&H_{T}^*(\text{Hilb}_{n}(\widehat{\mbb{C}^2/\mbb{Z}_r}))\arrow[d,"M_{\mc{L}_i}(z)"]\\
H_{T}^*(\text{Hilb}_{n}([\mbb{C}^2/\mbb{Z}_r]))\arrow[r,"F"]&H_{T}^*(\text{Hilb}_{n}(\widehat{\mbb{C}^2/\mbb{Z}_r}))
\end{tikzcd}.
\end{equation*}

It is conjectured that the quantum difference operator $\mbf{B}_{\mc{L}_i}^s(z)$ in the eigenbasis $H$ of $\mbf{B}_{\mc{L}_i}^s(\infty)$ has the degeneration limit as the quantum multiplication by $M_{\mc{L}_i}$ in the quantum cohomology of $H_{T}^*(\text{Hilb}_{n}(\widehat{\mbb{C}^2/\mbb{Z}_r}))$.

\section{\textbf{Appendix I: Basics of the abelian functions}}\label{sec:_textbf_appendix_i_basics_of_the_abelian_functions}

In this appendix we shall give a brief introduction on the theory of the function field of meromorphic functions over complex tori. The complex tori we consider is $\mc{A}:=\mbb{C}^n/L$ with $L\subset\mbb{C}^n$ a periodic full-rank lattice in $\mbb{C}^n$.

\subsection{Theta Functions}
Given $\Omega$ an $n\times n$ symmetric complex matrices with positive definite imaginary part, we define the Jacobi theta function $\theta(z,\Omega)$ associated to $z\in(\mbb{C}^*)^n$ and $\Omega$ as:
\begin{align*}
\theta(z;\Omega):=\sum_{\mbf{n}\in\mbb{Z}^n}z^{\mbf{n}}\exp(\frac{1}{2}\mbf{n}^t\Omega\mbf{n}).
\end{align*}

It is easy to check that:
\begin{align*}
\vartheta(e^{2\pi i}z;\Omega)=\theta(z;\Omega),\qquad\theta(\exp(\Omega\mc{L})z;\Omega)=\exp(-\frac{1}{2}\mc{L}^t\Omega\mc{L})z^{-\mc{L}}\theta(z;\Omega).
\end{align*}

Moreover if we take $z=e^{2\pi i\bm{\xi}}$, the Jacobi theta function can be written as:
\begin{align*}
\vartheta(z;\Omega):=\sum_{\mbf{n}\in\mbb{Z}^n}\exp(2\pi i\mbf{n}^t\bm{\xi})\exp(\frac{1}{2}\mbf{n}^t\Omega\mbf{n}).
\end{align*}

Also given $a,b\in\mbb{R}^g$, we define the theta function $\vartheta\begin{pmatrix}a\\b\end{pmatrix}(z;\Omega)$ as:
\begin{align*}
\vartheta\begin{pmatrix}a\\b\end{pmatrix}(z;\Omega)=\exp(\frac{1}{2}a^t\Omega a+a^t(\xi+b))\vartheta(\xi+\Omega\cdot a+b;\Omega)
\end{align*}
such that:
\begin{align*}
\vartheta\begin{pmatrix}a\\b\end{pmatrix}(z+\Omega\cdot m+n;\Omega)=\exp(2\pi i (a^tn-b^tm))\exp(-\frac{1}{2}m^t\Omega m-2\pi im^tz)\vartheta\begin{pmatrix}a\\b\end{pmatrix}(z;\Omega),\qquad\forall m,n\in\mbb{Z}^g.
\end{align*}

As the special example, if we take $n=1$ and $\Omega=\tau$, and set  $q=e^{2\pi i\tau}$, the Jacobi theta function can be written as:
\begin{align*}
\vartheta(z;q):=\sum_{n\in\mbb{Z}}q^{\frac{n^2}{2}}z^n=\prod_{n=1}^{\infty}(1-q^{2n})(1+q^{2n-1}z)(1+q^{2n-1}z^{-1}).
\end{align*}

It satisfies the following modular property:
\begin{align}\label{modular-property}
\vartheta(e^{2\pi i\xi};e^{2\pi i\tau})=\frac{i}{\sqrt{\tau}}e^{-\frac{i\pi\xi^2}{\tau}}\vartheta(e^{\frac{2\pi i\xi}{\tau}};e^{-\frac{2\pi i}{\tau}}).
\end{align}

In the context of \cite{AO21}, they use the odd Jacobi theta function:
\begin{align*}
\theta(z|q)=(z^{1/2}-z^{-1/2})\prod_{n\geq1}(1-q^nz)(1-q^nz^{-1}).
\end{align*}

The relation between the Jacobi theta function and the odd Jacobi theta function can be written as:
\begin{align*}
\vartheta(-zq^{1/2};q)=-z^{\frac{1}{2}}\varphi(q)\theta(z),\qquad\varphi(q)=\prod_{n=1}^{\infty}(1-q^{2n}).
\end{align*}

In our settings, We will use the Jacobi theta function $\vartheta(z;\Omega)$ and the odd Jacobi theta function $\theta(z|q)$ depending on the situation. It does not matter too much which kind of theta function we would use since they are only different up to a scaling transformation.

\subsection{Line bundles over the abelian varieties}

Fix the complex space $\mbb{C}^n=V$ and a rank $2n$ lattice $\Lambda\subset V$, and we define $X=V/\Lambda$. The line bundle $L$ on $X$ is uniquely determined by the pair $(H,\chi)$. the Hermitian form $H:V\times V\rightarrow\mbb{C}$ such that $\text{Im }H(\Lambda,\Lambda)\subset\mbb{Z}$. $\chi:\Lambda\rightarrow S^1\subset\mbb{C}$ is the semicharacter for $H$, which is defined as:
\begin{align*}
\chi(\lambda+\mu)=\chi(\lambda)\chi(\mu)\exp(\pi i\text{Im }H(\lambda,\mu)),\qquad\forall\lambda,\mu\in\Lambda.
\end{align*}

The construction is given by $(H,\chi)$, which defines the automorphy form $a=a_{(H,\chi)}:\Lambda\times V\rightarrow\mbb{C}^*$ by
\begin{align*}
a(\lambda,v):=\chi(\lambda)\exp(\pi H(v,\lambda)+\frac{\pi}{2}H(\lambda,\lambda))
\end{align*}
which satisfies the cocycle relation:
\begin{align*}
a(\lambda+\mu,v)=a(\lambda,v+\mu)a(\mu,v)
\end{align*}
In this way, we can define the line bundle $L(H,\chi):=V\times\mbb{C}/\Lambda$ such that $\Lambda$ acts on $V\times\mbb{C}$ by $\lambda\cdot(v,t)=(v+\lambda,a(\lambda,v)t)$.

Now suppose that given a line bundle $L$ with the automorphy form $f(\lambda,v)$, the elements in the global holomorphic section $H^0(L)$ of $L$ can be thought of as the set of holomorphic functions $\vartheta:V\rightarrow\mbb{C}$ such that
\begin{align}\label{theta-difference}
\vartheta(v+\lambda)=f(\lambda,v)\vartheta(v).
\end{align}

It is known that the global section of the line bundle can be written as $\vartheta(z^{\mbf{n}}e^{2\pi\Omega\cdot\mbf{m}};\Omega)$. We can choose vector $\mbf{n}\in\mbb{Z}^n$ and $\mbf{m}\in\mbb{R}^n$ such that $\vartheta(z^{\mbf{n}}e^{2\pi\Omega\cdot\mbf{m}};\Omega)$ satisfies the difference equation \eqref{theta-difference}.

Moreover, if we consider the space of meromorphic sections $\Gamma_{rat}(L)$, it is generated by the products of the theta functions of the form:
\begin{align*}
\prod_{i}\vartheta(z^{\mbf{n}_i}e^{2\pi\Omega\cdot\mbf{m}_i};\Omega)^{\lambda_i},\qquad\lambda_i\in\mbb{Z}.
\end{align*}

For the specific example that we use in this paper, we consider the space of meromorphic sections of the structure sheaf $\Gamma_{rat}(\mc{O}_X)$, i.e. the space of abelian functions over $X$. This means that the section $f(z)\in\Gamma_{rat}(\mc{O}_X)$ is meromorphic and $\Lambda$-periodic. This implies that $\Gamma_{rat}(\mc{O}_{X})$ is generated by the theta functions of the form $\prod_{i}\vartheta(z^{\mbf{n}_i}e^{2\pi\Omega\cdot\mbf{m}_i};\Omega)^{\lambda_i},\lambda_i\in\mbb{Z}$ such that $\sum_{i}\lambda_i\mbf{n}_i=\sum_{i}\Omega\cdot(\lambda_i\mbf{m}_i)=0$.

More precisely, if we choose the lattice of the form $\Lambda_{\Omega,D}=\Omega\cdot\mbb{Z}^g\oplus D\cdot\mbb{Z}^g$ with $D=\text{diag}(d_1,\cdots,d_g)$. Given $a,b\in\mbb{R}^g$, let $\{a_k|1\leq k\leq r^g\}$ be a system of representatives of $(a+D^{-1}\mbb{Z}^g)/\mbb{Z}^g$, and let $\{b_l|1\leq l\leq\prod_{i=1}^gd_i\}$ be a system of representatives of $(r^{-1}b+r^{-1}\mbb{Z}^g)/\mbb{Z}^g$. Then the family of theta functions $\{\vartheta\begin{pmatrix}a_k\\b_l\end{pmatrix}(z;r^{-1}\Omega)\}$ defines a projective embedding $\mbb{C}^g/\Lambda_{\Omega,D}\hookrightarrow\mbb{P}^{(\prod_{i=1}^gd_i)\cdot r^g-1}$. 

Using the projective embedding, we can show that the abelian function on $\mbb{C}^g/\Lambda_{\Omega,D}$ is generated by the following family of quotients of theta functions:
\begin{align*}
\frac{\vartheta\begin{pmatrix}a_{k_1}\\b_{l_1}\end{pmatrix}(z;r^{-1}\Omega)}{\vartheta\begin{pmatrix}a_{k_2}\\b_{l_2}\end{pmatrix}(z;r^{-1}\Omega)},\qquad1\leq k_1,k_2\leq r^g,1\leq l_1,l_2\leq\prod_{i=1}^gd_i.
\end{align*}

\subsection{Nodal Limit of the abelian functions}
In our setting we only consider the special case that $\Omega=\tau\text{Id}$, $D=\text{Id}$ with $\tau\in\mbb{H}$. In this case:
\begin{align*}
\vartheta\begin{pmatrix}a\\b\end{pmatrix}(z;\Omega)=\exp(\sum_{i=1}^g\pi i\tau a_i^2+2\pi ia_i(z_i+b_i))\prod_{i=1}^g\vartheta(z_i+a_i\tau+b_i|\tau).
\end{align*}

In this case, we can choose $b=0$, and the abelian function on $\mbb{C}^g/\Lambda_{\Omega,D}$ is generated by:
\begin{align*}
\exp(\sum_{i=1}^g\pi ir^{-1}\tau(a_i^2-b_i^2)+2\pi i(a_i-b_i)z_i)\prod_{i=1}^g\frac{\vartheta(z_i+r^{-1}a_i\tau|r^{-1}\tau)}{\vartheta(z_i+r^{-1}b_i\tau|r^{-1}\tau)}.
\end{align*}

Using the modular duality for the theta function:
\begin{align*}
\vartheta(e^{2\pi i\xi},e^{2\pi i\tau})=\frac{i}{\sqrt{\tau}}e^{-\frac{i\pi\xi^2}{\tau}}\vartheta(e^{\frac{2\pi i\xi}{\tau}},e^{-\frac{2\pi i}{\tau}}).
\end{align*}

Now we compute the limit $\tau\rightarrow0$ of the above generators:
\begin{equation*}
\begin{aligned}
&\lim_{\tau\rightarrow0}\exp(\sum_{i=1}^g\pi ir^{-1}\tau(a_i^2-b_i^2)+2\pi i(a_i-b_i)z_i)\prod_{i=1}^g\frac{\vartheta(z_i+r^{-1}a_i\tau|r^{-1}\tau)}{\vartheta(z_i+r^{-1}b_i\tau|r^{-1}\tau)}\\
=&\lim_{\tau\rightarrow0}\prod_{i=1}^g\frac{\vartheta(\frac{rz_i}{\tau}+a_i|-\frac{r}{\tau})}{\vartheta(\frac{rz_i}{\tau}+b_i|-\frac{r}{\tau})}=\lim_{q\rightarrow0}\prod_{i=1}^g\frac{\theta(q^{-z_i}e^{2\pi ia_i})}{\theta(q^{-z_i}e^{2\pi ib_i})},\qquad q=e^{-\frac{2\pi ir}{\tau}}
\end{aligned}
\end{equation*}
i.e. we have that in terms of the odd Jacobi theta function:
\begin{equation}\label{modular-transform}
\lim_{\tau\rightarrow0}\exp(\sum_{i=1}^g2\pi i(a_i-b_i)z_i)\prod_{i=1}^g\frac{\theta(e^{2\pi i(z_i+r^{-1}a_i\tau)|e^{2\pi ir^{-1}\tau}})}{\theta(e^{2\pi i(z_i+r^{-1}b_i\tau)|e^{2\pi ir^{-1}\tau}})}=\lim_{q\rightarrow0}\prod_{i=1}^{g}\frac{\theta(q^{-z_i}e^{2\pi ia_i})}{\theta(q^{-z_i}e^{2\pi ib_i})},\qquad q=e^{-\frac{2\pi ir}{\tau}}.
\end{equation}

This formula also appears in Lemma $9$ of \cite{S21}, which plays an important role in the computation of the degeneration limit of the connection matrix and the elliptic stable envelope of the Hilbert scheme of points $\text{Hilb}_{n}(\mbb{C}^2)$.

\section{\textbf{Appendix II: Distinct eigenvalues of $\mbf{B}_{\mc{O}(1)}^s(z)$}}\label{appendix-2}\label{sec:_textbf_appendix_ii_distinct_eigenvalues_of_mbf_b}

In this appendix we prove that $\mbf{B}_{\mc{O}^s(1)}(\infty)$ has distinct eigenvalues on the weighted modules $K_{T}(M(\mbf{v},\mbf{w}))$ for generic $q$, $t$ and $u_1,\cdots,u_{\mbf{w}}$.

In fact, we prove a stronger result:
\begin{prop}
For generic $q,t,u_1,\cdots,u_{\mbf{w}}$ and $z_1,\cdots,z_n$ including $z=0$ and $z=\infty$, $\mbf{B}_{\mc{O}(1)}^s(z)$ has distinct eigenvalues.
\end{prop}
Recall that $\mbf{B}_{\mc{O}(1)}^s(z)$ can be written as:
\begin{align*}
\mbf{B}_{\mc{O}(1)}^s(z)=\mc{O}(1)\prod^{\leftarrow}_{\mbf{m}\in[s,s-\mc{L})}\mbf{B}_{\mbf{m}}(z).
\end{align*}
with the path $[s,s-\mc{L})$ being generic.
It is well-known that $\mc{O}(1)$ has distinct eigenvalues with the eigenvectors given by the fixed point basis.

Our strategy is the following: We denote its degeneration limit as $\mbf{M}^{coh}_{\mc{O}(1)}(z)$. By computation we have that:
\begin{align*}
\mbf{M}^{coh}_{\mc{O}(1)}(z)=c_{1}(\mc{O}(1))\cup(-)+(\hbar_1+\hbar_2)\sum_{j>i}\frac{|[i,j)|}{1-z^{-[i,j)}}E_{[i,j)}E_{-[i,j)}.
\end{align*}

If we can prove that for the generic value of $\mbb{C}_{\hbar_1}\times\mbb{C}_{\hbar_2}\times t_{\mbf{w}}\times(\text{Pic}(X)\otimes\mbb{C}^*)$, $\mbf{M}^{coh}_{\mc{O}(1)}(z)$ has distinct eigenvalues, this implies that the corresponding quantum difference operator $\mbf{B}_{\mc{O}(1)}^s(z)$ has distinct eigenvalues.

We mimic the method in section $15.2$ in \cite{MO12}. Choose a path $\gamma:[0,1)\rightarrow\mbb{C}_{\hbar_1}\times\mbb{C}_{\hbar_2}\times t_{\mbf{w}}\times(\text{Pic}(X)\otimes\mbb{C}^*)$ such that $\gamma(0)=(\hbar_1,\hbar_2,a_1,\cdots,a_{\mbf{w}},0)$ and $\gamma(1)=(\hbar_1',\hbar_2',a_1',\cdots,a_{\mbf{w}}',z_1,\cdots,z_n)$ such that $(z_1,\cdots,z_n)$ lies on one wall $z^{[i,j)}=1$ intersecting on the path $[s,s-\mc{L})$. As $t\rightarrow1$, the factor $\mbf{B}_{\mbf{m}}(z)$ corresponding to the wall $w$ dominates the product and the rest of the terms in the product remain bounded and nonzero for generic $z\in w$.

Since $\mbf{M}_{\mc{O}(1)}^{coh}(\gamma(0))=c_1(\mc{O}(1))$ has distinct eigenvalues, it is diagonalisable. For $\mbf{M}_{\mc{O}(1)}^{coh}(z)$ with $z\in U$ in an open small neighborhood of $0$, $\mbf{M}_{\mc{O}(1)}^{coh}(z)$ has distinct eigenvalues and thus is diagonalisable. Also, each eigenvalue is a Laurent formal power series of $z$.

For $\mbf{M}_{\mc{O}(1)}^{coh}(\gamma(t))$ as $t\rightarrow1$, the term $\frac{(\hbar_1+\hbar_2)|[i,j)|}{1-z^{-[i,j)}}E_{[i,j)}E_{-[i,j)}$ dominates the term. While since $E_{[i,j)}E_{-[i,j)}$ is not an automorphism for generic $H_{T}^*(M(\mbf{v},\mbf{w}))$, the operator $E_{[i,j)}E_{-[i,j)}$ does not have full rank. We also set the Kahler variable $(z_1,\cdots,z_n)$ to be really large in some components on the wall such that the other terms of the form $\frac{1}{1-z^{-\alpha}}$ remain sufficiently small. In this case some of the eigenvalues of $\mbf{M}_{\mc{O}(1)}^{coh}(\gamma(t))$ as $t\rightarrow1$ are of the form $\frac{1}{|z-\gamma(1)|}$ with the eigenvalues of $c_1(\mc{O}(1))$. For the nonzero eigenvalues of $\frac{E_{[i,j)}E_{-[i,j)}}{1-z^{-[i,j)}}$, they contribute to the eigenvalues of the form $\frac{a}{|z-z_0|}$, the contribution from $c_1(\mc{O}(1))$ is independent of $z$, and this implies that we can choose generic $z$ such that the eigenvalues are still distinct. Since the eigenvalues are determined by the algebraic equations in $\mbb{C}_{\hbar_1}\times\mbb{C}_{\hbar_2}\times t_{\mbf{w}}\times(\text{Pic}(X)\otimes\mbb{C}^*)$, we can see that there exists Zariski open subset $U$ of $\mbb{C}_{\hbar_1}\times\mbb{C}_{\hbar_2}\times t_{\mbf{w}}\times(\text{Pic}(X)\otimes\mbb{C}^*)$ such that $\mbf{M}_{\mc{O}(1)}^{coh}(z)$ has distinct eigenvalues.

Recall that in \cite{MO12} the quantum multiplication of divisors of quiver varieties of different stability conditions are connected by the flop isomorphism:
\begin{equation*}
\begin{tikzcd}
H_{T}(M_{\theta,0}(\mbf{v},\mbf{w}))\arrow[r,"Q_{\theta,\mc{L}}(z)"]\arrow[d,"F_{\theta,\theta'}"]&H_{T}(M_{\theta,0}(\mbf{v},\mbf{w}))\arrow[d,"F_{\theta,\theta'}"]\\
H_{T}(M_{\theta',0}(\mbf{v},\mbf{w}))\arrow[r,"Q_{\theta',\mc{L}}(z)"]&H_{T}(M_{\theta',0}(\mbf{v},\mbf{w}))
\end{tikzcd}
\end{equation*}
up to a scalar operator. Thus we have:
\begin{align*}
Q_{\theta,\mc{L}}(z)=c_1(\mc{L}_{\theta})+(\hbar_1+\hbar_2)\sum_{\alpha\cdot\theta>0}\frac{z^{\alpha}\langle\alpha,\mc{L}\rangle}{1-z^{\alpha}}e_{\alpha}e_{-\alpha}.
\end{align*}

Now we take $\theta=(1,\cdots,1)$, $\theta'=-\theta$. In this case we have that:
\begin{equation*}
\begin{aligned}
&c_1(\mc{O}(1)_{-\theta})+(\hbar_1+\hbar_2)\sum_{j>i}\frac{|[i,j)|}{1-z^{[i,j)}}E_{-[i,j)}E_{[i,j)}\\
=&\text{Ad}_{F}(c_1(\mc{O}(1)_{\theta})-(\hbar_1+\hbar_2)\sum_{j>i}\frac{|[i,j)|}{1-z^{-[i,j)}}E_{[i,j)}E_{-[i,j)})+\text{Scalar operator}.
\end{aligned}
\end{equation*}

Since $c_1(\mc{O}(1)_{-\theta})$ also has distinct eigenvalues. (This follows from the fact that the affine type $A$ quiver varieties of stability condition $-\theta$ can be thought of as the subvariety of the Jordan quiver variety of the opposite stability condition). Using the similar procedure above, we can see that there is a Zariski open subset $U$ such that $Q_{-\theta,\mc{O}(1)}(z)$ has distinct eigenvalues. Also $\mbf{M}^{coh}_{\mc{O}(1)}(\infty)$ has distinct eigenvalues.

\end{document}